\documentclass[11pt,a4paper,reqno]{amsart}
\usepackage{amssymb}
\usepackage{amsmath}
\usepackage{amsthm}
\usepackage[dvipsnames]{xcolor}
\usepackage[shortlabels]{enumitem}
\usepackage{braket,comment,soul,cancel}
\usepackage{hyperref}
\usepackage{pdfpages,faktor}
\usepackage{graphicx,mathabx,bbm}
\usepackage{caption}
\usepackage{subcaption}
\usepackage{mathrsfs} 
\usepackage{tikz-cd} 
\usepackage{bm}
\usepackage{enumitem}
\usepackage{mathtools}
\usepackage{pgfplots}
\usepackage{import}
\usepackage{xifthen}
\usepackage{transparent}
\usepackage{xargs}

\usepackage[colorinlistoftodos,prependcaption,textsize=tiny]{todonotes}
\newcommandx{\unsure}[2][1=]{\todo[linecolor=red,backgroundcolor=red!25,bordercolor=red,#1]{#2}}
\newcommandx{\change}[2][1=]{\todo[linecolor=blue,backgroundcolor=blue!25,bordercolor=blue,#1]{#2}}
\newcommandx{\info}[2][1=]{\todo[linecolor=OliveGreen,backgroundcolor=OliveGreen!25,bordercolor=OliveGreen,#1]{#2}}
\newcommandx{\improvement}[2][1=]{\todo[linecolor=Plum,backgroundcolor=Plum!25,bordercolor=Plum,#1]{#2}}

\numberwithin{equation}{section}
\theoremstyle{plain}
\newtheorem{theorem}{Theorem}[section]

\newtheorem{prop}[theorem]{Proposition}
\newtheorem{lem}[theorem]{Lemma}
\newtheorem{cor}[theorem]{Corollary}

\newtheorem*{question*}{Question}

\newtheorem{conj}[theorem]{Conjecture}

\theoremstyle{definition}
\newtheorem{defn}[theorem]{Definition}
\newtheorem{remark}[theorem]{Remark}

\theoremstyle{definition}
\newtheorem{thmx}{Theorem}

\newcommand{\Hyp}{\mathbb{H}}
\newcommand{\R}{\mathbb{R}}

\newcommand{\C}{\mathbb{C}}

\newcommand{\Z}{\mathbb{Z}}

\newcommand{\N}{\mathbb{N}}
\newcommand{\D}{\mathbb{D}}
\newcommand{\E}{\mathscr{E}}

\newcommand{\cB}{\mathcal{B}}
\newcommand{\cE}{\mathcal{E}}
\newcommand{\cF}{\mathcal{F}}
\newcommand{\cG}{\mathcal{G}}
\newcommand{\cH}{\mathcal{H}}
\newcommand{\cN}{\mathcal{N}}
\newcommand{\cP}{\mathcal{P}}

\newcommand{\cK}{\mathcal{K}}
\newcommand{\cL}{\mathcal{L}}
\newcommand{\cD}{\mathcal{D}}
\newcommand{\cT}{\mathcal{T}}

\newcommand{\RT}{\mathscr{T}}

\renewcommand{\SS}{\mathbb S}
\newcommand{\bS}{\mathbb S}

\renewcommand{\epsilon}{\varepsilon}

\DeclareMathOperator{\PSL}{PSL}

\DeclareMathOperator{\Int}{Int}

\newcommand{\mate}{\bot \!\! \! {\bot}} 

\pgfplotsset{compat=1.18} 

\numberwithin{figure}{section}

\makeatletter
\DeclareFontFamily{U}{tipa}{}
\DeclareFontShape{U}{tipa}{m}{n}{<->tipa10}{}
\newcommand{\arc@char}{{\usefont{U}{tipa}{m}{n}\symbol{62}}}

\newcommand{\arc}[1]{\mathpalette\arc@arc{#1}}

\newcommand{\arc@arc}[2]{
  \sbox0{$\m@th#1#2$}
  \vbox{
    \hbox{\resizebox{\wd0}{\height}{\arc@char}}
    \nointerlineskip
    \box0
  }
}
\makeatother

\title[Mateability of Fuchsian lattices with polynomials]{Mateability of Fuchsian lattices with polynomials, and failure of algebraicity}

\begin{author}[Y.~Luo]{Yusheng Luo}
\address{Department of Mathematics, Cornell University, 212 Garden Ave, Ithaca, NY 14853, USA}
\email{yusheng.s.luo@gmail.com}
\thanks{Y.L. was partially supported by NSF Grant DMS-2349929.}
\end{author}

\begin{author}[M.~Mj]{Mahan Mj}
\address{School of Mathematics, Tata Institute of Fundamental Research, 1 Homi Bhabha Road, Mumbai 400005, India}
\email{mahan@math.tifr.res.in, mahan.mj@gmail.com}
\thanks{M.M. was partially supported by  the Department of Atomic Energy, Government of India, under Project Identification No. RTI 4014, a JC Bose Grant from ANRF, and an endowment of the Infosys Foundation.}
\end{author}

\begin{author}[S.~Mukherjee]{Sabyasachi Mukherjee}
\address{School of Mathematics, Tata Institute of Fundamental Research, 1 Homi Bhabha Road, Mumbai 400005, India}
\email{sabya@math.tifr.res.in, mukherjee.sabya86@gmail.com}
\thanks{S.M. was partially supported by the Department of Atomic Energy, Government of India, under Project Identification No. RTI 4014, an endowment of the Infosys Foundation, and ANRF research project grant ANRF/ARGM/2025/000089/MTR.}
\end{author}

\begin{document}

\begin{abstract}
We prove that all Fuchsian lattices, with the possible exceptions of those associated with compact orbifolds of type $(0,3)$ and $(1,1)$, can be mated with complex polynomials, providing the first examples of such combinations beyond genus-zero orbifolds. We also present the first construction of matings between simply degenerate Kleinian groups and polynomial maps. These matings are realized as holomorphic, though not necessarily algebraic, correspondences,  demonstrating that Fatou’s vision of unifying Kleinian group actions and polynomial dynamics can be fulfilled in the wider setting of transcendental  correspondences.
In fact, we show that, under a natural notion of mating, a Fuchsian lattice that gives rise to an algebraic correspondence is necessarily of genus zero.

We achieve this by generalizing the orbit equivalence mating framework introduced in our earlier work. Specifically, we use pinching and acceleration techniques to construct a new class of piecewise M{\"o}bius covering maps of the circle that are orbit equivalent to Fuchsian groups. By mating these circle coverings with polynomials, we construct degenerate polynomial-like maps whose analytic continuations yield the desired correspondences. The nature of this orbit equivalence plays a crucial role in the mating description, as well as in governing whether the resulting correspondences are algebraic.
\end{abstract}

\date{\today}

\maketitle

\setcounter{tocdepth}{1}
\tableofcontents

\section{Introduction}\label{sec-intro}
This paper is a contribution to an overarching program that tries to unify two kinds of  conformal dynamical systems: one coming from iteration of rational maps, the other from actions of Kleinian groups (see for instance \cite{BP94,BL1,LLMM23,MM1,MM2,LLM24,LMM26a,LM25a,LM25b}). The problem goes back to an  empirical observation by
Fatou \cite{Fatou29} who observed a striking similarity between these two classes of holomorphic dynamical systems. Fatou \cite{Fatou29} asked for a unification of these two theories in terms of a framework of iterated algebraic correspondences; i.e., multi-valued maps with holomorphic local branches globally generated by a polynomial equation $Q(z,w)=0$. 

\subsection*{Combining Fuchsian groups with polynomials: the genus-zero regime} In a pioneering work, Bullett-Penrose gave the first examples of algebraic correspondences realizing holomorphic matings between certain quadratic polynomials and the modular group \cite{BP94}. The resulting family of correspondences was studied extensively by Bullett and Lomonaco, who proved that it contains matings between the modular group and all quadratic parabolic rational maps \cite{BL1,BL2,BL3}. 

A general framework for mating antiholomorphic polynomials with reflection groups was developed in \cite{LLMM21,LLMM23,LMM22,LMM24,LLM24,LMMN}, through a dynamical theory of Schwarz reflections. Inspired by the anti-holomorphic framework, an orbit equivalence framework for conformal matings was initiated in \cite{MM1} and developed through \cite{MM2, LLM24, LMM26a,LMM26c}. This is done using the following broad scheme:
\begin{enumerate}[leftmargin=8mm]
 \item[(I)] Replace a Fuchsian group $G$ acting on the circle $\mathbb{S}^1$ by a single conformal Markov map $A$, typically a variant of the Bowen--Series map \cite{BS79}, with the following properties.
    \begin{enumerate}[leftmargin=8mm]
        \item $A$ is orbit equivalent to $G$; i.e., the grand orbits of $A$ equal the orbits of $G$. The map $A$ extends holomorphically to an open subset of the unit disk~$\D$.
        \item The action of the map $A$ on the circle is topologically conjugate to the conformal model $z^d$ of a complex polynomial $P$ on its Julia set,  via a topological conjugacy $\psi$.
    \end{enumerate}
    
    \item[(II)]  A \emph{topological mating}
    $P \mate A$ is then constructed as a map from an open subset of the sphere $\bS^2$ to $\bS^2$. Thus, $P \mate A$ combines dynamically the action of the conformal Markov $A$ on $\D$ and the polynomial $P$ on its filled Julia set by pasting them together along the Julia set respecting their topological dynamics.
    
    \item[(III)] The {topological mating}
    $P \mate A$ is then upgraded to a holomorphic map $F$ that is topologically conjugate to it. At this stage, we are interested in a holomorphic combination in the spirit of Bers' simultaneous uniformization  and polynomial mating. Equivalently, the map $F$ is a degenerate \emph{polynomial-like map} whose dynamics on the \emph{non-escaping set} is conformally conjugate to the polynomial $P$, and whose \emph{escaping set} dynamics is conformally conjugate to the Markov map~$A$.
    
    \item[(IV)] Finally, it turns out that these conformal matings lead to algebraic correspondences and recapture the Bullett--Penrose examples of \cite{BP94} through a very different complex analytic route.
\end{enumerate}

 The techniques  of \cite{MM1,MM2,LLM24} imposed 
 a severe topological restriction on the collection of Fuchsian groups that could be mated with polynomials:
 \begin{enumerate}[leftmargin=8mm]
     \item $\Sigma$, the underlying surface, had to have genus zero,
     \item $\Sigma$ could have at most one orbifold point of order greater than 2,
     \item $\Sigma$ was necessarily a punctured surface.
 \end{enumerate}

\subsection*{Breaking the genus-zero barrier}

A first aim of this paper is to broaden the above combination framework considerably to finite-area hyperbolic orbifolds $\Sigma = \Sigma_{g,n;\nu_1,\cdot, \nu_n}$ of genus $g$, $n$ marked points with orbifold signature $\nu_i \in \N_{\geq 2} \cup \{\infty\}$,
  so long as $\Sigma$ either has a \emph{puncture}, or has an \emph{essential separating simple closed curve}. Thus, the only excluded cases are compact surfaces $\Sigma$ with $g=0, n=3$ (sphere with three orbifold points), and $g=1, n=1$ (torus with one orbifold point). If $G$ is a Fuchsian lattice with $\Sigma\cong \D/G$, then the \emph{signature} of $G$ is defined to be the signature of $\Sigma$. 

\begin{thmx}\label{most_fuchsian_mateable_intro_thm}
    Let $G$ be a Fuchsian lattice. 
    \begin{itemize}[leftmargin=6mm]
        \item If $G$ is non-uniform, then there exists $d_0\in\N$ such that for any hyperbolic polynomial $P$ 
         of degree~$d_0$ with connected Julia set, there exists a mating between $G$ and $P$.
        \item If $G$ is uniform with signature $(g,n)\notin \{(0,3),(1,1)\}$, then there exists $d_0\in\N$ such that for any polynomial $P$ of degree~$d_0$ with connected Julia set, there exists a mating between $G$ and $P$.
    \end{itemize}
\end{thmx}
 Here, as in the genus-zero situation, a mating between $G$ and $P$ is a degenerate polynomial-like map whose non-escaping set dynamics is conformally conjugate to the polynomial $P$ (on its filled Julia set), and whose escaping set dynamics is conformally conjugate to the action of a piecewise M{\"o}bius Markov covering map that encodes the dynamics of $G$.
For the precise formulations, we refer the reader to \S \ref{conf_mating_def_subsec}, Theorems~\ref{uniform_latice_infused_mating_thm} and~\ref{uniform_lattice_holo_dyn_syst_thm} for the case of uniform lattices and to Theorems~
~\ref{non_uniform_latice_mating_thm} and~\ref{non_uniform_lattice_holo_dyn_syst_thm} for the non-uniform lattice case.

\subsection*{Accelerated Bowen–Series maps: weak encoding of Fuchsian lattices by conformal Markov coverings} There are serious challenges in implementing the scheme of the genus-zero regime to prove Theorem~\ref{most_fuchsian_mateable_intro_thm}. In fact, the strategy to replace a Fuchsian lattice with a conformal Markov map which is orbit equivalent to the group and topologically conjugate to a power map, is itself fraught with difficulty in the higher-genus situation. While the role of such maps was played by Bowen--Series maps (or their factors) in the genus-zero case, the discontinuities of Bowen--Series maps of closed surfaces are more obstinate. Nonetheless, we construct continuous avatars of Bowen--Series maps, paving the way for mating most Fuchsian lattices with polynomials.

We say that a piecewise Fuchsian Markov self-covering of the circle is \emph{weakly} (respectively, \emph{strongly}) orbit equivalent to $G$ if the orbit equivalence property holds only on a full-measure subset of the circle (respectively, everywhere on the circle).

\begin{thmx}\label{acc_bs_intro_thm}
Let $G$ be a Fuchsian lattice such that $\D/G$ is neither a sphere with three orbifold points, nor a torus with one orbifold point. Then, there exists a conformal Markov covering map $\cH_G$ of the circle, termed as the \emph{accelerated Bowen--Series map}, which is weakly orbit equivalent to $G$ and is topologically conjugate to $z^d\vert_{\bS^1}$, where $d=\deg\left(\cH_G:\bS^1\to\bS^1\right)$. If $\psi:\bS^1\to\bS^1$ is a homeomorphism conjugating $z^d$ to $\cH_G$, then
\begin{enumerate}[leftmargin=8mm]
    \item[(i)] for $G$ uniform, $\psi$ is quasisymmetric, and it extends continuously to a quasiconformal homeomorphism of $\D$ conjugating $z^d$ to $\cH_G$ on a (relative) pinched neighborhood of $\bS^1$, and
    \item[(ii)] for $G$ non-uniform, $\psi$ extends continuously to a David homeomorphism of $\D$.
\end{enumerate}
Further, if $G$ is a torsion-free non-uniform lattice, then $\cH_G$ is strongly orbit equivalent to $G$. If $G$ is uniform, then $\cH_G$ is not strongly orbit equivalent to $G$.
\end{thmx}

Theorem~\ref{acc_bs_intro_thm} yields the first examples, beyond the punctured sphere case, of piecewise Fuchsian, topologically expanding Markov covering maps of the circle that are strongly orbit equivalent to Fuchsian lattices (answering \cite[Question~5.3]{MM1} in the negative). An immediate consequence of Theorem~\ref{acc_bs_intro_thm} is the following generalization of \cite[Theorem~8.7]{MM1}.
\begin{cor}[Topological orbit equivalence rigidity fails dramatically]\label{all_nonuniform_lattice_top_oe_cor}
Let $G_1,G_2$ be torsion-free non-uniform lattices. Then the actions of $G_1,G_2$ on $\bS^1$ are topologically orbit equivalent; i.e., there exists a homeomorphism of the circle that carries $G_1$-orbits onto $G_2$-orbits. However, the two actions are not topologically conjugate when $\D/G_1\ncong\D/G_2$.
\end{cor}

We refer the reader to  Theorems~\ref{thm:accnon-uniformlattice},~\ref{thm:accuniformlattice},~\ref{strong_oe_thm},~\ref{thm:weakorbitequiv}, and Lemma~\ref{qc_david_conj_circle_lem} for more detailed statements. Corollary~\ref{all_nonuniform_lattice_top_oe_cor} is proved right after Theorem~\ref{strong_oe_thm}.

\subsection*{Promotion to holomorphic dynamical systems, and failure of algebraicity}
 We now turn to the last step of the scheme employed in the genus-zero regime.
Let $G$ be a Fuchsian lattice and $P$ a  polynomial. The output 
 of Theorem~\ref{most_fuchsian_mateable_intro_thm} is a degenerate polynomial-like map $F$ which is a single-valued holomorphic map. 
The map $F$ captures the full dynamics of $P$ on its filled Julia set,
but does not capture the full dynamics of the group $G$ on the unit disk
$\mathbb{D}$ (indeed, it only contains a conformal copy of the accelerated Bowen--Series map $\cH_G$). In \S \ref{holo_dyn_syst_from_conf_mating_subsec}, we describe a process to promote the degenerate polynomial-like maps arising from Theorem~\ref{most_fuchsian_mateable_intro_thm} to holomorphic dynamical systems $\mathscr{F}$ in the sense of McMullen--Sullivan \cite{McM96,MS98} that completely encode the action of $G$ on the disk. 
Roughly speaking, $\mathscr{F}$ consists of maximal analytic continuations of the branches of
$F^{-m}\circ F^{n}$, $m,n\geq 0$.

Somewhat fortuitously, in the genus-zero case this holomorphic dynamical system turns out to be generated by the iterates of an algebraic correspondence. However, in general, these dynamical systems are generated by holomorphic, but not necessarily algebraic, correspondences. This represents a fundamental departure of our broader combination theory from the genus-zero case.
(A holomorphic correspondence in the context of this paper will in general refer to multi-valued holomorphic maps on open subsets of the Riemann sphere.)

\begin{thmx}\label{holo_dyn_syst_intro_thm}
Let $F$ be a degenerate polynomial-like map realizing the mating between an accelerated Bowen--Series map $\cH_G$ (associated with a Fuchsian lattice $G$) and a polynomial $P$, produced by Theorem~\ref{most_fuchsian_mateable_intro_thm}. Then the maximal analytic continuations of $F^{-m}\circ F^n$, $m,n\geq 0$, yield a holomorphic dynamical system $\mathscr{F}$ with the following properties.
\begin{enumerate}[leftmargin=8mm]
    \item[(i)] The Riemann sphere admits an $\mathscr{F}$-invariant partition into a full, compact set $K$ and a simply connected domain $\Omega$.
    \item[(ii)] The restriction of $\mathscr{F}$ to $K$ is conformally conjugate to the action of the holomorphic dynamical system generated by the maps $P^{-m}\circ P^n$, $m,n\geq 0$, on the filled Julia set of $P$.
    \item[(iii)] The restriction of $\mathscr{F}$ to $\Omega$ is conformally conjugate to the action of $G$ on $\D$.
    \item[(iv)] The holomorphic dynamical system $\mathscr{F}$ is generated by the iteration of a holomorphic correspondence. This correspondence is transcendental when $G$ is uniform as well as when  $\D/G$ is a punctured 2-manifold of positive genus.
\end{enumerate}  
Moreover, if $G$ is either a uniform or a higher-genus non-uniform lattice, and $F$ is a mating between a polynomial $P$ and an arbitrary piecewise M{\"o}bius circle covering $\cH_G'$ that is weakly or strongly orbit equivalent to $G$, then the analytic continuations of $F$ \emph{cannot} give rise to an algebraic correspondence.
\end{thmx}
\noindent (See Theorems~\ref{uniform_lattice_holo_dyn_syst_thm},~\ref{non_uniform_lattice_holo_dyn_syst_thm}, and Propositions~\ref{non_alg_prop_uniform},~\ref{non_alg_prop_nonuniform},~\ref{alg_is_impossible_prop} for precise details.)

Thus, Fatou's vision from \cite{Fatou29} cannot literally be realized. We need to allow certain holomorphic correspondences over and above algebraic correspondences in order to combine general Fuchsian lattices with polynomial maps.

\subsection*{Bers slice of matings: analytic compactification, and boundary points}

For a fixed polynomial $P$ in the principal hyperbolic component (i.e., $P$ is a hyperbolic polynomial with a Jordan curve Julia set), the collection of all holomorphic dynamical systems realizing matings of $P$ with all Fuchsian lattices in the Teichm{\"u}ller space of a surface is called a \emph{Bers slice} through $P$ (by analogy with classical Bers slices in quasi-Fuchsian loci). For punctured spheres, the topological and analytic structures of such Bers slices were studied in \cite{MM2,LLM24,LMM26a}. Of particular importance is a result in \cite{LMM26a}, where compactifications of these Bers slices were constructed in the ambient \emph{character variety} of algebraic correspondences of a fixed bidegree (equipped with the algebraic topology of coefficient-wise convergence). 

In the current setting, lack of algebraicity of the holomorphic dynamical systems (arising as matings) forces us to equip these non-classical Bers slices with a weaker topology of Carath{\'e}odory convergence (see \S \ref{compactification_sec}). Precompactness of Bers slices in this topology is a consequence of standard results on normal families of holomorphic maps. 

With these structures in place, \S \ref{bers_bdry_blaschke_mating_sec} provides the following partial description of such a non-classical Bers boundary in the case of uniform lattices (see Lemma~\ref{bers_bdry_group_blaschke_mv_mating_lem}, Theorem~\ref{bers_bdry_group_blaschke_hds_mating_thm}, and Proposition~\ref{mating_cont_bers_bdry_prop} for precise statements).
 
\begin{thmx}\label{degenerate_group_mating_intro_thm}
Let $G_0$ be a uniform Fuchsian lattice with signature $(g,n) \notin \{(0,3), (1,1)\}$ and $P$ a polynomial in a principal hyperbolic component of appropriate degree. Then the boundary of the Bers slice of matings between $P$ and $G\in \mathrm{Teich}(G_0)$ contains holomorphic dynamical systems realizing matings of $P$ with all Kleinian groups on the (standard) Bers boundary of~$G_0$.
\end{thmx}

We note that Theorem~\ref{degenerate_group_mating_intro_thm} provides the first example of a mating between a polynomial and a simply degenerate group.  We expect that the matings described in this theorem exhaust the Bers boundary constructed in \S \ref{compactification_sec}; however we cannot establish this yet.

\subsection*{Proof techniques and discussions}

\subsubsection*{Pinching, pullback, and acceleration}

The construction of accelerated Bowen--Series maps goes through a sequence of steps.
Our starting point is the construction of \emph{Basilica Bowen--Series maps} and \emph{pullback Bowen--Series maps} in \cite{LMM26b}. 
\begin{enumerate}[leftmargin=6mm]
    \item[(i)] When $G$ is non-uniform, the existence of parabolics facilitates the construction of an appropriate ideal polygon (which is a fundamental domain when $G$ is torsion-free) such that the associated Bowen--Series map has controlled discontinuities (see \S \ref{ss:idealcorpolygon}). 
    \item[(ii)] For a uniform lattice $G$ satisfying the conditions of Theorem~\ref{most_fuchsian_mateable_intro_thm}, we first create a parabolic so that the problem can be partly reduced to the non-uniform case. More precisely, we go to the boundary of the Bers slice of $G$, where an essential separating simple closed curve $\gamma$ on the surface $\D/G$ is pinched (i.e., the associated element of $G$ becomes an accidental parabolic). The limit set of the resulting geometrically finite Kleinian surface group is topologically a Basilica, allowing us to use the technology of Basilica Bowen--Series maps developed in \cite{LMM26b}. The Basilica Bowen--Series map is a conformal Markov map on the Basilica limit set which is built by piecing together appropriate Bowen--Series maps of the components of the pinched surface. We then pull back the Basilica Bowen--Series map to the circle (the limit set of $G$ acting on the hyperbolic plane) via the associated Cannon-Thurston map, thereby obtaining the so-called pullback Bowen--Series map (see \S \ref{ss:bbspbs}). 
\end{enumerate}

A basic difficulty with pullback Bowen--Series maps is that they have discontinuities, though mild. To get rid of these discontinuities, we need to introduce the notion of \emph{acceleration} in \S \ref{acceleration_sec}. This gives rise to the accelerated Bowen--Series maps $\cH_G$ mentioned in Theorem~\ref{acc_bs_intro_thm} (obtained by accelerating the pullback Bowen--Series maps). The maps $\cH_G$ are, in fact, conformal Markov covering maps of the circle. We emphasize that for uniform lattices, the pullback Bowen--Series map is designed so that the one-sided limits at the points of discontinuity lie in the same $G$-orbit, ensuring continuity after suitable acceleration. (In contrast, the classical Bowen--Series maps of Fuchsian closed surface groups have irreparable discontinuities.) 

\subsubsection*{Quasiconformal surgery for uniform lattices, David for cusps}

The accelerated Bowen--Series maps $\cH_G$ constructed in Theorem~\ref{acc_bs_intro_thm} have the necessary analytic properties required to conformally mate them with complex polynomials. Indeed, the accelerated Bowen--Series map $\cH_G$ of a uniform Fuchsian lattice $G$ is quasiconformally conjugate to a power map on a pinched neighborhood of the circle. This allows us to apply standard quasiconformal surgery techniques (in particular, the measurable Riemann mapping theorem) to construct a degenerate polynomial-like map $F$ with external map $\cH_G$ and any prescribed polynomial $P$ (with connected Julia set) as its hybrid class (cf. \cite{DH85}, \cite[Appendix~A]{LLM24}). 

Due to the existence of parabolic fixed points, the accelerated Bowen--Series maps $\cH_G$ of non-uniform lattices are not quasisymmetrically equivalent to power maps; however, the David extension result of \cite{LMMN} (cf. \cite{LN26}) shows that the circle homeomorphism conjugating a power map to $\cH_G$ extends as a David homeomorphism of the circle. One then invokes the David integrability theorem to construct the desired degenerate polynomial-like maps $F$ with external map $\cH_G$.
However, due to certain analytic limitations, this process cannot be carried out for all polynomials with connected Julia set (cf. \cite{MM2,LLM24}).

The quasiconformal compatibility between accelerated Bowen--Series maps of uniform lattices and power maps also plays a pivotal role in the mating of the corresponding Bers boundary groups with polynomials (see  Theorem~\ref{degenerate_group_mating_intro_thm} and Lemma~\ref{bers_bdry_group_blaschke_mv_mating_lem}). The lack of such compatibility in the non-uniform setting is one of the main obstacles to constructing matings between polynomials and Kleinian groups on the Bers boundaries of non-uniform Fuchsian lattices, other than those uniformizing punctured spheres (cf. \S \ref{compare_sec}).

\subsubsection*{Strong vs weak orbit equivalence} 
\begin{figure}[htp]
 \captionsetup{width=0.96\linewidth}
    \centering
    \includegraphics[width=0.9\textwidth]{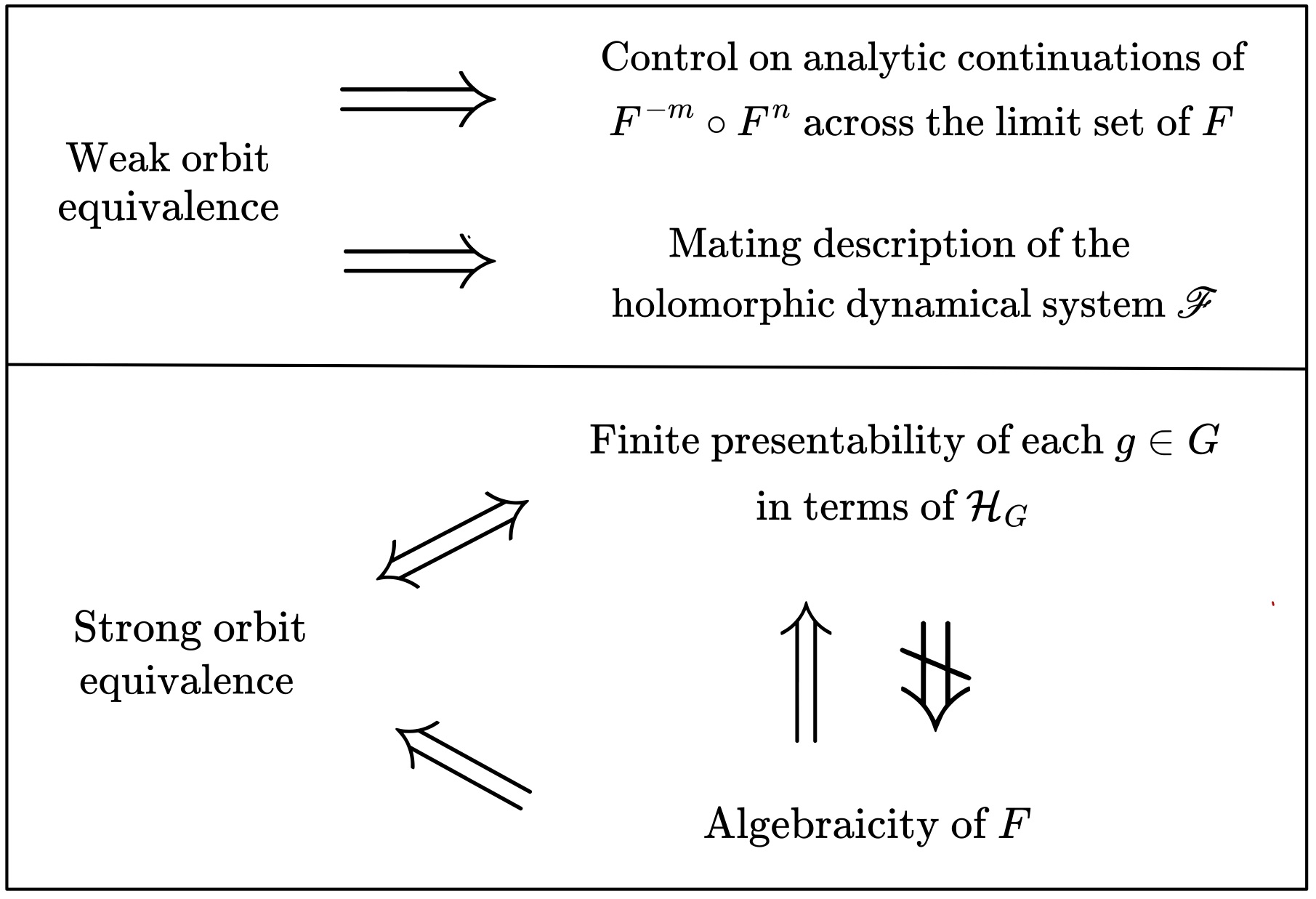}
    \caption{Implications of the two notions of orbit equivalence for the mating structure and algebraicity of the holomorphic dynamical system~$\mathscr{F}$.}
    \label{oe_chart_fig}
\end{figure}

The multi-step construction of an accelerated Bowen--Series map (involving pinching, pullback, and acceleration) comes at a cost: it is generally not orbit equivalent to the corresponding Fuchsian group on the whole circle. This forces us to slightly relax the orbit equivalence framework of \cite{MM1,MM2,LLM24}. We analyze the orbit equivalence properties of pullback and accelerated Bowen--Series maps in \S \ref{sec:orbit_equiv} and \S \ref{acceleration_sec}, establishing that accelerated Bowen--Series maps are always weakly orbit equivalent to the corresponding groups (i.e., orbit equivalent on full-measure subsets of the circle). Further, strong orbit equivalence (i.e., orbit equivalence at each point of the circle) holds for torsion-free non-uniform lattices, but fails for uniform ones (this failure is a consequence of the pinching and pullback process). 
This is central to our investigation of analytic continuations of the degenerate polynomial-like maps arising from Theorem~\ref{most_fuchsian_mateable_intro_thm}. We expound on this connection below (see Figure~\ref{oe_chart_fig}).

\subsubsection*{Weak orbit equivalence implies a mating description for holomorphic dynamical systems}
A degenerate polynomial-like map $F$ realizing the mating of a polynomial $P$ and an accelerated Bowen--Series map $\cH_G$ is promoted to a holomorphic dynamical system $\mathscr{F}$ by analytically continuing $F^{-m}\circ F^n$, $m,n\geq 0$. The weak orbit equivalence between an accelerated Bowen--Series map $\cH_G$ and the corresponding Fuchsian lattice $G$ ensures that each element of the group $G$ can be recovered from $\cH_G$ on an open, dense, full-measure subset of the circle (see Proposition~\ref{countable_presentation_prop}). This weak recoverability property enables us to control the analytic continuations of $F^{-m}\circ F^{n}$ thereby
\begin{enumerate}[leftmargin=8mm]
    \item[(i)] establishing an $\mathscr{F}$-invariant partition of the Riemann sphere, and
    \item[(ii)] ruling out the existence of \emph{rogue} analytic continuations by showing that all analytic continuations of $F^{-m}\circ F^{n}$ correspond either to the polynomial $P$ or to the group $G$.
\end{enumerate} 
Thus, the mating description of $\mathscr{F}$ (see Theorem~\ref{holo_dyn_syst_intro_thm}) relies heavily on the weak orbit equivalence between $\cH_G$ and $G$.

\subsubsection*{Strong orbit equivalence, periodicity of monodromy action on stabilizers, and algebraicity}

A necessary condition for algebraicity of the degenerate polynomial-like map $F$ is that its maximal analytic continuation has at most finitely many branch points each of which has finite local monodromy. On the other hand, strong orbit equivalence between $\cH_G$ and $G$ is characterized by a different finiteness property: for each $g\in G$, there is a finite partition of $\bS^1$ into intervals on each of which $g\equiv\cH_G^{-m}\circ\cH_G^n$ for some $(m,n)$ depending on the interval. We call this the \emph{finite presentability} of $G$ in terms of $\cH_G$; it is a strong form of recoverability of the group $G$ from the map $\cH_G$.

It turns out that the finite presentability property (equivalently, strong orbit equivalence between $\cH_G$ and $G)$ is a necessary condition for finiteness of the number of branch points for the maximal analytic continuation of $F$. Thus, for uniform lattices, failure of strong orbit equivalence gives infinitely many branch points and hence transcendentality of the degenerate polynomial-like map $F$ (see Theorem~\ref{thm:weakorbitequiv} and Proposition~\ref{non_alg_prop_uniform}).

However, the finite presentability property does not imply algebraicity.
This is because finiteness of local monodromies around branch points is a necessary condition for algebraicity, and it does not follow from finite presentability. Indeed, by definition, the map $\cH_G$ is also finitely expressible as elements of $G$. The finite presentations of $G$ in terms of $\cH_G$ and that of $\cH_G$ in terms of $G$ together induce a monodromy action on the stabilizers of the break-points of $\cH_G$. Finite-order local monodromy around the branch points of $F$ amounts to periodicity of these monodromy actions on stabilizers (see Propositions~\ref{non_alg_prop_nonuniform} and~\ref{alg_is_impossible_prop}). This periodicity condition breaks down for accelerated Bowen--Series maps $\cH_G$ of torsion-free non-uniform lattices of positive genus (although such a group $G$ is strongly orbit equivalent to $\cH_G$). Consequently, the holomorphic dynamical systems realizing matings between polynomials and Fuchsian groups uniformizing higher-genus punctured surfaces also fail to be algebraic.

In summary, in our construction involving accelerated Bowen--Series maps, the matings arising from closed surfaces are not algebraic as they have infinitely many branch points (Proposition~\ref{non_alg_prop_uniform}) and those arising from punctured surfaces of higher genus are not algebraic as they violate the finite local monodromy condition at some branch points (Proposition~\ref{non_alg_prop_nonuniform}).

More generally (as made explicit in the last statement of Theorem~\ref{holo_dyn_syst_intro_thm}), we show in Proposition~\ref{alg_is_impossible_prop} that this failure of algebraicity is not an artifact of our construction: for uniform or higher-genus non-uniform lattices, the above monodromy actions on stabilizers of break-points are always aperiodic (independent of whether the matings arise from accelerated Bowen--Series maps).

\begin{remark}
In the present paper and those leading up to it, we have focused on a number of variants of the Bowen--Series map. We do not know if there are other (possibly completely different) piecewise M{\"o}bius Markov circle coverings that are strongly orbit equivalent to uniform lattices (although this cannot be upgraded to algebraicity of matings by Proposition~\ref{alg_is_impossible_prop}).
\end{remark}

\begin{remark}
(1) The non-algebraicity of the mating of the polynomial $P$ and the Fuchsian group $G$, even though both are algebraic, is morally analogous to the group-theoretic phenomenon that two finitely presented groups can generate a group that is not finitely presentable.

(2) The algebraicity of matings in the punctured sphere examples of \cite{MM2,LLM24} can be traced back to the fact that the mating $F$ acts as an involution on the boundary of its domain of definition (see \cite{LLM24} for details on the theory of \emph{B-involutions}, which abstracts out this structure). In fact, the mating $F$ inherits this property from the involutive structure of Bowen--Series maps for Fuchsian punctured sphere groups. This forces the local monodromies of the analytic continuation of $F$ (corresponding to the parabolic break-points of the Bowen--Series map) to have order two, which is necessary for algebraicity according to Proposition~\ref{alg_is_impossible_prop}.
\end{remark}

\subsubsection*{Finite-dimensional Bers compactification}
While the Bers slices of matings are finite-dimensional (indeed, they are naturally homeomorphic to finite-dimensional Teichm{\"u}ller spaces of Fuchsian lattices), their compactifications constructed in \S \ref{compactification_sec} live in infinite-dimensional spaces of holomorphic dynamical systems. In the spirit of \cite{LMM26a}, it would be interesting to parametrize these Bers slices by finitely many coordinates, and establish their precompactness in the resulting finite-dimensional ambient spaces.

\subsection*{Acknowledgement} 
We thank Mikhail Lyubich and Curtis McMullen for valuable discussions.

\section{Pull-back Bowen--Series maps for Fuchsian lattices}\label{pbs_sec}
In this section, we summarize the construction of the pull-back Bowen--Series map in \cite{LMM26b}.

\subsection{Conformal Markov maps on $\SS^1$}\label{sec:cmm}
In this section, we summarize properties of conformal Markov maps on $\SS^1$. They are a special case of $(X, \mathscr{C})$-map studied in \cite{LMM26b}.
\begin{defn}\label{XCstrs}
    A finite collection $\mathcal{A} = \{A_1,\cdots, A_m\}$ of subsets of $\SS^1$ is called a {\em partition} if
    \begin{itemize}
        \item $A_i$ is a closed interval;
        \item $\Int A_i \cap \Int A_j = \emptyset$;
        \item $\bigcup A_i = \SS^1$.
    \end{itemize}
    We always assume the intervals $A_i$ are ordered and oriented counterclockwise on the circle.
    We define the set of {\em break-points} as 
    $$
    \mathrm{bk}(\mathcal{A}):= \{x \in X: x \in A_i\cap A_{i+1}\}.
    $$

    A finite collection $\mathcal{F}=\{f_1,\cdots, f_m\}$ is called a {\em conformal Markov map} if
    \begin{itemize}[leftmargin=8mm]
    \item for each $i$, there exists $j_i$ and $l_i$ so that 
    $$
    f_i: A_i \longrightarrow A_{j_i} \cup A_{j_i+1} \cup \ldots \cup A_{j_i+l_i}
    $$
    is a homeomorphism;
    \item $f_i$ extends to a conformal map on a neighborhood of $A_i$ in $\widehat{\C}$.
    \end{itemize}
\end{defn}
We do not require the collection $\mathcal{F}$ to induce a continuous map on $\SS^1$. In other words, if $x \in A_i \cap A_{i+1} \subseteq \mathrm{bk}(\mathcal{A})$, then $f_i(x)$ may be different from $f_{i+1}(x)$.
    It is convenient to regard $\mathcal{F}$ as a multi-valued map at the break-points. However, the Markov property implies that both $f_i(x)$ and $f_{i+1}(x)$ belong to $\mathrm{bk}(\mathcal{A})$.
We use the following notation:
$$
    \mathcal{F}(x^+) = \lim_{t \to x^+} \mathcal{F}(t), \text{ and } \mathcal{F}(x^-) = \lim_{t \to x^-} \mathcal{F}(t).
$$
Similarly, we denote $\mathcal{F}^n(x^\pm) := \lim_{t \to x^\pm} \mathcal{F}^n(t)$.

With this notation, two conformal Markov maps $\mathcal{F}, \cG$ are {\em topologically conjugate} if there exists a homeomorphism $h: \mathbb{S}^1 \longrightarrow \mathbb{S}^1$ so that
    $$
    \mathcal{F}(h(x)^\pm) = h(\cG(x^\pm).
    $$

\subsection*{Markov pieces and topological expansion}
A partition $\mathcal{B}$ of $\SS^1$ is a {\em refinement} of $\mathcal{A}$ if for each $B \in \mathcal{B}$, there exists $A \in \mathcal{A}$ so that $B \subseteq A$.
Given a conformal Markov map $\mathcal{F} = \{f_1,\cdots, f_m\}$, one can define the pull-back partition 
$$
\mathcal{A}^1 = \mathcal{F}^*(\mathcal{A}) = \{A_{ij} :=f_i^{-1}(A_j): f_i(A_i) \cap A_j \neq \emptyset\}.
$$
We call $\mathcal{A}^0 = \mathcal{A}$ (respectively, $\mathcal{A}^1$) the {\em level $0$} (respectively, {\em level $1$}) Markov partitions.
By construction, 
\begin{itemize}[leftmargin=8mm]
    \item(inclusion) $A_{ij} \subseteq A_i$, i.e., $\mathcal{A}^1$ is a refinement of $\mathcal{A}^0$; and
    \item(dynamics) $f_i: A_{ij} \longrightarrow A_j$ is a homeomorphism.
\end{itemize}
We define the level $n$ Markov partitions inductively, and denote them by $\mathcal{A}^n = \{A_{w}\}$, where $w$ is a word of length $n$ consisting of letters in $\{1,\cdots, m\}$. We will call a word $w$ of length $n$ {\em admissible} if $A_w \in \mathcal{A}^n$.

It is convenient to write a conformal Markov map as
$$
\mathcal{F} : \mathcal{A}^1 \longrightarrow \mathcal{A}^0.
$$
We say that a conformal Markov map $\mathcal{F}$ is {\em topologically expanding} if
$\bigcap_{w_n} A_{w_n}$ is a singleton set for each nested infinite sequence of level $n$ Markov pieces.

\subsection*{Markov matrix and primitivity}
Let $\mathcal{F} : \mathcal{A}^1 \longrightarrow \mathcal{A}^0$ be a conformal Markov map.
One can record the combinatorics of the Markov map by an $m \times m$ matrix, called the {\em Markov matrix} $M = M_\mathcal{F}$. More precisely, 
$$
M_{ij} = \begin{cases}
    1 & \text{ if }A_j \subseteq f_i(A_i)\\
    0 & \text{ otherwise.}
\end{cases}.
$$
Note that the Markov matrix requires a choice of an ordering of the elements in the partition $\mathcal{A}^0$. Thus, $M_{\mathcal{F}}$ is well-defined up to a cyclic conjugation.

The map $\cF$ is said to be {\em primitive} if there exists $k$ so that for each $A \in \mathcal{A}^1$, $\cF^k(A) = \SS^1$. This is equivalent to the fact the Markov matrix $M_\cF$ is {\em primitive}, i.e., there exists $k$ so that $M_\cF^k$ is a positive matrix.

\subsection*{Combinatorial and topological conjugation}
Two conformal Markov maps $\mathcal{F}: \mathcal{A}^1 \longrightarrow \mathcal{A}^0$ and  $\cG: \mathcal{B}^1 \longrightarrow \mathcal{B}^0$ are {\em combinatorially conjugate} if $M_\mathcal{F}$ is cyclically conjugate to $M_\mathcal{G}$. More precisely, after possibly cyclically permuting the indices, $M_\mathcal{F} = M_\mathcal{G}$.
The following lemma follows from \cite[Proposition~2.7]{LMM26b}.
\begin{lem}\label{lem:comb=>top}
    Let $\mathcal{F}: \mathcal{A}^1 \longrightarrow \mathcal{A}^0$ and  $\cG: \mathcal{B}^1 \longrightarrow \mathcal{B}^0$ be two conformal Markov maps.
    Suppose that they are topologically expanding and combinatorially conjugate. Then they are topologically conjugate.
\end{lem}

\subsection*{Break-points and Lyapunov exponent}
Let $\mathcal{F}$ be a topologically expanding conformal Markov map on $\bS^1$ with respect to the partition $\mathcal{A}$. 
Let $x \in A_i \cap A_{i+1}$ be a break-point for the dynamics of $\mathcal{F}$.
Denote by $x_n^\pm$ the right and left orbit of $x$; i.e., 
$$
x_n^\pm = \mathcal{F}^n(x^\pm).
$$
Since the partition is finite, we note that both right and left orbits are pre-periodic. Note that they may be eventually mapped to different periodic orbits.
Let us denote the right and left {\em Lyapunov exponent} by $\lambda(x^+)$ and $\lambda(x^-)$ respectively; i.e.,
$$
\lambda(x^\pm) = \lim_{n\to\infty} \frac{1}{n} \log \lim_{\epsilon \to 0} \frac{|\mathcal{F}^n(x\pm\epsilon)- x_n^\pm|}{\epsilon}.
$$
Since $\mathcal{F}$ is topologically expanding, we have $\lambda(x^\pm) \geq 0$.
We say that $x^+$ (respectively, $x^-$) is {\em parabolic} if $\lambda(x^+) = 0$ (respectively, $\lambda(x^-) = 0$). 
Similarly, we say that $x^+$ (respectively, $x^-$) is {\em hyperbolic} if $\lambda(x^+) > 0$ (respectively, $\lambda(x^-) > 0$).

Since $x_n^\pm$ are eventually periodic, $x^\pm$ is parabolic or hyperbolic if and only if the periodic tail of the sequence $x_n^\pm$ is a parabolic or hyperbolic cycle.
If $x^\pm$ is parabolic, we define its {\em parabolic multiplicity}, denoted by $N(x^\pm)$, as the parabolic multiplicity of the periodic cycle it is mapped to (cf. \cite[\S 10]{Mil06}).

We say that the break-point $x$ is {\em symmetric} if 
\begin{itemize}
    \item $q^+\lambda(x^+) = q^-\lambda(x^-)$, where $q^\pm$ are the periods of the eventually periodic orbits $x^\pm$,
    \item in addition, if  $x^\pm$ are parabolic, then 
    $N(x^+) = N(x^-)$.
\end{itemize} 
Note that in the symmetric case, we will call the break-point $x$ \emph{symmetrically parabolic/hyperbolic}.

\subsection*{Regular points and grand orbit}
Let $\mathcal{F} : \mathcal{A}^1 \longrightarrow \mathcal{A}^0$ be a conformal Markov map. 
We define the {\em orbit} of a point $x \in \mathbb{S}^1$ as the union of both right and left orbits of $x$, i.e., 
$$
\mathcal{O}_\mathcal{F}(x) := \{\mathcal{F}^n(x^\pm):n \geq 0\}.
$$
We call a point $x$ a {\em break-point of $\cF$} if its orbit intersects $\mathrm{bk}(\mathcal{A})$.
We denote this set by $\mathrm{bk}(\cF)$.
A point in $\mathbb{S}^1 \setminus \mathrm{bk}(\cF)$ is called a {\em regular} point.

Suppose that $x\in \mathbb{S}^1$ regular. Then we define the grand orbit of $x$ by
$$
\mathcal{GO}_{\mathcal{F}}(x) = \{y\in \mathbb{S}^1: \mathcal{O}_\mathcal{F}(x) \cap \mathcal{O}_\mathcal{F}(y) \neq \emptyset\}.
$$
This defines an equivalence relation on the regular set, and we write
$$
x\sim_{\mathcal{F}} y
  \quad\Longleftrightarrow\quad
  y\in \operatorname{GO}_{\mathcal{F}}(x).
$$
For simplicity of presentation, we shall not consider the grand orbits of break-points under discontinuous conformal Markov maps in this paper.

\subsection*{Puzzle structure}
Let $X \subseteq \SS^1$ be a closed interval. A subset $\widehat{X}$ of $\widehat{\C}$ is said to be a \emph{pinched neighborhood} of $X$, if
\begin{itemize}
    \item $X\subset\widehat{X}$;
    \item $\Int(\widehat{X})$ is an open set that contains the open interval $\Int(X)$;
    \item $\partial X \subseteq \partial \widehat{X}$.
\end{itemize}

\begin{defn}[Puzzle structure]
    Let $\mathcal{F}: \mathcal{A}^1 \longrightarrow \mathcal{A}^0$ be a conformal Markov map on $\SS^1$.
    We say that it admits a {\em puzzle structure} if there exist closed pinched neighborhoods $P_i$ of $A_i$, pinched at the two end points $\partial A_i$, so that for each level $1$ Markov piece $A_{ij} \subseteq A_i$, we have
    $$
        f_i: \Int P_{ij} = f_i^{-1}(\Int P_j) \longrightarrow \Int P_j
    $$
    is conformal, and $\Int P_{ij} \subseteq \Int P_{i}$.
\end{defn}

\subsection{Bowen--Series map for non-uniform Fuchsian lattice}\label{ss:idealcorpolygon}
In this subsection, we associate a special ideal polygon to each non-uniform Fuchsian lattice.

Let $\Sigma$ be a non-compact (finite area) hyperbolic surface with signature $(g, n; \nu_1,\ldots, \nu_n)$.
Let us arrange so that $\nu_1 = \nu_2 = \ldots = \nu_k = \infty$, and $\nu_{k+1}, \ldots, \nu_{n} < \infty$.
Since $\Sigma$ is non-compact, we have that $k \geq 1$.

\begin{figure}[htp]
 \captionsetup{width=0.96\linewidth}
    \centering
    \includegraphics[width=0.96\textwidth]{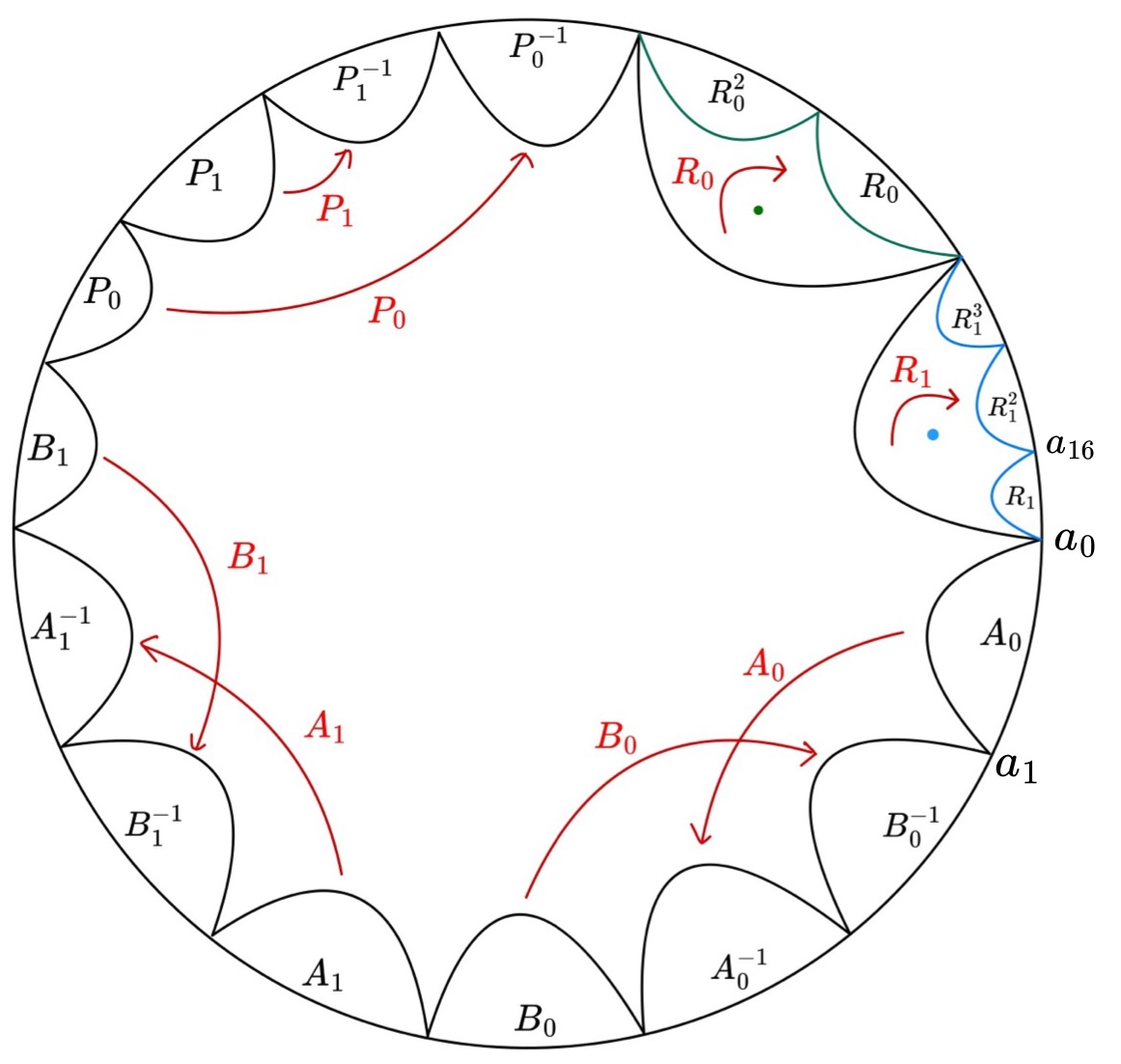}
    \caption{Depicted is a core polygon and the associated Bowen--Series map for a Fuchsian group of type $(2, 5; \infty, \infty, \infty, 3,4)$.}
    \label{fig:cuspBM}
\end{figure}

Let $G$ be a Fuchsian group with the given signature. This means that the quotient $S = \D / G$ is a hyperbolic surface homeomorphic to $\Sigma$.
We define a {\em core polygon} $\cP$ for $G$ as follows (see Figure~\ref{fig:cuspBM}). For $x,y\in\mathbb{S}^1$, we denote the hyperbolic geodesic in $\D$ connecting $x, y$ by $\arc{x,y}$.

\begin{itemize}[leftmargin=8mm]
    \item $\cP$ is an ideal hyperbolic $N$-gon in $\D$, where $N = 4g+2(k-1)+\sum_{i=k+1}^n(\nu_i-1)$.
    \item Denote the vertices of $\cP$ by $a_0,\ldots, a_{N-1} \in \mathbb{S}^1$ (labeled clockwise). Then each edge $\gamma_i,\ i\in\{0,\ldots, N-1\}$, is a hyperbolic geodesic in $\D$ with end-points $a_i, a_{i+1}$ (with indices mod $N$).
    \item ($g$ handles) There exist $2g$ loxodromic elements $A_i, B_i \in G,\  i\in\{0,\ldots, g-1\}$, such that for $i\in\{0,\cdots,g-1\}$, 
    $$
    A_i(\arc{a_{4i},a_{4i+1}}) = \arc{a_{4i+3}, a_{4i+2}}, \text{ and } B_i(\arc{a_{4i+3}, a_{4i+4}}) = \arc{a_{4i+2}, a_{4i+1}}.
    $$
    \item ($k$ punctures) There exist $k-1$ elements $P_0, \ldots, P_{k-2} \in G$ such that
    $$
    P_i(\arc{a_{4g+i}, a_{4g+i+1}}) = \arc{a_{4g+2(k-1)-i}, a_{4g+2(k-1)-i-1}}.
    $$
    Moreover, if $k\geq 4$, then $P_1,\ldots, P_{k-3}$ are loxodromic and $P_{k-2}$ is parabolic. The element $P_0$ is parabolic if $(g,n) = (0, k)$, i.e., if the surface is a punctured sphere, and is loxodromic otherwise. For $i\in\{0,\ldots, k-3\}$, the group element $P^{-1}_{i+1} \circ P_i$ is parabolic. (The cases $k\in\{1,2,3\}$ can be described similarly.)
    \item ($n-k$ orbifold points) There exist $n-k$ elliptic elements $R_0,\ldots, R_{n-k-1} \in G$ such that for $i\in\{0,\ldots, n-k-1\}$, we have $R_i^{\nu_{k+1+i}} = \text{id}$, and
    $$
    R_i(\arc{a_{i,j}, a_{i,j+1}}) = \arc{a_{i,j+1}, a_{i,j+2}},
    $$
    where $j\in\{0,\ldots, \nu_{k+1+i}-1\}\ \left(\textrm{mod}\ \nu_{k+1+i}\right)$. Here,
    \begin{align}
      a_{i,j} = 
      \begin{cases}
          a_{4g+2(k-1)+\sum_{s=k+1}^{k+i}(\nu_s-1)+j}\quad \mathrm{for}\ i\geq 1,\\
          a_{4g+2(k-1)+j}\hspace{2.4cm} \mathrm{for}\ i=0.
      \end{cases}  
      \label{a_i_j_formula}
    \end{align}
\end{itemize}

Note that if $G$ is torsion-free, then $\cP$ is a fundamental polygon for the action of $G$ on $\D$.

The vertex set induces a partition $\mathcal{A} = \{[a_i,a_{i+1}]:i\in\{0,\ldots, N-1\}\}$  of $\mathbb{S}^1$.
We define the {\em Bowen--Series Markov map} associated to the core polygon $\cP$ as a Markov map with respect to the partition $\mathcal{A}$:
$$
\mathcal{F} = \mathcal{F}_{G, \cP} = \{F_i:=\mathcal{F}|_{[a_i, a_{i+1}]}: i\in\{0,\ldots, N-1\}\},
$$ 
where $\mathcal{F}|_{[a_i, a_{i+1}]}$ is defined as follows (see Figure~\ref{fig:cuspBM}, cf. \cite{BS79}).

\begin{itemize}
    \item For $i\in\{0,\ldots, g-1\}$, 
    \begin{align*}
        \mathcal{F}|_{[a_{4i},a_{4i+1}]} &= A_i,\\
        \mathcal{F}|_{[a_{4i+1},a_{4i+2}]} &= B_i^{-1},\\
        \mathcal{F}|_{[a_{4i+2},a_{4i+3}]} &= A_i^{-1},\\
        \mathcal{F}|_{[a_{4i+3},a_{4i+4}]} &= B_i.
        \end{align*}
    \item For $i\in\{0,\ldots, k-1\}$,
    \begin{align*}
        \mathcal{F}|_{[a_{4g+i},a_{4g+i+1}]} &= P_i,\\
        \mathcal{F}|_{[a_{4g+2(k-1)-i-1},a_{4g+2(k-1)-i}]} &= P_i^{-1}.
    \end{align*}
    \item For $i\in\{0,\ldots, n-k-1\}$, $j\in\{0,1,\ldots, \nu_{k+1+i}-2\}$,
    $$
    \mathcal{F}|_{[a_{i,j}, a_{i,j+1}]} = R_i^{\nu_{k+1+i}-1-j} 
    $$
    where $a_{i,j}$ is defined by Formula~\eqref{a_i_j_formula}.
\end{itemize}
We remark that $\mathcal{F}$ is a conformal Markov map on $\mathbb{S}^1$ in the sense of \S \ref{sec:cmm}. It induces a continuous map on the circle if and only if $\Sigma$ has signature $(0,n;\infty,\cdots,\infty)$ or $(0,n;\infty,\cdots,\infty,2)$; i.e., the surface is a punctured sphere possibly with an order $2$ orbifold point (cf. \cite[\S 3]{MM1}).

The following properties follow from the construction.

\begin{prop}\label{prop:symFs}
    The Bowen--Series Markov map $\mathcal{F} = \mathcal{F}_{G, \cP}$ is topologically expanding, primitive, conformal, and admits a natural puzzle structure: for each Markov interval $[a_i, a_{i+1}] \in \mathcal{A}$, the corresponding level 0 puzzle $P_i$ is the disk bounded by the circle that meets $\mathbb{S}^1$ orthogonally at $a_i, a_{i+1}$. 
    
    Moreover, every break-point is symmetrically parabolic.
    More precisely, let $a$ be a periodic break-point. Then the stabilizer of $a$ in $G$ is generated by a parabolic element $g \in G$ so that $g$ is repelling on the circular arc $[a, a+\epsilon]$ and attracting on $[a-\epsilon, a]$.
    Further, if $q^\pm$ is the period of the right/left orbit $a_n^\pm$, respectively, then $\mathcal{F}^{q^\pm}|_{[a, a\pm \epsilon]} = g^{\pm 1}$ for all sufficiently small $\epsilon$.
\end{prop}

\subsection{Basilica Bowen--Series and pull-back Bowen--Series maps}\label{ss:bbspbs}
In this subsection, we give a natural puzzle structure for the limit set of any geometrically finite Bers boundary group.

Let $\Sigma$ be a surface with $\chi(\Sigma) < 0$ and let $\mathcal{C}$ be a pinching multicurve, with $(\Sigma, \mathcal{C}) = \Sigma_1+\Sigma_2+ \ldots + \Sigma_{j_0}$.
The components of $\Sigma \setminus \bigcup_{c \in \mathcal{C}} c$ are called  \emph{nodal components} \cite[\S 5]{LMM26b}. Each such nodal component is either hyperbolic or parabolic, where a parabolic nodal component necessarily has signature equal to $(0, 3; 2, 2, \infty)$. 
For simplicity of  presentation, let us first assume that all nodal components are hyperbolic.
 We will discuss necessary modifications to handle parabolic nodal components at the end in \S \ref{subsec:parabolicnodal}.

Let $G$ be a geometrically finite Bers boundary group associated to $(\Sigma, \mathcal{C})$.
Let $\Delta_1$ be a component of $\Omega(G) - \Delta_\infty$, where $\Delta_\infty$ is the $G$-invariant component.
Then $\Delta_1$ corresponds to some component $\Sigma_i$ which is hyperbolic. Without loss of generality, we assume that $\Delta_1$ corresponds to $\Sigma_1$.

\begin{figure}[htp]
 \captionsetup{width=0.96\linewidth}
    \centering
    \includegraphics[width=0.88\textwidth]{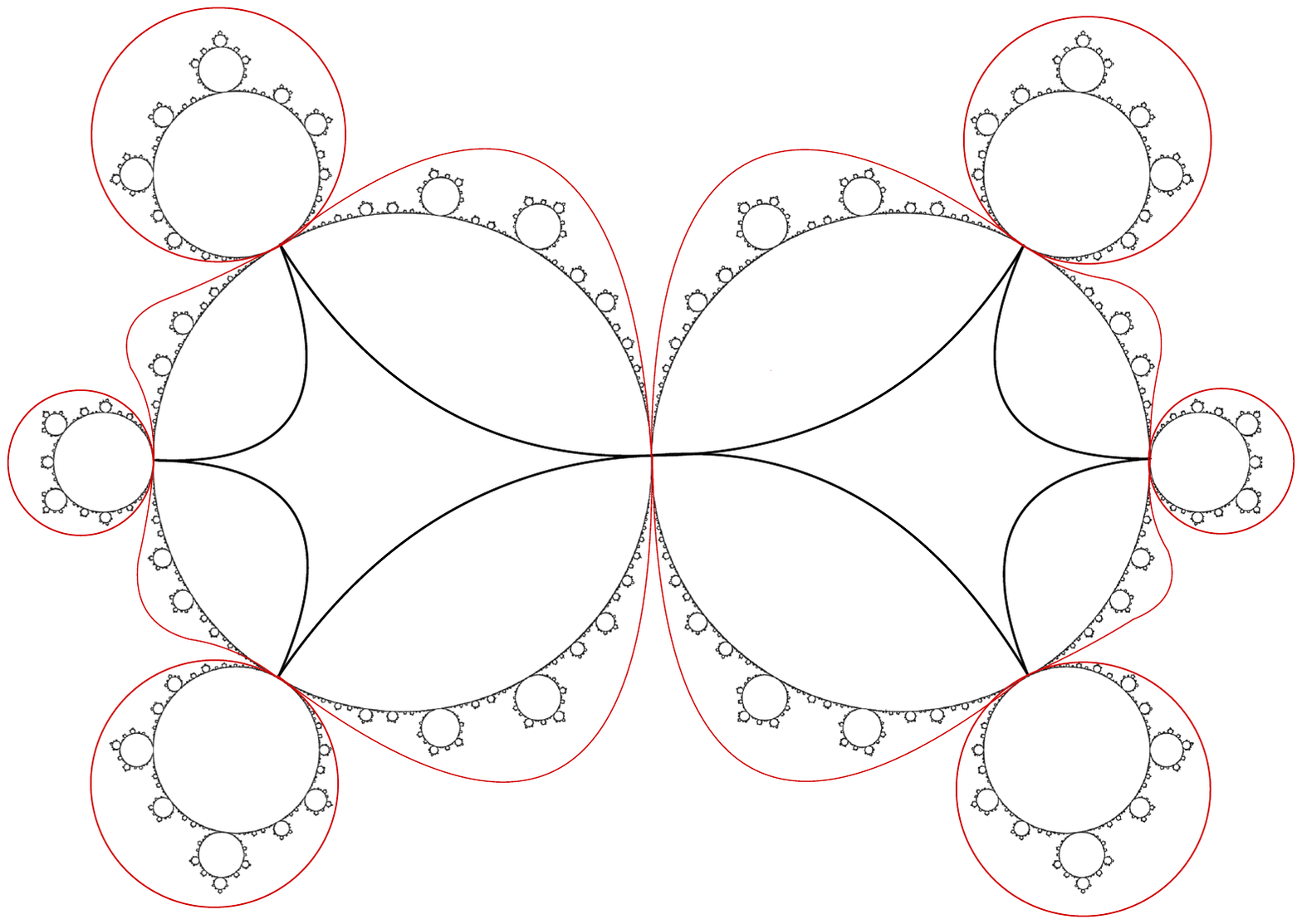}
    \caption{The puzzle structure for a genus $2$ Bers boundary group. The union of the two black ideal quadrilaterals $\overline{\cP_1} \cup \overline{\cP_2}$ is the pinched core polygon for $G$. The associated ideal polygon $\cP_\infty \subseteq \Delta_\infty$ is bounded by the union of red geodesic edges. The closures of the connected components of $\widehat{\C} - \left(\overline{\cP_1}\cup \overline{\cP_2} \cup \overline{\cP_\infty}\right)$ form the level $0$ puzzle pieces. A limb type puzzle piece is bounded by a red edge, while a torso type puzzle piece is bounded by the union of one black edge and one red edge.}
    \label{fig:PinchedPolygon}
\end{figure}

Since $\Sigma_1$ is non-compact, we can construct a core polygon $\cP_1 \subseteq \Delta_1$ following \S \ref{ss:idealcorpolygon}.
Let $G_{\Delta_1} \subseteq G$ be the stabilizer of $\Delta_1$.
Note that the $G_{\Delta_1}$-orbit of the ideal vertices of $\cP_1$ are in one-to-one correspondence with the cusps of $\Sigma_1$.
Thus, if $j_0>1$, there exists at least one component $\Delta_2 \left(\neq\Delta_1\right)$  of $\Omega(G) - \Delta_\infty$ associated to some surface in $\{\Sigma_2, \ldots, \Sigma_{j_0}\}$ that is adjacent to $\Delta_1$ at some vertex of $\cP_1$. Without loss of generality, we assume that $\Delta_2$ corresponds to the surface $\Sigma_2$.
Let $\cP_2 \subseteq \Delta_2$ be a core polygon.
Inductively, we can construct 
\begin{itemize}[leftmargin=8mm]
    \item a tree of components $\Delta_1,\ldots, \Delta_{j_0}$ corresponding to $S_1,\ldots, S_{j_0}$ (respectively) such that $\bigcup_{i=1}^{j_0} \overline{\Delta_i}$ is connected, and
    \item a collection of core polygons $\cP_1,\ldots,\cP_{j_0}$ so that $\cP_i \subseteq \Delta_i$ and $\bigcup_{i=1}^{j_0} \overline{\cP_i}$ is connected.
\end{itemize}
We will call this collection $\{\cP_1,\ldots, \cP_{j_0}\}$ a {\em pinched core polygon} for $G$ (see Figure~\ref{fig:PinchedPolygon}).

We remark that if $G$ is torsion-free, then one can realize this collection of ideal polygons as the limit of some fundamental polygon for the Fuchsian model $G_0$ of $\Sigma$. Specifically, pinching the corresponding multicurve $\mathcal{C}$ on $\Sigma$ yields a sequence of quasi-Fuchsian deformations of $G_0$ that degenerates to~$G$.

Recall that we denote the ideal boundary of $\Delta_\infty$ by $\partial^I \Delta_\infty$.
Each vertex $x$ of $\cP_i$ corresponds to a cusp of $\Sigma_i$ which is either a nodal singularity of $\overline{\Sigma_i}$ or a cusp of $\Sigma_i$ coming from the original surface $\Sigma$.
In the former case, $x$ corresponds to $2$ ideal boundary points for $\Delta_\infty$, and there is a component of $\Omega(G) - \Delta_\infty$ attached to $\Delta_i$ at~$x$.
In the latter case, $x$ corresponds to a unique ideal boundary point for $\Delta_\infty$.
We consider the ideal polygon $\cP_\infty \subseteq \Delta_\infty$ whose vertices are the ideal boundary points (for $\Delta_\infty$) associated to the vertices of $\bigcup_{i=1}^{j_0}\cP_i$. (In Figure~\ref{fig:PinchedPolygon}, the polygon $\cP_\infty$ in $\Delta_\infty$ is shown in red. Here, $\Sigma$ is a compact surface of genus $2$, and the multi-curve $\mathcal{C}$ has a single component; hence, pinching $\mathcal{C}$ produces two genus $1$ surfaces $\Sigma_1$ and $\Sigma_2$ each of which has a unique node and no other cusp. The sides of the polygon $\cP_\infty$ are lifts of the curve $\mathcal{C}$. In particular, the second possibility discussed above does not arise in this example.)

The complement $\widehat\C - \left(\bigcup_{i=1}^{j_0}\overline{\cP_i} \cup \overline{\cP_\infty}\right)$ consists of finitely many Jordan domains, whose closures are denoted by $P_1,\ldots, P_r$. We will refer to these as level 0 puzzle pieces, and denote the collection of such puzzle pieces by
$$
\mathfrak{P}^0 := \{P_1, \ldots, P_r\}.
$$

Note that there are two types of puzzle pieces.
We say that a puzzle piece $P_i$ is of {\em limb type} if $\Int P_i \cap \bigcup_{i=1}^{j_0}\Delta_i = \emptyset$, and call it {\em torso type} otherwise. 
If $P_i$ is of limb type, then it is bounded by some geodesic $\widetilde{\alpha}_i \subseteq \Delta_\infty$ together with its accumulation point, which is the fixed point of an accidental parabolic element in $G$.
If $P_i$ is of torso type, then it is bounded by two geodesics and their accumulation points, with one geodesic in $\Delta_\infty$, and the other in $\Delta_s$ for some $s \in \{1,\ldots, j_0\}$ (see Figure~\ref{fig:PinchedPolygon}).

\subsubsection{Modification for parabolic nodal components}\label{subsec:parabolicnodal}
Suppose that 
$$
(\Sigma, \mathcal{C}) = \Sigma_1+\Sigma_2+ \ldots + \Sigma_{l_0} + \Sigma_{l_0+1} + \ldots + \Sigma_{j_0},
$$
where $\Sigma_i$ is hyperbolic for $i \leq l_0$ and parabolic for $i \geq l_0+1$.
Note that each parabolic nodal component is attached to a unique hyperbolic nodal component via a nodal singularity.
These parabolic nodal components do not correspond to any components in the domain of discontinuity.
Thus, the same construction gives  a connected pinched core polygon
$$
\overline{\cP_1} \cup \overline{\cP_2} \cup \ldots \cup \overline{\cP_{l_0}},
$$
where $\cP_i \subseteq \Delta_i$ is the core polygon associated to $\Sigma_i$.
The level $0$ puzzle pieces $\mathfrak{P}^0 = \{P_1,\ldots, P_r\}$ are constructed exactly the same as the closure of components $\widehat{\C} - \left( \overline{\cP_1} \cup \overline{\cP_2} \cup \ldots \cup \overline{\cP_{l_0}} \cup \overline{\cP_\infty}\right)$.

\subsubsection{The construction of Basilica Bowen--Series maps}\label{subsec:constructionBBS}
In this section, we explain how the natural puzzle constructed induces a Markov $(\Lambda, G)$-map in the sense of \cite{LMM26b} (the details are not important for this paper, Proposition~\ref{prop:MarkovProp} is all that we shall need).

Let $\mathfrak{P}^0 = \{P_1,\ldots, P_r\}$ be the level $0$ puzzle associated to a pinched core polygon $\{\cP_1,\ldots, \cP_{j_0}\}$.
For each $\cP_i$, we choose a subset $\mathcal{PV}_i$ of the vertex set of $\cP_i$ so that it contains exactly one element in the $G_{\Delta_i}$-orbit of each vertex of $\cP_i$, where $G_{\Delta_i}$ is the stabilizer of $\Delta_i$.
We call $\mathcal{PV}_i$ the {\em representative vertex set}.

We are going to define a collection of maps
$$
\mathcal{F} = \{f_i: i\in\{1,\ldots, r\}\}
$$
where $f_i : P_i \longrightarrow \widehat{\C}$.

Suppose that $P_i$ is of limb type. Let $a$ be the unique point in $\partial P_i \cap \Lambda(G)$, and $\Delta \subseteq P_i$ be the component of $\Omega(G)-\Delta_\infty$ with $a \in \partial \Delta$.
Then $\Delta$ corresponds to some hyperbolic nodal component, say $\Sigma_s$.
Let $\Delta_s$ be the component of $\Omega(G)$ associated to $\Sigma_s$ that contains $\cP_s$ in the pinched core polygon.
Let $g \in G$ be an element so that $g(\Delta) = \Delta_s$ and $g(a) \in \mathcal{PV}_s$.
We define the map $f_i: P_i \longrightarrow \widehat{\C}$ as $f_i(x) = gx$. 
We remark that the choice of the element $g$ is not unique: any two choices differ by post-composition with an element in the stabilizer of $g(a)$, which is isomorphic to $\Z$. Note that this group is generated by the (primitive) accidental parabolic element with its fixed point at $g(a)$.

Suppose that $P_i$ is torso type. Let $s \in \{1,\ldots, j_0\}$ be the unique index such that $\Int P_i \cap \Delta_s \neq \emptyset$.
Note that $\cP_s$ is the corresponding core polygon for $\Sigma_s$ as defined in \S \ref{ss:idealcorpolygon}, and hence $P_i\cap \partial \Delta_s$ is some Markov piece $[a,b]$ for the Bowen--Series Markov map on $\Delta_s$ associated to $\cP_s$ as defined in \S \ref{ss:idealcorpolygon}. This Bowen--Series map is defined by some $h \in G_{\Delta_s} \subseteq G$ on $[a,b]$. We define the map $f_i: P_i \longrightarrow \widehat{\C}$ as $f_i(x) = hx$.

\begin{defn}
    Let $G$ be a geometrically finite Bers boundary group, with pinched core polygon $\{\cP_1, \ldots, \cP_{j_0}\}$ and representative vertex set $\mathcal{PV}_i$.
    Let $\mathfrak{P}^0 = \{P_1,\ldots, P_r\}$ be the corresponding level $0$ puzzle.
    We call $\mathcal{F} = \{f_1,\ldots, f_r\}$ constructed above a \emph{Basilica Bowen--Series map} for $G$.
\end{defn}

\begin{remark}\label{bbs_non_unique_rem}
    The construction of a Basilica Bowen--Series map is not canonical. There are several steps where choices need to be made.
    \begin{enumerate}
        \item The pinched core polygon $\{\cP_1,\ldots, \cP_{j_0}\}$ requires a choice. Note that the level $0$ puzzle depends only on such a choice.
        \item The subset $\mathcal{PV}_i$ of the vertex set of $\cP_i$ requires choices of representatives from the $G_{\Delta_i}$-orbit of the vertex set.
        \item Finally, the map on a limb type puzzle piece requires a choice of the element in the stabilizer of the corresponding point.
    \end{enumerate}
\end{remark}

The following proposition follows from \cite[Proposition 5.6 and Lemma 5.8] {LMM26b}.
\begin{prop}\label{prop:MarkovProp}
    A Basilica Bowen--Series map $\mathcal{F} = \{f_1,\ldots, f_r\}$ induces a topologically expanding Markov $(\Lambda,G)$-map on the limit set $\Lambda(G)$.
    More precisely, for each $i$, there exists a subset $I_i \subseteq \{1,\ldots, r\}$ so that
    $$
    f_i(P_i \cap \Lambda(G)) = \bigcup_{j\in I_i} P_j \cap \Lambda(G).
    $$
    Further, $\bigcup_{j\in I_i} P_j \subseteq f_i(P_i)$.
\end{prop}

By Proposition~\ref{prop:MarkovProp}, we inductively define level $n$ puzzle pieces 
$$
\mathfrak{P}^n = \{P_w\}_{w\in \mathcal{I}^n}
$$ 
by pulling back the level $n-1$ puzzle piece.

\subsubsection{The pull-back Bowen--Series maps}\label{pullback_bs_subsec}
Let $\Phi_\infty: \D \longrightarrow \Delta_\infty$ be a Riemann uniformization.
Let $\mathcal{F}|_{\Delta_\infty}$ be the restriction of the Basilica Bowen--Series map on $\Delta_\infty$ (recall that strictly speaking, $\cF$ is only defined on a subset of $\Delta_\infty$), and let $\mathcal{F}_\infty = \{f_{1,\infty},\ldots, f_{r_\infty, \infty}\}$ be the conjugate of $\mathcal{F}|_{\Delta_\infty}$ via $\Phi_\infty$.
By the Schwarz reflection principle, $\mathcal{F}_\infty$ extends to a conformal Markov map on $\mathbb{S}^1 = \partial \D$ (see \S \ref{sec:cmm}).
We call $\mathcal{F}_\infty$ the {\em pull-back Bowen--Series map}.

The following proposition is proved in \cite[Proposition 5.9]{LMM26b}.
\begin{prop}\label{prop:basilicabssymmetric}
    Let $\mathcal{F}$ be a Basilica Bowen--Series map for a geometrically finite Bers boundary group, and let $\mathcal{F}_\infty$ be the pull-back Bowen--Series map.
    Then $\mathcal{F}_\infty$ is topologically expanding, primitive, conformal, and admits a puzzle structure. 
    
    Moreover, all break-points are symmetric. More precisely, a break-point is symmetrically parabolic if it corresponds to a cusp of the surface $\Sigma$; and it is symmetrically hyperbolic otherwise.
    In particular, if $\Sigma$ is compact, then any break-point for $\mathcal{F}_\infty$ is symmetrically hyperbolic.
\end{prop}

\section{Orbit equivalence}\label{sec:orbit_equiv}
In this section, we investigate the orbit equivalence between lattices (both uniform or non-uniform) and certain associated pull-back Bowen--Series maps on the circle (see \S \ref{pullback_bs_subsec}). Recall the notion of regular points for conformal Markov maps from \S \ref{sec:cmm}.

\begin{defn}\label{oe_def}
Let $G<\operatorname{PSL}_2(\mathbb R)$ be a lattice and let $\mathcal{F}$ be a conformal Markov map on $\SS^1$. We say that $\mathcal{F}$ is orbit equivalent to $G$ at a regular point $x\in \SS^1$ if
$$
   \mathcal O_G(x)=\mathcal{GO}_\mathcal{F}(x).
$$
We say that $\mathcal{F}$ is \emph{{regular} orbit equivalent} to $G$ if this equality holds for every regular $x\in \bS^1$. More generally, if $E\subset \bS^1$, we say that $\mathcal{F}$ is {regular} orbit equivalent to $G$ on $E$ if the equality holds for every regular $x\in E$.
\end{defn}

We now state the main result of this section, whose proof will be divided into several cases. 

\begin{theorem}\label{thm:orbitequiv}
    Let $G$ be a lattice in $\PSL_2(\R)$. Then we have the following.
    \begin{itemize}
        \item Suppose that $G$ is non-uniform. Then for any core polygon $\mathcal{P}$ of $G$, the Bowen--Series Markov map $\mathcal{F} = \mathcal{F}_{G, \mathcal{P}}$ is {regular} orbit equivalent to $G$.
        \item Suppose that $G$ is uniform with signature $(g,n) \notin \{(0,3), (1,1)\}$. Then there exists a pull-back Bowen--Series Markov map $\mathcal{F}$ which is {regular} orbit equivalent to $G$ away from one $G$-orbit.

        More precisely, there exists a unique $G$-orbit $\mathcal{O}_G(x), x \in \SS^1 \setminus \mathrm{bk}(\cF)$, {and $x_1, x_2 \in \SS^1$} so that 
    $$
    \mathcal{O}_G(x) = \mathcal{GO}_\mathcal{F}(x_1) \cup \mathcal{GO}_\mathcal{F}(x_2).
    $$
    For any other $G$-orbit $\mathcal{O}_G(y) \neq \mathcal{O}_G(x), y \in \SS^1 \setminus \mathrm{bk}(\cF)$, we have
    $$
    \mathcal{O}_G(y) = \mathcal{GO}_\mathcal{F}(y).
    $$
    \end{itemize}
\end{theorem}

\subsection{Orbit inclusion}
We start with the following simple lemma.
\begin{lem}\label{lem:GO_FsubsetO_G}
    Let $\cF$ be either 
    \begin{itemize}
        \item a Bowen--Series Markov map of a non-uniform lattice $G$; or
        \item a pull-back Bowen--Series Markov map of a uniform lattice $G$ with signature $(g,n) \notin \{(0,3), (1,1)\}$.
    \end{itemize}
    Then for each $x \in \SS^1 \setminus \mathrm{bk}(\cF)$, $\mathcal{GO}_\mathcal{F}(x) \subseteq \mathcal{O}_G(x)$.
\end{lem}
\begin{proof}
    This follows directly from the fact that each branch of $\cF$ is an element of $G$.
\end{proof}

To study the opposite inclusion, we introduce the following local version of orbit equivalence for an element $g \in G$.
\begin{defn}
    Let $G$ be a lattice of $\PSL_2(\R)$, and let $\mathcal{F}: \mathbb{S}^1 \longrightarrow \mathbb{S}^1$ be a map (not necessarily continuous).
    We say that an element $h \in G$ is {\em good} (with respect to $\mathcal{F}$) if 
   \[
  x\in \mathbb{S}^1\setminus\mathrm{bk}(\mathcal{F})
  \quad\Longrightarrow\quad
  hx\in\mathbb S^1\setminus\mathrm{bk}(\mathcal{F})
  \ \text{and}\ x\sim_{\mathcal{F}} hx.
    \]
    We say $h$ is {\em good on $X$} or {\em good away from $\mathbb{S}^1 \setminus X$} if
    \[
  x\in X\setminus\mathrm{bk}(\mathcal{F})
  \quad\Longrightarrow\quad
  hx\in\mathbb S^1\setminus\mathrm{bk}(\mathcal{F})
  \ \text{and}\ x\sim_{\mathcal{F}} hx.
    \]
\end{defn}
The following property follows from the definition and the symmetric nature of the grand-orbit equivalence relation. 
\begin{lem}\label{lem:goodprop}
    If $h$ is good, then $h^{-1}$ is good.
    If $h$ is good away from $I$, then $h^{-1}$ is good away from $h(I)$.
\end{lem}

\subsection{Non-uniform lattice}
\begin{prop}\label{prop:oenonuniform}
    Let $G$ be a non-uniform lattice, and $\mathcal{P}$ be a core polygon for $G$. 
    Then the Bowen--Series Markov map $\mathcal{F} = \mathcal{F}_{G, \mathcal{P}}$ is {regular} orbit equivalent to $G$ on $\SS^1$.
\end{prop}

\begin{proof}
    Suppose that $G$ is a non-uniform lattice associated to a surface of type $(g, n; \nu_1,\ldots, \nu_n)$ with $k \geq 1$ punctures.
    Let $\mathcal{P}$ be a core polygon for $G$ as in \S \ref{ss:idealcorpolygon} (see Figure~\ref{fig:cuspBM}). Then 
    $$
    S = \{A_i, B_i, 0 \leq i \leq g-1\} \cup \{P_i: 0 \leq i < k-1\} \cup \{R_i: 0 \leq i < n-k\}
    $$ 
    generates the group $G$.

    By Lemma~\ref{lem:GO_FsubsetO_G}, $\mathcal{GO}_\mathcal{F}(x) \subseteq \mathcal{O}_G(x)$. 
    To prove the other direction, we first claim that each $l \in S$ is good.
    
    \begin{proof}[Proof of the claim]
    Let $l \in S$, and let $x \in \mathbb{S}^1 \setminus \mathrm{bk}(\cF)$ and $y = lx$. 
    
    Suppose $l \in  \{A_i, B_i, 0 \leq i \leq g-1\} \cup \{P_i: 0 \leq i < k-1\}$. Then there exist a pair of intervals $I, J$ so that $\mathcal{F}(t) = lt$ for $t\in I$ and $\mathcal{F}(t) = l^{-1}t$ for $t \in J$. Moreover, $\mathcal{F}(J) = \mathbb{S}^1 - \Int(I)$. Thus, if $x \in I$, we have $\mathcal{F}(x) = y$. If $x\notin I$, we have $\mathcal{F}(y) = x$. In both cases, $y \in \mathbb{S}^1 \setminus \mathrm{bk}(\cF)$ and $y \in \mathcal{GO}_\mathcal{F}(x)$.

    Suppose that $l \in \{R_i: 0 \leq i < n-k\}$. Then there exists a partition $\mathbb{S}^1 = I_0\cup I_1\ldots\cup I_{s-1}$ with $s=\nu_i$, $l(I_i) = I_{i+1}$ (where the indices are considered mod $s$) so that $\mathcal{F}(t) = l^{s-j}t$ if $t \in I_j, j \neq 0$. Thus, if $x \notin I_0$, then $\mathcal{F}(x) = \mathcal{F}(y)$. Otherwise, $x \in I_0, y \in I_1$, and $\mathcal{F}(y) = x$. In both cases, $y \in \mathbb{S}^1 \setminus \mathrm{bk}(\cF)$ and $y \in \mathcal{GO}_\mathcal{F}(x)$.
    \end{proof}
    Write $g = g_rg_{r-1}\ldots g_1$ where $g_i \in S \cup S^{-1}$. By the claim, we have $g_k\ldots g_1x \in \mathcal{GO}_{\mathcal{F}}(g_{k-1}\ldots g_1x)$. Therefore, $gx \in \mathcal{GO}_{\mathcal{F}}(x)$. This shows that $\mathcal{O}_G(x) \subseteq \mathcal{GO}_\mathcal{F}(x)$, and Proposition~\ref{prop:oenonuniform} follows. 
\end{proof}

\begin{remark}
As mentioned earlier, the Bowen--Series Markov map $\cF$ is a covering map of the circle when the surface $\Sigma$ is a punctured sphere, possibly with an orbifold point of order two. In this case, the grand orbits of $\cF$ are well-defined at the break-points as well, and the map $\cF$ is orbit equivalent to the corresponding Fuchsian lattice on the whole circle, including the break-points (cf. \cite{MM1}).
\end{remark}

\subsection{Uniform lattice without torsion}
In this subsection, we prove an orbit equivalence theorem for closed surface groups. 

Recall the notion of pull-back Bowen--Series maps from \S \ref{pullback_bs_subsec}. 
\begin{prop}\label{prop:closedsurfaceOE}
    Let $G \leq \PSL_2(\R)$ be a uniform lattice such that $\D/G$ is a closed genus $g\geq 2$ surface
    (with no orbifold points). There is a pull-back Bowen--Series Markov map $\mathcal{F}:\mathbb{S}^1 \longrightarrow \mathbb{S}^1$ and a unique exceptional $G$-orbit \[
  E\subset\mathbb S^1\setminus\mathrm{bk}(\mathcal{F})
\]
with the following properties.
\begin{enumerate}[leftmargin=8mm]
\item If $x\in\mathbb S^1\setminus(\mathrm{bk}(\mathcal{F})\cup E)$, then
\[
  \mathcal O_G(x)=\operatorname{GO}_\mathcal{F}(x).
\]
\item There are points $x_v,x_s\in E$ in distinct $\mathcal{F}$-grand orbits such
that
\[
  E=\operatorname{GO}_\mathcal{F}(x_v)\cup\operatorname{GO}_\mathcal{F}(x_s).
\]
\end{enumerate}
In particular, $\mathcal{F}$ is {regular} orbit equivalent to $G$ on $\mathbb S^1$ away from the single $G$-orbit $E$.
\end{prop}
 \begin{figure}[htp]
 \captionsetup{width=0.96\linewidth}
    \centering
    \includegraphics[width=0.75\textwidth]{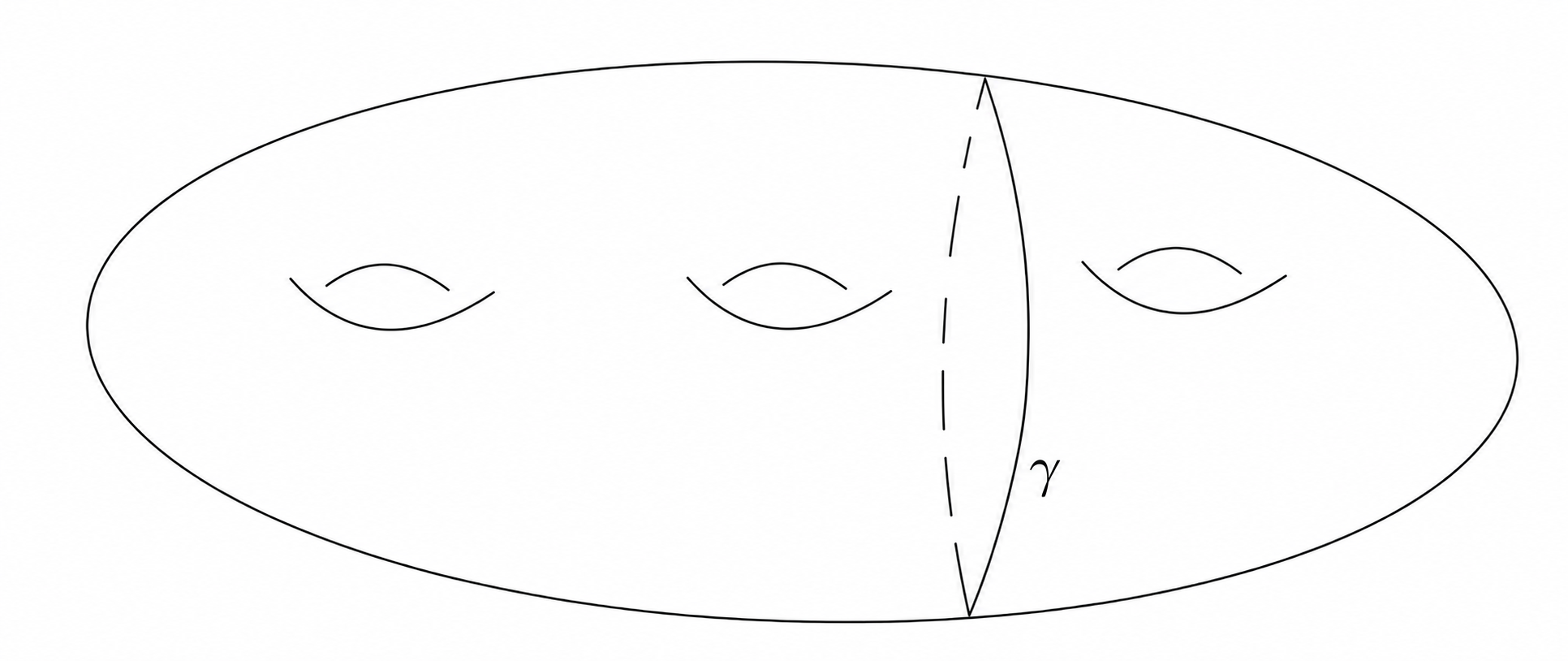}
    \includegraphics[width=1\textwidth]{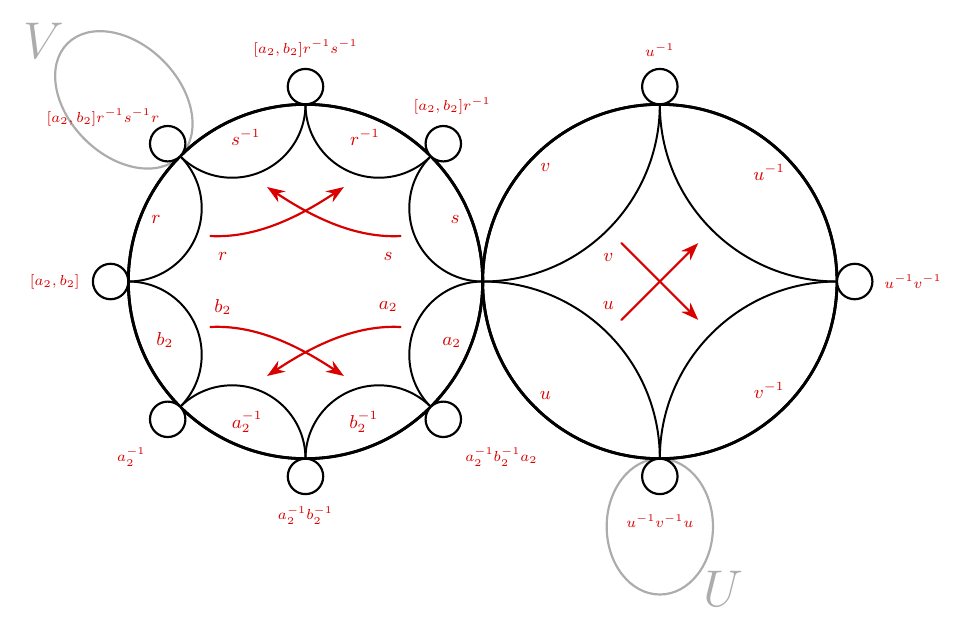}
    \caption{The Basilica Bowen--Series map for $\widetilde{G}$ on the boundary of the Bers slice of type $(3,0)$.}
    \label{fig:g3}
\end{figure}

\subsubsection{The pinching construction and notation}\label{sss:closed-pinching-notation}
Choose a separating simple closed curve $\gamma$ which cuts off a
one-holed torus.  Pinching $\gamma$ gives a geometrically finite boundary
group $\widetilde{G}$ with two hyperbolic nodal components; see Figure~\ref{fig:g3} (alternately, this is a geometrically finite Kleinian surface group with one accidental parabolic).  Let
$\Delta_\infty \subseteq \Omega(\widetilde{G})$ be the invariant component.  We choose the marking so that the
action of $\widetilde{G}$ on $\Delta_\infty$ is conformally conjugate to the
action of $G$ on $\mathbb D$.  Construct a Basilica Bowen--Series map
$\widetilde{\mathcal{F}}$ for this pinched surface, and let $\mathcal{F}$ be its pull-back to
$\mathbb S^1=\partial\mathbb D$. We remark that this is the simplest possible choice in the definition of $\widetilde{\cF}$; see Remark~\ref{bbs_non_unique_rem}.

Let $\Phi\colon\mathbb D\to\Delta_\infty$ be a Riemann map.  As in the
definition of the pull-back Bowen--Series map, $\Phi$ induces a
$G$-equivariant identification of ideal
boundaries
\[
  \partial^I\Phi\colon\mathbb S^1\longrightarrow
  \partial^I\Delta_\infty.
\]
By definition, a torso (respectively, limb)
interval of $\mathbb S^1$ is the pull-back, under this ideal-boundary identification, of the prime end represented by a torso (respectively, limb) puzzle piece.

Choose standard generators
\(a_1,b_1,\ldots,a_g,b_g
\) for $G$.  
Abusing  notation slightly, we  use the same notation for  the generators for $\widetilde{G}$ under the natural isomorphism between $G$ and $\widetilde{G}$.

Put
\[
  u:=a_1,\qquad v:=b_1,\qquad r:=a_g,\qquad s:=b_g,
\]
and the orientation of the side-pairings is chosen so that
\[
  K:=[a_{g-1},b_{g-1}]\cdots[a_2,b_2],
  \hspace{1mm}
  W:=[u,v]=K[r,s],
  \hspace{1mm} \text{where } [p,q]:=p^{-1}q^{-1}pq.
\]
Here $K=\mathrm{id}$ when $g=2$, and $W$ is a primitive generator of the stabilizer of the central node. The pinched core polygon of $\widetilde{G}$ consists of two core polygons which touch at an (accidental) parabolic fixed point $\xi$. The pull-back of $\xi$ on the ideal boundary $\partial^I\Delta_\infty$ consists of the two fixed points $\xi_\pm$ of $W$.
See Figure~\ref{fig:g3} for an illustration for $g = 3$.

We partition the ideal boundary of $\Delta_\infty$ into two closed intervals $I_R, I_L$ as follows. Let $I_R\subset\partial^I\Delta_\infty$ be the closed interval 
corresponding to the component of $\partial\Delta_\infty\setminus\{\xi\}$ containing the one-holed torus core
polygon, and let $I_L$ be the interval corresponding to the other
component. Evidently, the common endpoints of $I_L,I_R$ are the two fixed points $\xi_\pm$ of $W$. Let $U, V$ be the limbs at ideal vertices of the one-holed torus polygon and the other core polygon which are mapped to the central node by $u^{-1}v^{-1}u$ and $K r^{-1}s^{-1}r$ respectively; see Figure~\ref{fig:g3}.  Denote the
corresponding intervals on the ideal boundary $\partial\Delta_\infty$ by $I_U$ and $I_V$.  Thus
\[
  I_U\subset I_R,\qquad I_V\subset I_L,
  \qquad v(I_L)=I_U,\qquad s(I_R)=I_V.
\]
We first record how the break-points of $\cF$ interact with the
group action. 

\begin{lem}\label{lem:closed-regular-set-invariant}
Let $\xi_\pm$ be as above. Then
\[
  \mathrm{bk}(\mathcal{F})
  =\mathcal O_G(\xi_+)\cup\mathcal O_G(\xi_-).
\]
Consequently, both $\mathrm{bk}(\mathcal{F})$ and its complement are
$G$-invariant.
\end{lem}
\begin{proof}
    By construction, $\mathrm{bk}(\mathcal{F})$ corresponds to the union, over all cut-points $x\in\Lambda(\widetilde{G})$, of the two prime ends of $\Delta_\infty$ associated with $x$. Since $\widetilde{G}$ is obtained by pinching a single simple closed curve, all cut-points of $\Lambda(\widetilde{G})$ lie in a single $\widetilde{G}$-orbit. The lemma follows.
\end{proof}

\subsubsection{The good generators}\label{sss:mostgood}
\begin{lem}\label{lem:closed-most-generators-good}
The generators
\[
  u,r,\quad\text{and}\quad a_j,b_j,\ j\in\{2,\cdots, g-1\},
\]
are good. The generator $v$ is good away from $I_L$, and $s$ is good
away from~$I_R$.  

Consequently, $u^{-1},r^{-1}$ and $a_j^{-1},b_j^{-1}$, $j\in\{2,\cdots, g-1\}$, are good, and $v^{-1}$ is good away from $I_U$ and $s^{-1}$ is good away from $I_V$.
\end{lem}
\begin{proof} We start by assigning  \emph{exceptional sets} $E_h$ to generators $h$ as in the following table:

\[
\begin{array}{c|c|c}
 h & E_h & h(E_h)\\ \hline
 u,r,a_j,b_j\ (2\leq j\leq g-1) & \varnothing & \varnothing\\
 v & I_L & I_U\\
 s & I_R & I_V.
\end{array}
\]
For a generator $h$ and a regular point $x$ outside the exceptional set $E_h$, we claim that one of the following holds:
\[
  \mathcal{F}(x)=hx,\qquad \mathcal{F}(hx)=x,\qquad\text{or}\qquad \mathcal{F}(x)=\mathcal{F}(hx).
\]

Suppose first that $x$ projects to a torso-type puzzle piece on which $\widetilde{\mathcal{F}}=h$.
Then the first  identity holds. If
 $hx$ projects to a torso-type puzzle piece on which $\widetilde{\mathcal{F}}=h^{-1}$, then the  second identity holds.

Next suppose that $x, hx \in I_L$ or $x, hx \in I_R$. Then  $x, hx$ project to two limb-type puzzle pieces whose forward branches have the same image; this is the third identity.

The only remaining possibilities are
\[
h\in\{u,v\}\ \text{and}\ x\in I_L,
\qquad\text{or}\qquad
h\in\{a_2,s\}\ \text{and}\ x\in I_R.
\]
Since $x\notin E_h$, the cases $h = v$ or $h = s$ are excluded. Thus either $h = u, x\in I_L$ or $h = a_2, x\in I_R$. In either case, $hx$ projects to a limb-type puzzle piece where $\widetilde{\mathcal{F}} = h^{-1}$, and thus the second identity holds.
This proves the claim.

Since $x$ is a regular point, each of the three relations implies that $hx \in \mathbb{S}^1 \setminus \mathrm{bk}(\cF)$ and $x\sim_{\mathcal{F}} hx$. Thus, the lemma follows for the generators. The analogous statements for the inverses of the generators follow from Lemma~\ref{lem:goodprop}.
\end{proof}

Our strategy now is to augment the domain where $v, s$ are good, culminating in Lemma~\ref{lem:closed-exceptional-points} below. We need to pass to higher iterates of $\mathcal{F}$ to implement this strategy. It is here that the commutator $W = [u,v]$ comes into play via the following straightforward preliminary computation.
\begin{lem}\label{lem:closed-node-branches}
If $x\in I_L\setminus\mathrm{bk}(\mathcal{F})$, then
\[
  \mathcal{F}(vx)=u^{-1}v^{-1}u(vx)=Wx\in I_L.
\]
If $x\in I_R\setminus\mathrm{bk}(\mathcal{F})$, then
\[
  \mathcal{F}(sx)=K r^{-1}s^{-1}r(sx)=Wx\in I_R.
\]
\end{lem}

\subsubsection{The two exceptional transitions}\label{subsubsec-except}
Define
\[
  \sigma:=v^{-1}uv=uW,
  \quad
  \tau:=s^{-1}rs=rK^{-1}W,
  \quad
  m:=4(g-2)+1=4g-7.
\]
The above relations for $\sigma$ and $\tau$ follow immediately from the definition\\
$W=[u,v]=K[r,s]$. We also recall that
$$
x\in I_L\iff Wx\in I_L,\quad \mathrm{and}\quad x\in I_R\iff Wx\in I_R.
$$

\begin{lem}[Exceptional transitions]\label{lem:closed-defect-transitions} $ $
\begin{enumerate}[leftmargin=8mm]
\item Suppose that $x\notin\mathrm{bk}(\mathcal{F})$ and $x,Wx\in I_R$.
Then either 
\begin{itemize}[leftmargin=6mm]
    \item $Wx\sim_{\mathcal{F}} x$, or
    \item $uvx=\mathcal{F}^2(x)\in I_U$.
\end{itemize}
In the latter case,
\[
  \mathcal{F}^3(x)=W\sigma x,\qquad \mathcal{F}(Wx)=\sigma x,
\]
and both points belong to $I_L$.

\item Suppose that $x\notin\mathrm{bk}(\mathcal{F})$ and $x,Wx\in I_L$.
Then either 
\begin{itemize}
    \item $Wx\sim_{\mathcal{F}} x$, or
    \item $rsx=\mathcal{F}^2(x)\in I_V$.
\end{itemize}
In the latter case,
\[
  \mathcal{F}^3(x)=W\tau x,\qquad \mathcal{F}^m(Wx)=\tau x,
\]
and both points belong to $I_R$.
\end{enumerate}
\end{lem}
\begin{proof}
For the first part, write
\[
  Wx=u^{-1}v^{-1}uvx.
\]
Since $x$ is regular, it is not a common endpoint of $I_L$ and $I_R$.
Thus the generator $v$ is good at $x$.
The generator $u$ is good at every regular point.  If $uvx\notin I_U$, then $v^{-1}$ is also good at $uvx$ (by Lemma~\ref{lem:closed-most-generators-good}), and hence $u^{-1}$ is good at $v^{-1}uvx$.  Thus $Wx\sim_{\mathcal{F}} x$.
Therefore failure of the relation $Wx\sim_{\mathcal{F}} x$ can occur only when $uvx\in I_U$. 
Note that if $uvx\in I_U$, then $x$ lies in $v^{-1} u^{-1}I_U$ which is contained in the torso interval on which $\cF$ acts as $v$. Further, on $v^{-1} u^{-1}I_U$, 
the successive iterates of $\mathcal{F}$ are given by the following. 
\[
  \mathcal{F}(x)=vx,\qquad \mathcal{F}^2(x)=uvx,\qquad
  \mathcal{F}^3(x)=u^{-1}v^{-1}u(uvx)=W\sigma x \in I_L.
\]
Observe also that if $uvx\in I_U$, then $Wx \in u^{-1}v^{-1}I_U$ and the map $\cF$  on $u^{-1}v^{-1}I_U$ equals $u$, and hence
\[
  \mathcal{F}(Wx)=uWx=\sigma x \in I_L.
\]

For the second part, use
\[
  Wx=K r^{-1}s^{-1}rsx.
\]
Again, regularity of $x$ ensures that it is not a  common endpoint of $I_L$ and $I_R$.  Hence $s$ is good at $x$, and $r$ is good
everywhere.  If
$rsx\notin I_V$, then $s^{-1}$ is good at $rsx$ (by Lemma~\ref{lem:closed-most-generators-good}); and then $r^{-1}$ and
every letter in $K$ is good at the subsequent images of $x$. Thus $Wx\sim_{\mathcal{F}} x$ unless
$rsx\in I_V$. Note that $s^{-1} r^{-1}I_V$ is contained in a torso interval on which $\cF$ acts as $s$. Further, on $s^{-1}r^{-1}I_V$, we have,
\[
  \mathcal{F}(x)=sx,\qquad \mathcal{F}^2(x)=rsx,\qquad
  \mathcal{F}^3(x)=K r^{-1}s^{-1}r(rsx)=W\tau x\in I_R.
\]
It remains to prove the statement for the iterates  of $Wx$ under $\cF$ when $rsx\in I_V$. In this case, $Wx \in Kr^{-1}s^{-1}I_V$, and we have $\cF^{4(g-2)} (Wx) = K^{-1}Wx$. Also, $K^{-1}Wx \in r^{-1}s^{-1}I_V$, so the next iterate of $\cF$ is given by $r$. Hence
\[
  \mathcal{F}^m(Wx)=rK^{-1}Wx=\tau x \in I_R.
\]
\end{proof}

\subsubsection{The two exceptional periodic points}
Recall the definitions of $\sigma, \tau$ from \S \ref{subsubsec-except}. Let $H_v:=\sigma\tau$ and $H_s:=\tau\sigma$, and let $x_v, x_s$ be the repelling fixed points of $H_v, H_s$ respectively. Evidently, $x_v,x_s\notin\{\xi_+,\xi_-\}$; i.e., they are not the common endpoints of $I_R$ and $I_L$. Note also that $x_s=\tau x_v$ and $x_v=\sigma x_s$; in particular, $\mathcal O_G(x_v)=\mathcal O_G(x_s)$.

\begin{lem}\label{lem:closed-exceptional-points}

For every regular $x$,
\[
  vx\sim_{\mathcal{F}} x\quad\Longleftrightarrow\quad x\neq x_v,\quad \mathrm{and} 
  \quad
  sx\sim_{\mathcal{F}} x\quad\Longleftrightarrow\quad x\neq x_s.
\]
The periods of $x_v,x_s$ under $\mathcal{F}$ are
\[
  \operatorname{per}_{\mathcal{F}}(x_v)=4,
  \qquad
  \operatorname{per}_{\mathcal{F}}(x_s)=m+3=4(g-1),
\]
and their $\mathcal{F}$-grand orbits are distinct.  Finally,
\[
  \mathcal{F}^{m+1}(vx_v)=x_s,
  \qquad
  \mathcal{F}^2(sx_s)=x_v.
\]
\end{lem}
\begin{proof}
Let us first assume that $vx\not\sim_{\mathcal{F}} x$. Then, Lemma \ref{lem:closed-most-generators-good} gives $x\in I_L$, and Lemma \ref{lem:closed-node-branches} gives
\begin{equation}
  \mathcal{F}(vx)=Wx \in I_L.
 \label{vxv_f_image_eqn} 
\end{equation}
Thus $x\not\sim_{\mathcal{F}} Wx$.  
Applying part (2) of Lemma \ref{lem:closed-defect-transitions} to $x\in I_L$ and then part (1) of the same lemma to $\tau x\in I_R$ shows that the $\cF$-orbits of $x$ and $Wx$ are forced to return to $I_L$ eventually, and their first returns to $I_L$ are given by
\begin{equation}
  \mathcal{F}^4(x)=\sigma\tau x=H_vx, \quad \mathcal{F}^{m+3}(Wx) = W\sigma\tau x = W H_vx.
  \label{period_4_eqn}
\end{equation}
By transitivity of the grand orbit equivalence relation, $\mathcal{F}^4(x) \not\sim_\cF \mathcal{F}^{m+3}(Wx)$, so $H_vx \not\sim_\cF WH_v x$.
Thus, we can apply the above argument to $H_vx$ and $WH_vx$, and conclude that 
$$
\cF^8(x) = \cF^4(H_vx) = H_v^2x \not\sim_\cF \cF^{2(m+3)}(Wx) = \mathcal{F}^{m+3}(WH_vx) = WH_v^2x.
$$
Inductively, we conclude that $\mathcal{F}^{4n}(x) = H_v^nx \in I_L$ for all $n \geq 0$.

Therefore, $x$ must be the repelling fixed point $x_v$ of $H_v$. It has minimal period $4$ under $\cF$ by Relation~\eqref{period_4_eqn}.

We now look at the repelling fixed point $x_s = \tau x_v$ of $H_s$.
A similar computation as above shows that $sx\not\sim_{\mathcal{F}} x \implies x =x_s$, and that the minimal period of $x_s$ is $m+3$.

Now we show that the $\cF$-grand orbits of $x_v, x_s$ are distinct. Indeed, if $g \geq 3$, then $m+3 = 4(g-1) \geq 8$. Since $x_v, x_s$ are periodic with different periods, they lie in distinct cycles, and hence do not belong to the same grand orbit under $\cF$.
If $g = 2$, then a direct computation shows that the periodic orbit $\{\cF^n(x_v)\}_{n= 0}^3$ does not contain the periodic point $x_s$, so they are not in the same grand orbit.

It remains to prove that $vx_v \not\sim_\cF x_v$ and $sx_s \not\sim_\cF x_s$. It is readily verified from the construction that $x_v\in s^{-1} r^{-1}I_V\subset I_L$ and $x_s\in v^{-1} u^{-1}I_U\subset I_R$. It now follows from Lemma~\ref{lem:closed-node-branches} and Lemma~\ref{lem:closed-defect-transitions} that
\[
  \mathcal{F}(vx_v)=Wx_v,\qquad \mathcal{F}^m(Wx_v)=\tau x_v=x_s,
\]
and hence $\mathcal{F}^{m+1}(vx_v)=x_s$.  Likewise,
\[
  \mathcal{F}(sx_s)=Wx_s,\qquad \mathcal{F}(Wx_s)=\sigma x_s=x_v,
\]
and hence $\mathcal{F}^2(sx_s)=x_v$.
Thus, the $\cF$-orbit of $vx_v$ eventually lands in the cycle containing $x_s$, and the $\cF$-orbit of $sx_s$ eventually lands in the cycle containing $x_v$. As the two cycles containing $x_v,x_s$ are disjoint, we conclude that $vx_v\not\sim_{\mathcal{F}} x_v$ and $sx_s\not\sim_{\mathcal{F}} x_s$.
\end{proof}

\begin{proof}[Proof of Proposition~\ref{prop:closedsurfaceOE}]
    Let $x\in \mathbb{S}^1 \setminus\mathrm{bk}(\cF)$. 
    By Lemma~\ref{lem:GO_FsubsetO_G}, $\mathcal{GO}_\mathcal{F}(x) \subseteq \mathcal{O}_G(x)$.

    Let $E:=\mathcal O_G(x_v)=\mathcal O_G(x_s)$, and suppose $x\notin E$.
Write
\[
  h=h_k\cdots h_1,
  \qquad
  h_i\in\{a_1^{\pm1},b_1^{\pm1},\ldots,a_g^{\pm1},b_g^{\pm1}\}.
\]
By Lemma~\ref{lem:closed-regular-set-invariant}, $h_j\cdots h_1x$ is regular for each $j = 1,\cdots, k$. Since $E$ is a $G$-orbit, $h_j\cdots h_1x$ also remains outside $E$ for all $j$. Lemmas
\ref{lem:closed-most-generators-good} and
\ref{lem:closed-exceptional-points}, together with Lemma~\ref{lem:goodprop}, show inductively that
\[
  h_j\cdots h_1x\sim_{\mathcal{F}} h_{j-1}\cdots h_1x.
\]
Thus $hx\sim_{\mathcal{F}} x$.  This proves
$\mathcal O_G(x)\subseteq\operatorname{GO}_{\mathcal{F}}(x)$ and hence item (1) of Proposition~\ref{prop:closedsurfaceOE} follows.

Finally, we claim that 
    $$
    \mathcal{O}_G(x_v) = \mathcal{GO}_\mathcal{F}(x_v) \cup \mathcal{GO}_\mathcal{F}(x_s).
    $$
    Every grand orbit of $\cF$ is preserved by any signed generator other than $v^{\pm 1}, s^{\pm 1}$. 
    By Lemma~\ref{lem:closed-exceptional-points}, the signed generators $v^{\pm 1}, s^{\pm 1}$ switch the two grand orbits at $x_v, x_s$ respectively. Every other grand orbit is preserved by  $v^{\pm 1}, s^{\pm 1}$. Therefore, the right hand side is invariant under $G$. The claim and thus the proposition follow.
\end{proof}

\subsection{Uniform lattice with torsion}
In this subsection, we state the result for general uniform lattices (possibly with torsion). The proof is similar to Proposition~\ref{prop:closedsurfaceOE}, but involves a detailed case-by-case analysis. The proof is postponed to Appendix~\ref{sec:oeorbifold}.

\begin{prop}\label{prop:nonclosedsurfaceOESt}
    Let $G$ be a uniform lattice corresponding to a compact hyperbolic surface of type $(g, n; \nu_1,\ldots,\nu_n)$ where $n \geq 1$ and 
    $$
    (g,n) \notin\{(0,3), (1,1)\}.
    $$
    There is a pull-back Bowen--Series Markov map $\mathcal{F}:\mathbb{S}^1 \longrightarrow \mathbb{S}^1$ and a unique exceptional $G$-orbit \[
  E\subset\mathbb S^1\setminus\mathrm{bk}(\mathcal{F})
\]
with the following properties.
\begin{enumerate}
\item If $x\in\mathbb S^1\setminus(\mathrm{bk}(\mathcal{F})\cup E)$, then
\[
  \mathcal O_G(x)=\operatorname{GO}_\mathcal{F}(x).
\]
\item There are points $x_v,x_s\in E$ in distinct $\mathcal{F}$-grand orbits such
that
\[
  E=\operatorname{GO}_\mathcal{F}(x_v)\cup\operatorname{GO}_\mathcal{F}(x_s).
\]
\end{enumerate}
In particular, $\mathcal{F}$ is regular orbit equivalent to $G$ on $\mathbb S^1$ away from the single $G$-orbit $E$.
\end{prop}

\subsection{Proof of Theorem~\ref{thm:orbitequiv}}
Theorem~\ref{thm:orbitequiv} follows immediately  from the cases covered in Proposition~\ref{prop:oenonuniform}, Proposition~\ref{prop:closedsurfaceOE} and Proposition~\ref{prop:nonclosedsurfaceOESt}. \qed

\section{Acceleration of pull-back Bowen--Series maps}\label{acceleration_sec}
{In this section, we construct covering maps of the circle, called \emph{accelerated Bowen--Series maps}, by applying suitable accelerations to 
\begin{enumerate}[leftmargin=8mm]
    \item[(i)] the Bowen--Series maps constructed in \S \ref{ss:idealcorpolygon} (for non-uniform lattices), and
    \item[(ii)] the pullback Bowen--Series maps constructed in \S \ref{ss:bbspbs} (for uniform lattices).
\end{enumerate}
The main results are Theorem~\ref{thm:accnon-uniformlattice} for non-uniform lattices and Theorem~\ref{thm:accuniformlattice} for uniform lattices. 

 In order to study orbit equivalence properties of these (accelerated) continuous piecewise M\"obius maps, we need to upgrade the notion of regular orbit equivalence (see Definition~\ref{oe_def}). If $\cG$ is a continuous conformal Markov map on $\bS^1$, the grand orbit of each $x\in\bS^1$ (including the break-points) under $\cG$ is well-defined. 
\begin{defn}\label{defn:strongOEweakOE}
We say that a continuous conformal Markov map $\cG$ is \emph{strongly orbit equivalent} (respectively, \emph{weakly orbit equivalent}) to a Fuchsian lattice $G$ (on $\bS^1$) if $\cG$ is orbit equivalent to $G$ at each point of $\bS^1$ (respectively, on a full-measure subset of $\bS^1$).
\end{defn}

In Theorems~\ref{strong_oe_thm} and~\ref{thm:weakorbitequiv}, we record the orbit equivalence properties of accelerated Bowen--Series maps (also see Proposition~\ref{countable_presentation_prop}).

\begin{remark}
We emphasize that while a classical Bowen--Series map is regular orbit equivalent to the corresponding Fuchsian lattice, it is typically discontinuous. The combination of continuity and weak orbit equivalence along with its piecewise M\"obius nature
is what makes accelerated Bowen--Series maps valuable for the mating constructions carried out in the subsequent sections. In particular, for closed surfaces, they are quite different from the classical Bowen--Series maps of \cite{BS79}.
\end{remark}}

\subsection{Acceleration of a Markov map}\label{acc_subsec}

\begin{defn}\label{acc_def}
    Let $\mathcal{F}$ be a topologically expanding conformal Markov map on $\mathbb{S}^1$, with respect to a partition $\mathcal{A}$. We say that a Markov map $\mathcal{G}$ with respect to a partition $\mathcal{B}$ is an {\em acceleration} of $\mathcal{F}$ if $\mathcal{B}$ is a refinement of $\mathcal{A}$, and for each $B\in \mathcal{B}$, there exists $n = n(B)$ so that 
    $$
    \mathcal{G}|_{B} = \mathcal{F}^n|_{B}.
    $$
    We set $\displaystyle \cN = \max_{B\in \mathcal{B}} n(B)$, and call it the {\em rate} of the acceleration.
\end{defn}

{The rate of acceleration quantitatively controls the failure of orbit equivalence between $\cF$ and $\cG$, as the next result shows. }
Recall that a point $x \in \mathbb{S}^1$ is {\em regular} if its orbit is disjoint from the break-points. 
\begin{prop}\label{prop:virtualorbitequiv}
    Let $\mathcal{F}$ be a topologically expanding conformal Markov map on $\mathbb{S}^1$ and $\mathcal{G}$ an acceleration of $\mathcal{F}$ with rate of acceleration $\cN$.
    Then for any regular point $x\in \mathbb{S}^1$, there exist $k=k(x) \leq \cN$ {and $\xi_1,\ldots, \xi_k\in\bS^1$} such that
    $$
    \mathcal{GO}_\mathcal{F}(x) = \bigcup_{i=1}^k\mathcal{GO}_\mathcal{G}(\xi_i),
    $$
    where $\mathcal{GO}_\mathcal{F}(t)$, $\mathcal{GO}_\mathcal{G}(t)$ denote the grand orbit of $t\in\bS^1$ under $\mathcal{F}$, $\mathcal{G}$, respectively. 
\end{prop}
The assumption that $x$ is a regular point is mild, as only finitely many grand orbits contain break-points.
\begin{proof}

    We first claim that for any regular $x\in \mathbb{S}^1$, $\mathcal{GO}_\mathcal{G}(x) \subseteq \mathcal{GO}_\mathcal{F}(x)$.
    Indeed, let $y \in \mathcal{GO}_\mathcal{G}(x)$. Then there exists $n, m$ so that $\mathcal{G}^n(y) = \mathcal{G}^m(x)$. Since $\mathcal{G}$ is an acceleration of $\mathcal{F}$, there exists $n', m'$ so that $\mathcal{G}^n(y) = \mathcal{F}^{n'}(y)$ and $\mathcal{G}^m(x) = \mathcal{F}^{m'}(x)$.
    Therefore, $y \in \mathcal{GO}_\mathcal{F}(x)$. This proves the claim.

    Let $\cN$ be the rate of the acceleration.
    We now claim that for each regular $x \in \mathbb{S}^1$, there exist $\xi_1,\ldots, \xi_k\in\bS^1$, $k \leq \cN$, so that
    $$
    \mathcal{GO}_\mathcal{F}(x) = \bigcup_{i=1}^k\mathcal{GO}_\mathcal{G}(\xi_i).
    $$
    Note that by the previous claim, each grand orbit of $\mathcal{G}$ is contained in a grand orbit of $\mathcal{F}$.
    By way of contradiction, suppose that $\mathcal{GO}_\mathcal{F}(x)$ contains $\cN+1$ distinct grand orbits of $\mathcal{G}$, and denote them by $\mathcal{GO}_\mathcal{G}(\xi_i)$, $i\in\{1,\ldots, \cN+1\}$.
    Set $x_n = \mathcal{F}^n(x)$.
    Since $\mathcal{GO}_\mathcal{G}(\xi_i) \subseteq \mathcal{GO}_\mathcal{F}(x)$, by replacing $\xi_i$ with {an iterated image of it under $\mathcal{G}$ if necessary}, we may assume that $\xi_i = x_{\sigma(i)}$ for some $\sigma(i) \in \N$.
    Choose $K$ sufficiently large so that $K \geq \sigma(i)$ for all $i\in\{1,\ldots, \cN+1\}$. Since $\cN$ is the rate of the acceleration, {the $\cG$-orbit of $\xi_i$ cannot miss $\cN$ consecutive points of the sequence $(x_n)$, and hence, each $\xi_i$ falls into the set $\{x_K, x_{K+1},\ldots, x_{K+\cN-1}\}$ under some iterate of $\mathcal{G}$.}
    Since there are only $\cN$ elements in this list, two of the points $\xi_i$, $i\in\{1,\ldots, \cN+1\}$, say $\xi_1, \xi_2$, must be mapped to the same point. Thus, $\mathcal{GO}_\mathcal{G}(\xi_1) = \mathcal{GO}_\mathcal{G}(\xi_2)$, which is a contradiction. This proves the claim and hence the proposition.
\end{proof}

{The following general criterion for orbit equivalence between $\cF$ and $\cG$ will be useful in our study of accelerated Bowen--Series maps.}

\begin{prop}\label{prop:oe}
     Let $\mathcal{F}$ be a topologically expanding conformal Markov map on $\mathbb{S}^1$ and $\mathcal{G}$ an acceleration of $\mathcal{F}$.
     Let $\cN$ be the rate of the acceleration.
     Suppose there exists an interval $I$ so that for all {regular} $x\in I$ and for all $k \leq \cN$, $\mathcal{G}^k(x) = \mathcal{F}^k(x)$.
     Let 
     $$
     \mathcal{K}(I) = \{x\in \mathbb{S}^1: |\{n: \cF^n(x)\in I\}| = \infty\}.
     $$
     Then $\mathcal{F}$ and $\mathcal{G}$ are orbit equivalent for any regular point in $\mathcal{K}(I)$.
\end{prop}
\begin{proof}
    The proof is similar to the proof of Proposition~\ref{prop:virtualorbitequiv}.
    Suppose $x$ is regular and in $\mathcal{K}(I)$. 
    Denote its $\cF$-orbit by $x_n = \cF^n(x)$.
    By Proposition~\ref{prop:virtualorbitequiv} {and its proof}, there exist $\xi_1,\ldots, \xi_k$, $k\leq \cN$ so that 
    $$
    \mathcal{GO}_\mathcal{F}(x) = \bigcup_{i=1}^k\mathcal{GO}_\mathcal{G}(\xi_i),
    $$
    {where $\xi_i = x_{\sigma(i)}$ for some $\sigma(i) \in \N$, $i\in\{1,\cdots,k\}$.}
    Choose $K$ sufficiently large so that $K \geq \sigma(i)$ for each $i$. {Since $x$ is a regular point of $\mathcal{K}(I)$, we may assume that $x_K \in I$ and $x_K$ is regular}. Again, since $\cN$ is the rate of acceleration, each $\xi_i$ is mapped by some iterate of $\mathcal{G}$ into the set $\{x_K, x_{K+1},\ldots, x_{K+\cN-1}\}$.
    Finally, as $\cF^j = \cG^j$ {at all regular points in} $I$ for all $j\leq \cN$, each $\xi_i$ is mapped to $x_{K+\cN}$ by some iterate of $\cG$.
    Thus, $k = 1$ and the proposition follows.
\end{proof}

\subsection{Acceleration of Bowen--Series map for non-uniform lattice}
\begin{theorem}\label{thm:accnon-uniformlattice}
    Let $G$ be a non-uniform lattice. Let $\mathcal{F}_G$ be a Bowen--Series map associated to a core polygon with vertices $a_0,\ldots, a_{N-1}$, as in \S \ref{ss:idealcorpolygon}. Then there exists an acceleration $\mathcal{H}_G$ of $\mathcal{F}_G$ so that
    \begin{itemize}
        \item $\mathcal{H}_G$ is continuous;
        \item every break-point is symmetric.
    \end{itemize}

    More precisely, there exists a subset $S \subseteq \{a_0,\ldots, a_{N-1}\}$ consisting of {unique} representatives of the $G$-orbits and a map $\Phi: \{a_0,\ldots, a_{N-1}\} \longrightarrow S$ sending $a_i$ to {the} point in $S$ that is in the same $G$-orbit so that
    \begin{itemize}[leftmargin=8mm]
        \item $\mathcal{H}_G|_{[a_i, a_{i+1}]}$ is continuous;
        \item $\mathcal{H}_G(a_i) = \Phi(a_i)$;
        \item Each periodic break-points of $\mathcal{H}_G$ are symmetric;
        \item Moreover, let $a \in S$ and let $P_a$ be a generator of the stabilizer of $a$ in $G$. Note that $a$ is a parabolic fixed point of $P_a$, and suppose $P_a$ is repelling in the direction $[a, a+\epsilon)$.
        Then there exists $\epsilon > 0$ so that $\mathcal{H}_G|_{[a, a+\epsilon)} = P_a$ and $\mathcal{H}_G|_{(a-\epsilon, a]} = P^{-1}_a$.
    \end{itemize}
\end{theorem}

\begin{figure}
\captionsetup{width=0.98\linewidth}
    \centering
    \includegraphics[width=0.64\linewidth]{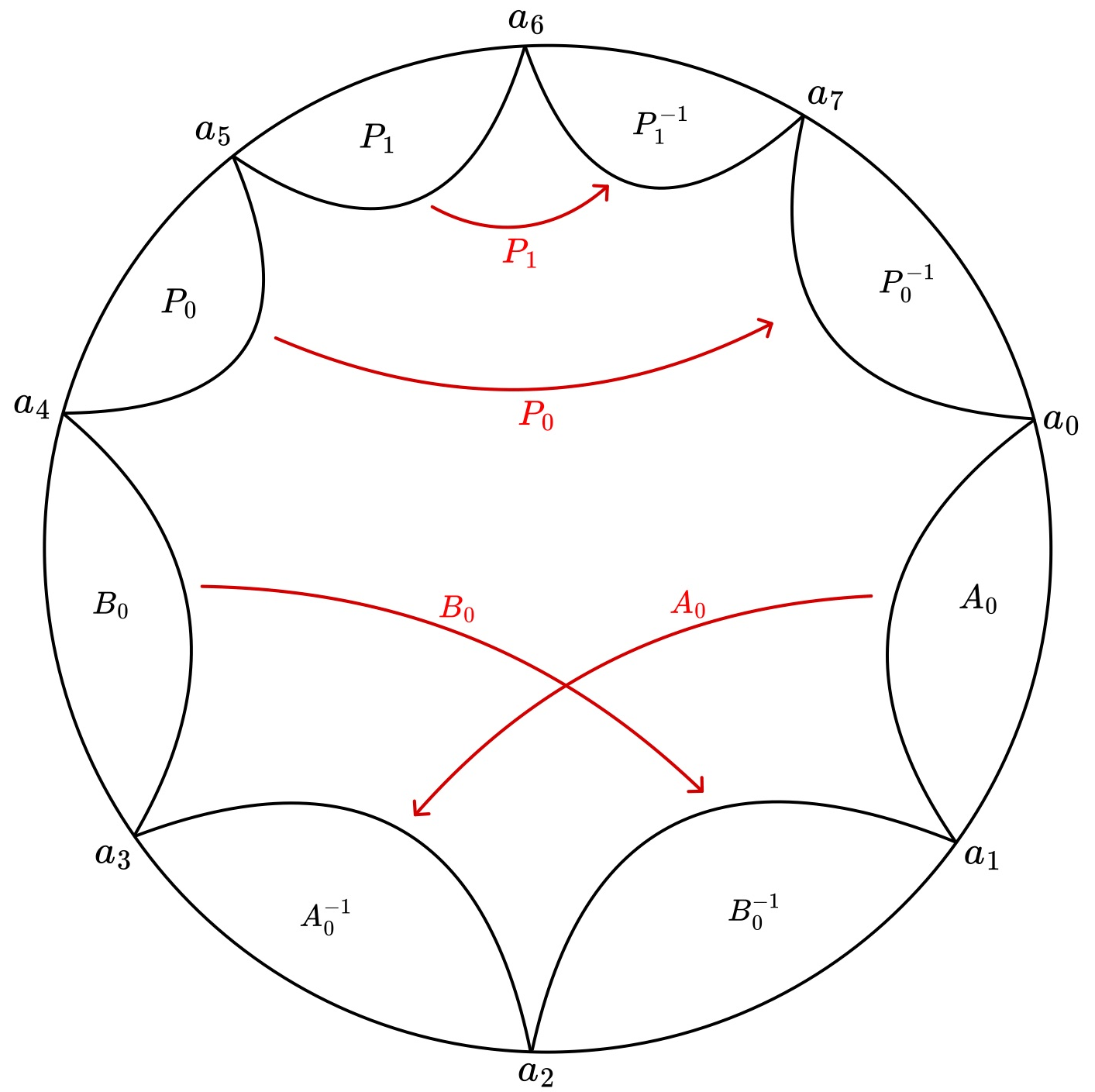}\\
    \includegraphics[width=1\linewidth]{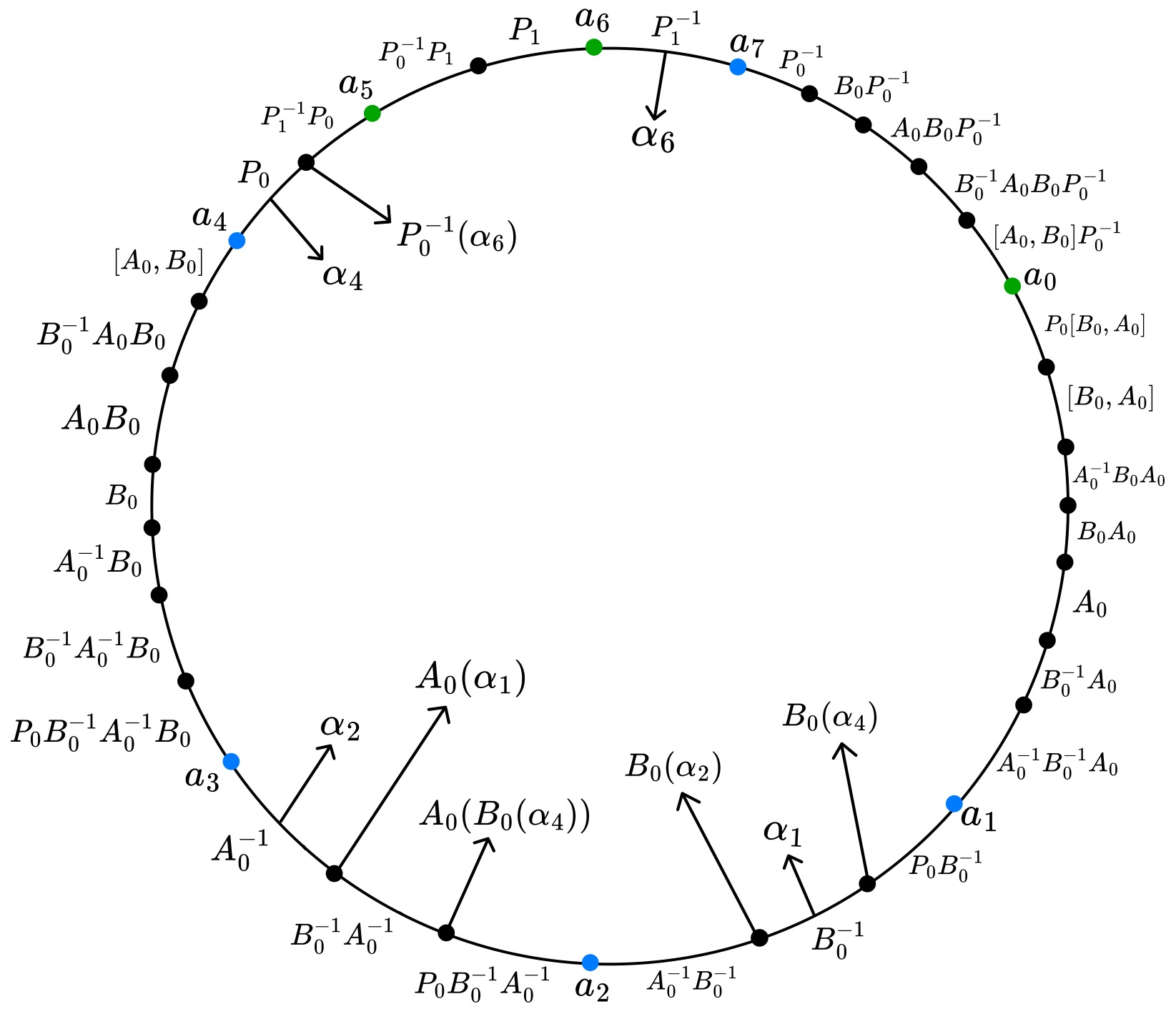}
    \caption{Top: a core polygon and the associated Bowen--Series map $\cF_G$ for a torsion-free Fuchsian group of signature $(1,3)$. Bottom: the corresponding accelerated Bowen--Series map $\cH_G$ with some of its break-points labeled. The green points on $\bS^1$ constitute the set $S$ defined in Theorem~\ref{thm:accnon-uniformlattice}. Note that either $\cH_G(x) = \cF_G(x)$ or $\cH_G(x) = \cH_G\circ \cF_G(x)$. For example, on the clockwise interval $[a_1, a_2]$, $\cH_G(x) = \cF_G(x)$ on $[B_0(\alpha_4), B_0(\alpha_2)]$ and $\cH_G(x) = \cH_G\circ \cF_G(x)$ elsewhere. {Similarly, the maps $\cF_G$, $\cH_G$ act as $A_0^{-1}$, $P_0B_0^{-1}A_0^{-1}$, respectively, on the arc $[A_0(B_0(\alpha_4)),a_2]$, and $\cH_G$ acts as $P_0B_0^{-1}$ on the image arc $\cF_G([A_0(B_0(\alpha_4)),a_2])=[B_0(\alpha_4),a_1]$, yielding the relation $\cH_G(x) = \cH_G\circ \cF_G(x)$ on $[A_0(B_0(\alpha_4)),a_2]$.}}
    \label{acc_bs_fig}
\end{figure}

\begin{proof}
\noindent{\bf Case (1): $G$ is torsion-free.}
    Let us assume that $G$ has signature $(g, n; \infty,\ldots,\infty)$. Let $\mathcal{P}$ be a core polygon. Let $\mathcal{F}_G$ be the Bowen--Series map associated to $\mathcal{P}$ with the notation introduced in \S \ref{ss:idealcorpolygon}.

    Note that the vertices $a_i$ form $n$ orbits under $G$. More precisely, they are 
    \begin{equation*}
    \begin{split}
    \{a_0,\ldots, a_{4g}\}, \{a_{4g+1}, a_{4g+2n-3}\},\ldots, \{a_{4g+n-i}, a_{4g+n+i-2}\}, \ldots,\\
    \{a_{4g+n-2}, a_{4g+n}\}, \{a_{4g+n-1}\}.
    \end{split}
    \end{equation*}
    We pick one representative in each orbit. Concretely, let 
    $$
    S = \{a_0, a_{4g+1}, a_{4g+2},\ldots, a_{4g+n-1}\}.
    $$
    Let $\Phi: \{a_0,\ldots, a_{4g+2n-1}\} \longrightarrow S$ be the map that sends $a_i$ to the unique element in $S$ that is in the same $G$-orbit as $a_i$.

    Note that  $\mathcal{F}_G|_{[a_i, a_{i+1}]}$ is given by a loxodromic element if $i \notin\{4g+n-2, 4g+n-1\}$ and a parabolic element if $i\in\{4g+n-2, 4g+n-1\}$.
    Let $\alpha_i \in [a_i, a_{i+1}]$ be the unique fixed point of $\mathcal{F}_G|_{[a_i, a_{i+1}]}$.
    If $i \notin \{4g+n-2, 4g+n-1\}$, then $\alpha_i \in (a_i, a_{i+1})$ and if $i \in \{4g+n-2, 4g+n-1\}$, then $\alpha_i = a_{4g+n-1}$.
    
    We construct $\mathcal{H}_G$ as follows {(see Figure~\ref{acc_bs_fig}). Let us fix $i\in\{0,\cdots,N-1\}$, and define} $I^0:=[a_i, a_{i+1}]$. We define $\mathcal{H}_{G,0}|_{I^0} = \mathcal{F}_{G}|_{I^0}$, {where the values of $\cF_G(a_i),\cF_{G}(a_{i+1})$ are chosen as limits of $\cF_G(x)$ as $I^0\ni x\to a_i,a_{i+1}$, respectively}. Note that $\mathcal{H}_{G,0}|_{I^0}$ is continuous.
    Suppose that $\mathcal{H}_{G,m-1}|_{I^{0}}$ is defined and continuous.
    Let 
    $$
    I^m = \mathcal{H}_{G,m-1}|_{I^{0}}(I^0) =: [x^{m,-}, x^{m,+}].
    $$
    Abusing notation slightly,  we view $[x^{m,-}, x^{m,+}]$ as an interval in the universal cover
    of the circle, so that its projection to 
    the circle may wrap around the circle multiple times {(we identify sub-intervals of $[x^{m,-}, x^{m,+}]$ with arcs of $\bS^1$ whenever possible)}. Suppose that $x^{m,\pm} \in S$. Then we set $\mathcal{H}_G|_{I^0} = \mathcal{H}_{G,m-1}|_{I^{0}}$.
    Otherwise, suppose that $x^{m,-} \not\in S$. Let $[x^{m,-}, \widetilde{x}^{m,+}] \subseteq I^m$ so that $[x^{m,-}, \widetilde{x}^{m,+}]$ is one of the level $0$ intervals; i.e., $[x^{m,-}, \widetilde{x}^{m,+}] = [a_j, a_{j+1}]$ for some $j$.
    Let $[a_i, \widetilde{\alpha}^{m,+}] = (\mathcal{H}_{G,m-1}|_{I^{0}})^{-1}([x^{m,-}, \alpha_j])$. Similarly, we can define $[\widetilde{\alpha}^{m,-}, a_{i+1}]$ with the understanding that $\widetilde{\alpha}^{m,-} = a_{i+1}$ if $x^{m,+} \in S$.
    We define
    $$
    \mathcal{H}_{G,m}|_{I^{0}}(x) = \begin{cases}
        \mathcal{F}_G \circ \mathcal{H}_{G,m-1}|_{I^{0}}(x)& \text{ if } x\in [a_i, \widetilde{\alpha}^{m,+}] \cup [\widetilde{\alpha}^{m,-}, a_{i+1}],\\
        \mathcal{H}_{G,m-1}|_{I^{0}}(x)& \text{ if } x\in [\widetilde{\alpha}^{m,+}, \widetilde{\alpha}^{m,-}],
    \end{cases}
    $$
    {where the values of $\cF_G(x^{m,\pm})$ are chosen as limits of $\cF_G(x)$ as $I^m\ni x\to x^{m,\pm}$, respectively.}
    Note that $\mathcal{H}_{G,m}|_{I^{0}}$ is a continuous map as $\mathcal{H}_{G,m-1}|_{I^{0}}$ is continuous by induction hypothesis, 
    and the points $\mathcal{H}_{G,m-1}|_{I^{0}}(\widetilde{\alpha}^{m,\pm})$ are fixed points of $\mathcal{F}_G$.
    Note that by construction, we have that either 
    $$
    \mathcal{H}_{G,m}(a_i) = \Phi(a_i),\quad \mathrm{or}\quad \mathcal{H}_{G,m}(a_i) = \mathcal{F}_G^{m+1}|_{I^0}(a_i).
    $$
    Since there exists $K>0$ so that $\mathcal{F}_G^{K}|_{I^0}(a_i) = \Phi(a_i)$, we conclude that $\mathcal{H}_{G,m} = \mathcal{H}_{G,m'}$ for all sufficiently large $m,m'$. Let $\mathcal{H}_G = \mathcal{H}_{G,m}$ for $m$ sufficiently large.
    Then $\mathcal{H}_G|_{[a_i, a_{i+1}]}$ is continuous, and $\mathcal{H}_G(a_i) = \Phi(a_i)$.
    Note that the break-points are eventually mapped to $\alpha_i$ or to $S$, and that each $\alpha_i$ is a fixed point of $\mathcal{F}_G$, and thus a fixed point of $\mathcal{H}_G$. Moreover, all $\alpha_i$ are symmetric break-point for $\cH_G$.
    Finally, we note that by construction there exists $\epsilon > 0$, so that 
    \begin{align*}
    \mathcal{H}_G|_{[a_0, a_0+\epsilon)} = &\ \mathcal{F}^{4g+1}_G|_{[a_0, a_0+\epsilon)}\\
    = &\ P_0[B_{g-1},A_{g-1}][B_{g-2},A_{g-2}]\ldots[B_0,A_0] =: P_{a_0},
    \end{align*}
    where $\mathcal{F}^{4g+1}_G$ is the first return map at $a_0$ {(in the special case when $n=1$, the first return map is given by $\cF_G^{4g}=[B_{g-1},A_{g-1}][B_{g-2},A_{g-2}]\ldots[B_0,A_0]$)}. Thus by Proposition~\ref{prop:symFs}, $P_{a_0}$ is the generator of the stabilizer of $a_0$ which is repelling in $[a_0, a_0+\epsilon)$. A similar computation also gives $\mathcal{H}_G|_{(a_0-\epsilon, a_0]} = P_{a_0}^{-1}$. This verifies the last condition for $a_0 \in S$. The condition for other points in $S$ can be verified by a similar computation. 
    Therefore, all periodic break-points are symmetric. This proves the second statement starting with `more precisely' and thus the theorem for the torsion-free case.
\medskip

\noindent{\bf Case(2): $G$ has torsion points.}
    The idea of the proof is similar, but the computation becomes more complicated, and there are several subcases to consider. We refer the appendix \S~\ref{sec:torsionproof} for more details.
\end{proof}

\subsection{Acceleration of Bowen--Series map for uniform lattice}
\begin{theorem}\label{thm:accuniformlattice}
    Let $G$ be a uniform lattice for a compact hyperbolic surface $\Sigma$ of type $(g, n; \nu_1,\ldots,\nu_n)$. 
    Suppose that $(g,n) \notin\{(0,3), (1,1)\}$.
    Then there exists an acceleration $\cH_{G}$ of a pull-back Bowen--Series map $\mathcal{F}_G$ so that
    \begin{itemize}
        \item $\cH_G$ is continuous;
        \item every break-point is symmetrically hyperbolic.
    \end{itemize}
\end{theorem}
\begin{proof} Since $(g,n) \notin\{(0,3), (1,1)\}$,
     there is a separating essential simple closed curve $\mathcal{C}$ on $\Sigma$ so that after pinching, we obtain two nodal surfaces that are both hyperbolic.
    Let $\widetilde{G}$ be a geometrically finite Bers boundary group associated to $(\Sigma, \mathcal{C})$. Note that the action of $\widetilde{G}$ on the unique invariant component $\Delta_\infty \subseteq \Omega(\widetilde{G})$ is conformally equivalent to $G$ on $\D^*:=\widehat{\C}-\overline{\D}$.

    To construct the Basilica Bowen--Series map, we need to specify the pinched polygon and representative vertex set.
    To this end, we choose two core polygons $\cP_1, \cP_2$ for the two corresponding  nodal surface components. Let $\Delta_i \subseteq \Omega(\widetilde{G})$ be the  components that contain $\cP_i$, and $a \in \partial \Delta_1 \cap \partial \Delta_2$ be the touching point of the two components. 
    Since the original surface is compact, and $\mathcal{C}$ consists only of a single curve, we choose the representative vertex set $\mathcal{PV}_i = \{a\}$ for $i \in\{1, 2\}$. Let $\cF$ be the corresponding Basilica Bowen--Series map. Note that $\cF_G:= \cF_\infty$ is a pull-back Bowen--Series map.
    
    Let $P$ be a limb type puzzle and $x$ be the unique point in $\partial P \cap \Lambda(\widetilde{G})$. Since there is only one pinched curve, we have $\cF(x) = a$.
    Note that $\cF|_{\Delta_i}$ is the induced Bowen--Series map for the corresponding nodal surface. Let $\cH_i$ be the acceleration of this Bowen--Series map as in Theorem~\ref{thm:accnon-uniformlattice}.
    By replacing $\cF|_{\Delta_i}$ with $\cH_i$, we obtain a new $\widetilde{G}$-map $\cH$. The pull-back of $\cH_G:=\cH_\infty$ is an acceleration of $\cF_G$.
    By construction, $\cH_G$ is continuous. Note that since $a\in \partial \Delta_1 \cap \partial \Delta_2$, it corresponds to two ideal boundary points of $\Delta_\infty$. Denote them by $a_{\pm}$. By construction, $a_{\pm}$ are the only potential periodic break-points of $\cH_G$, and they are both fixed. Let $L$ be the generator of the stabilizer of $a_\pm$ in $G$, which is attracting at $a_+$ and repelling at $a_-$. Then by construction, $\cH_G = L$ in a neighborhood of $a_-$ and $\cH_G = L^{-1}$ in a neighborhood of $a_+$. In particular, $\cH_G$ is symmetrically hyperbolic.
\end{proof}

{
\subsection{Strong and weak orbit equivalence}\label{weak_oe_subsec}

In this subsection, we put together the preparatory results from \S \ref{sec:orbit_equiv} and \S \ref{acc_subsec} to study orbit equivalence (or lack thereof) between Fuchsian lattices and the associated accelerated Bowen--Series maps.

\begin{theorem}\label{strong_oe_thm}
Let $G$ be a torsion-free non-uniform Fuchsian lattice with an accelerated Bowen--Series map $\cH_G$. Then, $\cH_G$ is strongly orbit equivalent to $G$.
\end{theorem}
\begin{proof}
Here, a `regular point' refers to a regular point for $\cF_G$. It is easy to verify from the construction of Theorem~\ref{thm:accnon-uniformlattice}, Case (1) that for a regular point $x \in \bS^1$, we have either $\cH_G(x) = \cF_G(x)$ or $\cH_G(x) = \cH_G\circ \cF_G(x)$ (see Figure~\ref{acc_bs_fig}). In either case, $x, \cF_G(x)$ belong to the same grand orbit of $\cH_G$. Therefore, $\mathcal{GO}_{\cF_G}(x) \subseteq \mathcal{GO}_{\cH_G}(x)$. On the other hand, since $\cH_G$ is an acceleration of $\cF_G$, by Proposition~\ref{prop:virtualorbitequiv}, $\mathcal{GO}_{\cH_G}(x) \subseteq \mathcal{GO}_{\cF_G}(x)$. Therefore, $\mathcal{GO}_{\cH_G}(x) = \mathcal{GO}_{\cF_G}(x)$. By Theorem~\ref{thm:orbitequiv}, $\cF_G$ is regular orbit equivalent to $G$.
 Thus, for each regular point $x$ of $\cF_G$, we have $\mathcal{GO}_{\cH_G}(x) = \mathcal{O}_G(x)$.
 On the other hand, each break-point $x \in \mathrm{bk}(\cF_G)$ is eventually mapped (by $\cH_G$) to a fixed point of $\cH_G$ which corresponds to a unique cusp of $\D/G$. Therefore, the theorem follows.
\end{proof}

We are now in a position to prove Corollary~\ref{all_nonuniform_lattice_top_oe_cor}.

\begin{proof}[Proof of Corollary~\ref{all_nonuniform_lattice_top_oe_cor}] We refer the reader to \cite[\S 8]{MM1} for background on topological orbit equivalence and terminology. The proof follows Lemma 8.6 and Theorem 8.7 of \cite{MM1}.

By Theorem~\ref{strong_oe_thm}, $\cH_{G_i}$ is strongly orbit equivalent to $G_i$ on $\bS^1$, $i\in\{1,2\}$ (for $G_1, G_2$ as in the statement of the Corollary). On the other hand, $\cH_{G_i}$ is topologically conjugate to 
$z \mapsto z^{d_i}$ on the circle, for some $d_1,d_2>1$. By \cite[Lemma~8.6]{MM1}, $z^{d_1}\vert_{\bS^1}$ and $z^{d_2}\vert_{\bS^1}$ are topologically orbit equivalent to each other, and hence so are $\cH_{G_1}$ and $\cH_{G_2}$.
It follows that the actions of $G_1, G_2$ on $\bS^1$ are topologically orbit equivalent.

Next, if $\mathbb{D}/G_1, \mathbb{D}/G_2$ are different topologically, then
\begin{enumerate}[leftmargin=8mm]
    \item either $G_1, G_2$ are abstractly non-isomorphic, in which case a topological conjugation is impossible;
    \item  or $G_1, G_2$ are abstractly isomorphic, but the number of cusps (and hence the number of orbits of parabolics in $\bS^1$) are different. This again prevents a topological conjugation.
\end{enumerate}
\end{proof}

\begin{theorem}\label{thm:weakorbitequiv}
    Let $G$ be a lattice of $\PSL_2(\R)$ with an accelerated Bowen--Series map $\cH_G$. Then, $\cH_G$ is weakly orbit equivalent to $G$. If $G$ is uniform, then $\cH_G$ is not strongly orbit equivalent to $G$. 
\end{theorem}}
\begin{proof}
    Let $\cN$ be the rate of acceleration of $\cH_G$ with respect to the Bowen--Series (for non-uniform lattices) or the pull-back Bowen--Series (for uniform lattices) map.
    By construction, there exists an {open} interval $I$ so that $\cH_G = \cF_G$ on $I, \cF_G(I),\cdots, \cF_G^\cN(I)$. By Proposition~\ref{prop:symFs} and Proposition~\ref{prop:basilicabssymmetric}, $\cF_G$ is topologically expanding and primitive. 
    By Proposition~\ref{prop:oe}, $\cH_G$ and $\cF_G$ are orbit equivalent at all regular points in 
    $$
    \mathcal{K}(I) = \{x\in \mathbb{S}^1: |\{n: \cF_G^n(x)\in I\}| = \infty\}.
    $$
    Let 
    $$
    \Lambda_0:= \bigcap_{n\geq 0} \cF_G^{-n}(\mathbb{S}^1\setminus I);\quad \text{and}\quad \Lambda_\infty:= \bigcup_{n\geq 0} \cF_G^{-n} (\Lambda_0).
    $$
    Note that $\bS^1=\mathcal{K}(I) \cup \Lambda_\infty$ is an $\cF_G$-invariant partition of the circle. Since $\cF_G$ is regular orbit equivalent to $G$ on a full-measure subset of $\bS^1$ (by Theorem~\ref{thm:orbitequiv}), it now suffices to show that $\Lambda_\infty$ (or $\Lambda_0$) has zero Lebesgue measure.
    
    Suppose that $G$ is a uniform lattice. Then $\cF_G$ is expanding, and hence $\Lambda_0$ has zero Lebesgue measure.
    {Otherwise, $G$ is a non-uniform lattice. By \cite[Theorems~3,~4]{Sta04} (cf. \cite{BS79}), the Bowen--Series map $\cF_G$ is conservative and ergodic with respect to a (sigma-finite, infinite) invariant measure which is equivalent to the Lebesgue measure on $\bS^1$. Since $I$ is a non-degenerate interval, it follows by ergodicity and conservativity that the fully invariant set $\Lambda_\infty$ (under $\cF_G$) has zero Lebesgue measure.

   Now suppose that $G$ is uniform. Let $\cH_G$ be the acceleration of the pullback Bowen--Series map $\cF_G$. By Propositions~\ref{prop:closedsurfaceOE} and~\ref{prop:nonclosedsurfaceOESt}, there exists a $G$-orbit in $\bS^1$ that is the union of two disjoint $\cF_G$-grand orbits.  By the proof of Proposition~\ref{prop:virtualorbitequiv}, for a point $x\in\bS^1$ (which is regular for $\cF_G$), we have $\mathcal{GO}_{\cH_G}(x) \subseteq \mathcal{GO}_{\cF_G}(x)$. Thus, the above $G$-orbit splits into multiple $\cH_G$-grand orbits. It follows that $\cH_G$ is not strongly orbit equivalent to~$G$.}
\end{proof}

\begin{prop}\label{countable_presentation_prop}
Let $G$ be a lattice of $\PSL_2(\R)$ with an accelerated Bowen--Series map $\cH_G$. Then for each $g \in G$, there exists an open set $\displaystyle\mathscr{J}_g:=\bigcup_{\ell} J_{\ell,g} \subseteq \mathbb{S}^1$ such that
\begin{enumerate}[leftmargin=8mm]
    \item $\mathscr{J}_g$ has full measure, in particular, $\mathscr{J}_g$ is open and dense;
    \item there exist $m_\ell,n_\ell\geq 0$, such that $g=\cH_G^{-m_\ell}\circ\cH_G^{n_\ell}$ on $J_{\ell,g}$.
\end{enumerate}
{Further, if $G$ is uniform, then there exists $g\in G$ such that the associated pairs of non-negative integers $(m_\ell,n_\ell)$ (given by Item~(2) above) form an infinite set.}
\end{prop}
\begin{proof}
  Let $g \in G$. We may assume $g \neq \mathrm{id}$.
  Let $\mathrm{bk}(\cH_G) = \bigcup \cH_G^{-n}\left(\mathrm{bk}\left(\mathcal{A}_0\right)\right)$ be the union of iterated preimages of the break-points and let $\mathrm{fix}(G)$ be the set of fixed points of $G\setminus \{\mathrm{id}\}$.
  Then $\mathrm{bk}(\cH_G) \cup \mathrm{fix}(G)$ is a countable set.
  By Theorem~\ref{thm:weakorbitequiv}, there exists a full-measure set $\mathcal{X} \subseteq \SS^1$ where we have orbit equivalence between $G$ and $\cH_G$.
  For $x\in \mathcal{X} \setminus (\mathrm{bk}(\cH_G) \cup \mathrm{fix}(G))$, since $gx, x$ are in the same $\cH_G$ grand orbit, there exists $n_x, m_x \geq 0$ so that $\cH_G^{m_x}(gx) = \cH_G^{n_x}(x)$. Thus, $gx = \cH_G^{-m_x} \circ \cH_G^{n_x}(x)$.
  Since $x \not\in \mathrm{bk}(\cH_G)$, there exists an open neighborhood $I_x$ so that $\cH_G^{-m_x} \circ \cH_G^{n_x} = h \in G$ on $I_x$.
  Since $hx = gx$, we conclude that $h^{-1}gx = x$. But $x \notin \mathrm{fix}(G)$, so $g = h$, and $\cH_G^{-m_x} \circ \cH_G^{n_x} = g$ on $I_x$.
  Thus, we obtain an open set $\mathscr{J}_g$ containing $\mathcal{X}\setminus (\mathrm{bk}(\cH_G) \cup \mathrm{fix}(G))$ so that the second condition holds.
  Since $\mathcal{X}\setminus (\mathrm{bk}(\cH_G) \cup \mathrm{fix}(G))$ has full measure, so does $\mathscr{J}_g$, and the proposition follows.

 {For the last part, let us assume that $G$ is uniform. By way of contradiction, we suppose that the statement is false; i.e., for each $g\in G$, there exists a finite set $\mathscr{Z}_g$ of pairs of non-negative integers such that $g$ can be expressed as $\cH_G^{-m}\circ\cH_G^{n}$, with $(m,n)\in\mathscr{Z}_g$, on an open, dense subset of $\bS^1$. Since there are only finitely many analytic branches of $\cH_G^{-m}\circ\cH_G^{m}$ for each fixed $(m,n)$, it follows that there exists a finite partition of the circle such that $g$ can be expressed as $\cH_G^{-m}\circ\cH_G^{n}$ on each interval of this partition (with $(m,n)\in\mathscr{Z}_g$ depending on the interval). This clearly implies strong orbit equivalence between $G$ and the accelerated Bowen--Series map $\cH_G$, contradicting Theorem~\ref{thm:weakorbitequiv}.}
\end{proof}

\section{Matings of Fuchsian lattices with polynomial maps}\label{mating_sec}
The following general notion of pinched neighborhoods will be useful in the discussion below. Let $X$ be a closed set in $\widehat{\C}$. A subset $\widehat{X}\subset\widehat{\C}$ is said to be a \emph{pinched neighborhood} of $X$, pinched at finitely many points $x_1,\cdots,x_n\in\partial X$, if $X\subset\widehat{X}$, each point of $X-\{x_1,\cdots,x_n\}$ is an interior point of $\widehat{X}$, but no $x_i$ is an interior point of $\widehat{X}$. If $\widehat{X}$ is a pinched neighborhood of $X$ pinched at $x_1,\cdots,x_n$, then the set $\widehat{X}-\{x_1,\cdots,x_n\}$ is called a \emph{deleted pinched neighborhood} of $X$. Finally, if $X\subset Y\subset\widehat{\C}$ for some closed set $Y$, and $\widehat{X}$ is a pinched neighborhood of $X$ pinched at $x_1,\cdots,x_n$, then $\widehat{X}\cap Y$ is called a \emph{relative pinched neighborhood} of $X$ in $Y$, and $\left(\widehat{X}-\{x_1,\cdots,x_n\}\right)\cap Y$ is called a \emph{relative deleted pinched neighborhood} of $X$ in $Y$.

\subsection{Weakly mateable maps}\label{weak_mateable_subsec}

Recall that an expansive covering map $A: \mathbb{S}^1 \to \mathbb{S}^1$ of degree $d\geq 2$ is  topologically conjugate to the polynomial $z^d\vert_{\mathbb{S}^1}$.  

Recall the notion of conformal Markov maps from \S \ref{sec:cmm}. For a Fuchsian lattice $G$, an $(\bS^1,G)$-map is a conformal Markov map $A:\bS^1\to\bS^1$ such that the maps involved in the piecewise definition of $A$ are elements of $G$. In particular, an $(\bS^1,G)$-map $A$ is piecewise M{\"obius}, but it is allowed to be discontinuous at the break-points of its piecewise definition (cf. \cite[\S 2]{LMM26b}, \cite{MM1}).

\begin{defn}\label{def-pwfm}
	Let $G$ be a Fuchsian lattice. An $(\bS^1,G)$-map $A: \mathbb{S}^1 \to \mathbb{S}^1$ is said to be an $(\bS^1,G)$-covering if $A$ is an expansive covering map of degree greater than one.
	
	Further, an $(\bS^1,G)$-covering $A$ is said to be \emph{mateable} (respectively, \emph{weakly mateable}) if
	\begin{enumerate}[leftmargin=8mm]
		\item\label{cond_orbit_equiv} $A$ is {strongly} orbit equivalent
		(respectively, weakly orbit equivalent) to $G$, and
		
		\item\label{cond-asymhyp} none of the periodic break-points of $A$ is asymmetrically hyperbolic.
	\end{enumerate} 
\end{defn} 

We remark that the conditions imposed on a (weakly) mateable map in the above definition guarantees conformal mateability of $A$ with polynomials (cf. Theorems~\ref{uniform_latice_infused_mating_thm} and~\ref{non_uniform_latice_mating_thm}).

By Theorem~\ref{thm:accnon-uniformlattice}, Theorem~\ref{thm:accuniformlattice} and Theorem~\ref{thm:weakorbitequiv}, an accelerated Bowen--Series map $\cH_G$ associated with a Fuchsian lattice $G$ is an example of a weakly mateable map.

\subsubsection{Compatible weakly mateable maps}\label{compatible_wm_subsubsec}
Let $G_0$ be a Fuchsian lattice, and $\cH_{G_0}$ be an associated (weakly) mateable map.
The map $\cH_{G_0}$ admits an analytic extension to a canonical domain of definition  (cf. \cite[Definition~2.6]{MM1})
given by a pinched neighborhood of $\bS^1$.  The  boundary of this domain is 
given by the union of
\begin{enumerate}[leftmargin=8mm]
\item the finitely many  break-points of $\cH_{G_0}$ on $\bS^1$,
    \item a finite collection of hyperbolic geodesics in $\D$
    connecting the break-points cyclically, and 
    \item a finite collection of hyperbolic geodesics in 
    $\D^*=\widehat{\C}-\overline{\D}$ connecting the break-points cyclically.
\end{enumerate}  
Given such a map $\cH_{G_0}$, one can define compatible (weakly) mateable maps for all groups in $\mathrm{Teich}(G_0)$.
\smallskip

\noindent\textbf{Fuchsian deformations.} For any Fuchsian group $G\in\mathrm{Teich}(G_0)$, there exists a quasiconformal homeomorphism that preserves $\D$ and conjugates the marked group $G_0$ (i.e., the group $G_0$ equipped with a generating set) to the marked group $G$. Such a quasiconformal map is not unique, but all such maps agree on $\bS^1$. We define 
a (weakly) mateable map $\cH_G:\bS^1\to\bS^1$ by conjugating $\cH_{G_0}:\bS^1\to\bS^1$ by this quasiconformal map, and then extend the conjugated map piecewise to its canonical domain of definition. By construction, such a quasiconformal map conjugates $\cH_{G_0}$ to $\cH_G$ on a pinched neighborhood of the circle. We call $\cH_G$ the \emph{(weakly) mateable map of $G$ compatible with $\cH_{G_0}$}. {We denote the intersection of the canonical domain of definition of $\cH_G$ with $\overline{\D}$ (respectively, with $\overline{\D^*}$) by $\cD_G$ (respectively, by~$\cD_G^*$). Note that $\cD_G^*=\eta(\cD_G)$, where $\eta(z)=1/\overline{z}$.}
\smallskip

\noindent\textbf{Bers slice deformations.} {For any marked group $G$ in the Bers slice $\mathcal{B}(G_0)$, there exists a quasiconformal homeomorphism $\xi_G:\widehat{\C}\to\widehat{\C}$ that is conformal on $\D^*$ and conjugates the marked group $G_0$ to the marked group $G$. We set 
$$
\Delta_\infty(G)=\xi_G(\D^*),\quad \Delta_0(G)=\xi_G(\D).
$$ 
The procedure described in the above paragraph defines the (weakly) mateable map $\cH_G$ associated with $G$ and compatible with $\cH_{G_0}$. This map $\cH_G$ is defined on a pinched neighborhood of $\Lambda(G)$.  It is quasiconformally conjugate to $\cH_{G_0}$ in some pinched neighborhood of the limit set such that the conjugacy is conformal on $\Delta_\infty(G)$. As before, we denote the intersection of the canonical domain of definition of $\cH_G$ with $\overline{\Delta_0(G)}$ (respectively, with $\overline{\Delta_\infty(G)}$) by $\cD_G$ (respectively, by~$\cD_G^*$).}

\subsection{Conformal matings as single-valued maps}\label{conf_mating_def_subsec}

In this subsection, we will explicate what it means for a single-valued map to be a conformal mating of a polynomial and a Fuchsian group.

Let $\cH_G:\cD_G (\subset \overline{\D})\longrightarrow\overline{\D}$ be a weakly mateable map associated with a Fuchsian lattice~$G$. Set $d:=\deg\left(\cH_G:\bS^1\to\bS^1\right)$. We will also be working with the symmetrized map
$$
\cH_G:\cD_G\cup\cD_G^*\to\widehat{\C}.
$$

Let $P$ be a degree $d$ monic, centered complex polynomial with a connected Julia set $\mathcal{J}(P)$. We denote its Julia set and filled Julia set by $\mathcal{J}(P)$ and $\cK(P)$, respectively.
The (normalized) \emph{B{\"o}ttcher coordinate} of $P$ is the unique conformal map 
$$
\psi_P:\D\to\mathcal{B}_\infty(P):=\widehat{\C}-\cK(P)
$$ 
conjugating $z^{d}$ to $P$, and tangent to $1/z$ near the origin (cf. \cite[\S 9]{Mil06}).

\subsubsection{Special case: locally connected Julia set}\label{lc_case_subsubsec}

Note that if $\mathcal{J}(P)$ is locally connected, $\psi_P$ extends continuously to $\mathbb{S}^1$ yielding a semi-conjugacy between $z^d\vert_{\mathbb{S}^1}$ and $P\vert_{\mathcal{J}(P)}$. The map $\psi_P:\bS^1\to\mathcal{J}(P)$ is called the (dynamical) \emph{Carath{\'e}odory loop} of $\mathcal{J}(P)$. In this case, 
one can use this Carath{\'e}odory loop to sew a copy of $\overline{\D}$ (equipped with the $\cH_G$-action) outside the filled Julia set $\cK(P)$ (equipped with the $P$-action), giving rise to a \emph{topological mating} between $P$ and $\cH_G$. To do so, we will use the normalized homeomorphism $\mathfrak{h}_G:\mathbb{S}^1\to\mathbb{S}^1$ that conjugates $z^d$ to~$\cH_G$. 

As an intermediate step, we form the disjoint union $\cK(P)\sqcup \overline{\D}$ and define the map 
\begin{center}
	$P\sqcup \cH_G: \cK(P)\sqcup \cD_G\to \cK(P)\sqcup \overline{\D},$\\
	$\left(P\sqcup \cH_G\right)\vert_{\cK(P)}=P,\quad  \left(P\sqcup \cH_G\right)\vert_{\cD_G}=\cH_G.$
\end{center}
Let $\sim_{\mathrm{m}}$ (where `m' denotes mating) be the equivalence relation on $\cK(P)\bigsqcup \overline{\D}$ generated by 
\begin{equation}
\psi_{P}(z)\sim_{\mathrm{m}} \mathfrak{h}_G(\overline{z}),\ \textrm{for all}\ z\in\mathbb{S}^1.
\label{conf_mat_equiv_rel}
\end{equation} 
The map $P\sqcup \cH_G$ descends to a continuous map $P\mate_{\mathrm{top}}\ \cH_G$ on the quotient 
$$
\cK(P)\ \mate\ \overline{\D}\ :=\ \left(\cK(P)\sqcup \overline{\D}\right)/\sim_{\mathrm{m}}\ \cong\ \mathbb{S}^2,
$$ 
and we call this map the \emph{topological mating} $P\mate_{\mathrm{top}} \cH_G$ of $P$ and $\cH_G$.
The dynamical systems $P$ and $\cH_G$ are said to be \emph{conformally mateable} if the topological $2$-sphere $\cK(P)\mate \overline{\D}$ admits a complex structure that turns the topological mating $P\mate_{\mathrm{top}} \cH_G$ into a holomorphic~map, which is referred to as a \emph{conformal mating} 
$$
P\mate\cH_G\equiv P\mate_{\mathrm{conf}} \cH_G
$$
of $P$ and~$\cH_G$.

We refer the reader to \cite[\S 2.3]{MM1} or \cite[\S 3.1]{MM2} for more details and equivalent definitions of conformal matings (in the locally connected case).

\subsubsection{General case via degenerate polynomial-like maps}\label{degen_poly_like_subsubsec}

Since Theorem~\ref{most_fuchsian_mateable_intro_thm} deals with polynomials $P$ whose Julia set $\mathcal{J}(P)$ is connected, but not necessarily locally connected, the above construction of topological matings is not adequate for our purposes. To address this issue, we will work with a more general framework that is modeled on the theory of polynomial-like maps (cf. \cite{DH85}). The following definitions are taken from \cite{LLM24}.

A \textit{pinched polygon} is a set in $\widehat{\C}$ which is homeomorphic to a closed disk quotiented by a finite geodesic lamination, and which has a piecewise smooth boundary. The cut-points of a pinched polygon will be called its \emph{pinched points}, and the non-cut-points where two smooth local boundary arcs meet are called the \emph{corners} of the pinched polygon. (The word `corner' carries no geometric meaning in this context; for instance, two smooth edges of a pinched polygon are allowed to meet at a zero/full angle at a corner.) 

\begin{defn}\label{degenerate_poly_def}
Let $\mathfrak{P}_1, \mathfrak{P}_2\subset \widehat{\C}$ be pinched polygons such that 
\begin{enumerate}[leftmargin=8mm]
\item $\mathfrak{P}_1\subset \mathfrak{P}_2$,
\item the pinched points (respectively, corners) of $\mathfrak{P}_2$ are also pinched points (respectively, corners) of $\mathfrak{P}_1$, and
\item $\partial \mathfrak{P}_1\cap \partial \mathfrak{P}_2$ consists of the corners and pinched points of $\mathfrak{P}_2$.
\end{enumerate}
A \emph{degenerate polynomial-like map} of degree $d$ is a triple $(h, \mathfrak{P}_1, \mathfrak{P}_2)$ such that 
\begin{enumerate}[leftmargin=8mm]
    \item $h:\mathfrak{P}_1 \to \mathfrak{P}_2$ is a degree $d$ continuous map,
    \item $h$ is a proper holomorphic map from each component of $\Int{\mathfrak{P}_1}$ onto some component of $\Int{\mathfrak{P}_2}$, and 
    \item  the singular points of $\partial \mathfrak{P}_1$ are the $h$-preimages of the singular points of $\partial \mathfrak{P}_2$.
\end{enumerate}
\end{defn}

The \textit{filled Julia set} of a degenerate polynomial-like map is defined as 
$$
\cK(h)=\bigcap_{n\geq 0} h^{-n}( \mathfrak{P}_1).
$$ 
The complement of the filled Julia set is called the \emph{escaping set} of $h$, and is denoted by $\cE(h)=\widehat{\C}- \cK(h)$. The common boundary of $\cK(h)$ and $\cE(h)$ is the \emph{Julia set} $\Lambda(h)$.

One can follow \cite[Chapter I.1]{DH85} to construct the \emph{external map} of a degenerate polynomial-like map with connected filled Julia set. The external map of a degree $d$ degenerate polynomial-like map is a piecewise real-analytic, degree $d$, covering map of $\mathbb{S}^1$ that encodes the action of $F$ on the escaping set $\cE(F)$. 

The following notion of conjugacy between two degenerate polynomial-like maps is standard in the theory of polynomial-like maps. In our setting, this will imply a stronger interaction between the dynamics of the polynomial $P$ and the weakly mateable map $\cH_G$ in the \emph{conformal mating} plane.

\begin{defn}\label{hybrid_conj_def}
Two degenerate polynomial-like maps $(h_1,\mathfrak{P}_1, \mathfrak{P}_2)$ and $(h_2,\mathfrak{Q}_1,\mathfrak{Q}_2)$ are said to be \textit{hybrid conjugate} if there exists a quasiconformal map $\Phi\colon \widehat{\C}\to \widehat{\C}$ such that
\begin{enumerate}[leftmargin=8mm]
	\item $\Phi$ sends the pinched points (respectively, corners) of $\mathfrak{P}_2$ to the pinched points (respectively, corners) of $\mathfrak{Q}_2$,
	\item $\Phi$ conjugates $h_1$ to $h_2$, restricted to pinched neighborhoods of $\cK(h_i)$ (pinched at the singular points of $\partial \mathfrak{P}_1, \partial \mathfrak{Q}_1$),
	\item $\overline \partial \Phi\equiv 0$ almost everywhere on $\cK(h_1)$.
\end{enumerate}
\end{defn}

With this preparation, we now come to the notion of single-valued conformal matings that will be used in this paper.

\begin{defn}\label{degen_poly_like_mating_def}
Let $(F, \mathfrak{P}_1, \mathfrak{P}_2)$ be a degenerate polynomial-like map.
\smallskip

\noindent $\bullet$ We say that $(F, \mathfrak{P}_1, \mathfrak{P}_2)$ is a \emph{conformal mating} of $P$ and $\cH_G$ (or simply a \emph{mating} of $P$ and $\cH_G$ for brevity), denoted by $P\mate\cH_G$, if 
\smallskip

\begin{enumerate}[leftmargin=8mm]
    \item[(I)]\label{external_cond}(external group dynamics) the external map of $F$ is given by $\cH_G$; specifically, the map $F:\mathfrak{P}_1\cap\cE(F)\to\cE(F)$ is conformally conjugate to $\cH_G:\D\cap\cD_G\to\D$, and
    \item[(II)]\label{internal_cond}(internal polynomial dynamics) $F\vert_{\cK(F)}$ is topologically conjugate to $P\vert_{\cK(P)}$, where the conjugacy is conformal on $\Int{\cK(F)}$.
\end{enumerate}
\smallskip

\noindent $\bullet$ The map $F$ is said to be a \emph{quasiconformally $P$-infused conformal mating} of $P$ and $\cH_G$ (or simply a \emph{$P$-infused mating} of $P$ and $\cH_G$) if it satisfies Condition~(I) and the following strengthening of Condition~(II):
\smallskip

\begin{enumerate}
\item[(III)]\label{P_infusion_cond}(hybrid polynomial dynamics)  $F$ is hybrid conjugate to a degenerate polynomial-like restriction of $P$ (with filled Julia set~$\cK(P)$).
\end{enumerate}
\smallskip

\noindent $\bullet$ In the special case when $P$ lies in the principal hyperbolic component of monic degree $d$ polynomials, the map $F$ is said to be a \emph{doubly quasiconformally infused conformal mating} of $P$ and $\cH_G$ (or simply a \emph{doubly infused mating} of $P$ and $\cH_G$) if it satisfies Condition~(III) and the following strengthening of Condition~(I):
\smallskip

\begin{enumerate}
\item[(IV)]\label{double_infusion_cond}(hybrid group dynamics) there exists a quasiconformal homeomorphism $\Phi_G:\widehat{\C}\to\widehat{\C}$ such that $\Phi_G(\D)=\cE(F)$, $\overline{\partial}\Phi_G=0$ on $\D$, and $\Phi_G$ conjugates $\cH_G$, restricted to a pinched neighborhood of $\overline{\D}$, to $F$, restricted to a pinched neighborhood of $\overline{\cE(F)}$.
\end{enumerate}
\end{defn}

\begin{remark}\label{hybrid_conj_rem}
Any reasonable notion of conformal mating of $P,\cH_G$ must mandate the mating $F$ to behave like $P$ and $\cH_G$ in two disjoint subsets of the sphere. The notion of mating introduced in \S \ref{lc_case_subsubsec} (in the locally connected case) has the further requirement that the action of $F$ on the common boundary (of these two subsets) is a factor of $\cH_G\vert_{\bS^1}$, in addition to being topologically conjugate to $P\vert_{\mathcal{J}(P)}$. In the non-locally connected case, this requirement needs to be relaxed:
Condition~(I) in Definition~\ref{degen_poly_like_mating_def} only necessitates the  action of $F$ on ideal boundary $\partial^i\cE(F)$ to be given by $\cH_G\vert_{\bS^1}$. If $\mathcal{J}(P)$ is locally connected, then Condition~(I) implies that the action of $F$ on its Julia set $\Lambda(F)=\partial \cK(F)$ is a factor of $\cH_G\vert_{\bS^1}$.

On the other hand, the notion of $P$-infused mating (in Definition~\ref{degen_poly_like_mating_def}) makes up for this deficiency by requiring $F$ to be conjugate to $P$ not only on the filled Julia set $\cK(F)$, but also on a pinched neighborhood of it  encroaching into $\cE(F)$. Differently put, in the dynamical plane of a $P$-infused mating, the dynamical systems $P$ and $\cH_G$ are not just topologically stitched together along the Julia set; rather, the $P$-dynamics quasiconformally extends into the region of the $\cH_G$-dynamics.

A conjugacy between $F$ and $P$ as in Item~(III) (respectively, a conjugacy between $F$ and $\cH_G$ as in Item~(IV)) is called a \emph{hybrid conjugacy}.
\end{remark}

\subsection{Constructing conformal matings of polynomials with accelerated Bowen--Series maps}\label{mating_construct_subsec}

{We now prove the existence of sufficiently regular circle homeomorphisms conjugating accelerated Bowen--Series maps $\cH_G$ to the power map $z^d$ (or to more general Blaschke products). This will play a key role in the construction of degenerate polynomial-like maps as conformal matings. Further, to produce infused matings in the uniform lattice case (see Theorem~\ref{uniform_latice_infused_mating_thm}), we will need a stronger equivariant map; namely, a conjugacy between $\cH_G$ and $z^d$ not only on the circle, but also on a pinched neighborhood of it. 

We refer the reader to \cite{AIM09,LMMN} for background on David homeomorphisms}

\begin{lem}\label{qc_david_conj_circle_lem}
\noindent\begin{enumerate}[leftmargin=8mm]
    \item Let $G$ be a uniform lattice in $\PSL_2(\R)$, and $\cH_G$ be an associated accelerated Bowen--Series map with $d=\deg\left(\cH_G:\bS^1\to\bS^1\right)$. Further, let $B$ be an expanding degree $d$ Blaschke product. Then, there exists a quasiconformal homeomorphism $\mathfrak{h}:\widehat{\C}\to\widehat{\C}$, preserving $\D$, such that $\mathfrak{h}$ conjugates $B$ to $\cH_G$ on a pinched neighborhood of $\bS^1$.

    \item Let $G$ be a non-uniform lattice in $\PSL_2(\R)$, and $\cH_G$ be an associated accelerated Bowen--Series map with $d=\deg\left(\cH_G:\bS^1\to\bS^1\right)$. Then, there exists a homeomorphism $\mathfrak{h}:\bS^1\to\bS^1$ that conjugates $B$ to $\cH_G$, and extends continuously to $\D$ as a David homeomorphism. 
\end{enumerate}
\end{lem}
\begin{proof}
(1) {Let us first note that \cite[Theorem~4.10]{LMMN} provides us with a quasisymmetric conjugacy between $\cH_G\vert_{\bS^1}$ and $B\vert_{\bS^1}$. However, this conjugacy does not automatically admit an equivariant extension to a pinched neighborhood of $\bS^1$. 
Our strategy is to show that the actions of the maps $\cH_G$, $B$ near the unit circle admit quasiconformally compatible `fundamental domains', which will give rise to a quasiconformal homeomorphism $\mathfrak{h}:\D\to\D$, extending quasisymmetrically to the unit circle, that conjugates $\cH_G$ to $B$ on a relative pinched 
neighborhood of $\bS^1$ in $\overline{\D}$. Finally, the desired global quasiconformal map is produced by symmetrizing $\mathfrak{h}\vert_{\overline{\D}}$.

\begin{figure}[htp]
 \captionsetup{width=0.98\linewidth}
    \centering
    \includegraphics[width=0.96\textwidth]{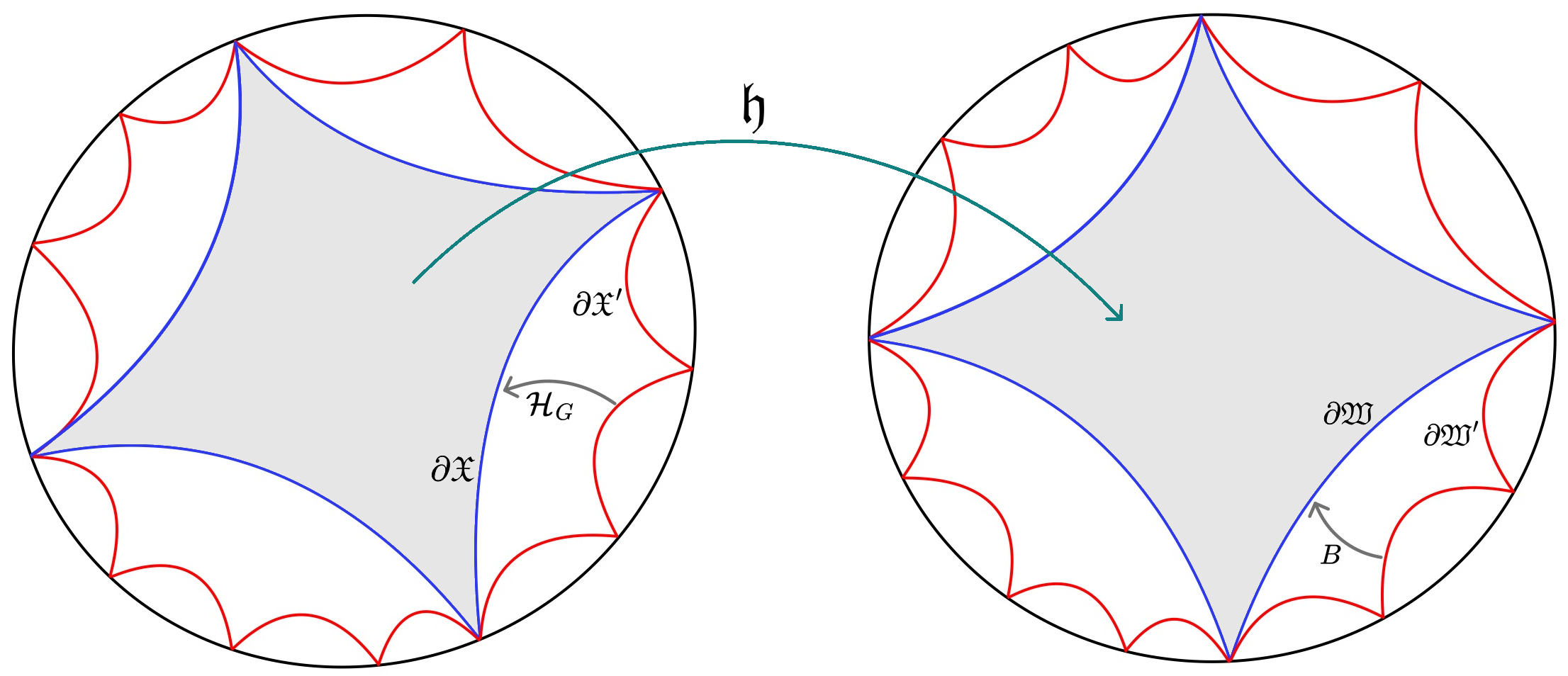}
    \caption{Depicted are the ideal polygons $\mathfrak{X},\mathfrak{X}'$ in the $\cH_G$-plane and the ideal polygons $\mathfrak{W},\mathfrak{W}'$ in the $B$-plane. The quasiconformal homeomorphism $\mathfrak{h}$ carries the gray polygon $\mathfrak{X}$ onto the gray polygon $\mathfrak{W}$.}
    \label{equiv_qc_fig}
\end{figure}

We denote the set of break-points of the piecewise definition of $\cH_G$ by $\mathfrak{p}$. As $\cH_G$ is Markov, we have $\mathfrak{p}\subset\cH_G^{-1}(\mathfrak{p})$. Let $\mathfrak{X},\mathfrak{X}'\subset\D$ be respectively the geodesic polygons having $\mathfrak{p}$, $\cH_G^{-1}(\mathfrak{p})$ as their ideal vertex sets. Then, the polygon $\mathfrak{X}$ is contained in $\mathfrak{X}'$. Since $\cH_G$ is piecewise M{\"o}bius, the edges of $\mathfrak{X}'$ are mapped to those of $\mathfrak{X}$ by $\cH_G$. Further, these polygons have no edges in common (see the construction in \S \ref{ss:bbspbs}, where the polygons have no common edges before acceleration), and their common ideal vertices coincide with $\mathfrak{p}$. After perturbing the geodesic edges of $\mathfrak{X},\mathfrak{X}'$, we can assume in addition to the above properties that the edges of these polygons are smooth, they meet at positive angles, and that $\cH_G(\partial\mathfrak{X})\subset\Int{\mathfrak{X}}\cup\mathfrak{p}$. The smoothness and the positive angle assumptions imply that $\partial\mathfrak{X}, \partial\mathfrak{X}'$ are quasi-circles, a fact that will allow us to apply the Ahlfors-Beurling extension theorem later. (See Figure~\ref{equiv_qc_fig}; the angle condition can be guaranteed by adopting the Klein model of the hyperbolic plane, instead of the Poincar\'e disk model for instance.)

Our next goal is to justify the existence of a homeomorphism $\mathfrak{h}:\overline{\mathfrak{X}'}\hookrightarrow\overline{\D}$, quasiconformal on the interior, with the following properties (see Figure~\ref{equiv_qc_fig}).
\begin{enumerate}[leftmargin=*]
    \item $\mathfrak{h}\left(\cH_G^{-1}(\mathfrak{p})\right)\subset\bS^1$,
    \item The interiors of $\mathfrak{W}:=\mathfrak{h}(\mathfrak{X})$, $\mathfrak{W}':=\mathfrak{h}(\mathfrak{X}')$ are quasi-disks (with $\mathfrak{W}\subset\mathfrak{W}'$),
    \item $B:\partial\mathfrak{W}'\to\partial\mathfrak{W}$ is a degree $d$ covering that maps the edges (respectively, vertices) of the first polygon to the edges (respectively, vertices) of the second.
    \item $\mathfrak{h}\circ\cH_G=B\circ\mathfrak{h}$ on $\partial\mathfrak{X}'$.
\end{enumerate}
Once the map $\mathfrak{h}:\overline{\mathfrak{X}'}\to\overline{\mathfrak{W}}'$ (with the above properties) is at our disposal, lifting $\mathfrak{h}:\mathfrak{X}'-\mathfrak{X}\to\mathfrak{W}'-\mathfrak{W}$ iteratively via $\cH_G, B$ would give us a quasiconformal homeomorphism $\mathfrak{h}:\D\to\D$ that conjugates $\cH_G\vert_{\D-\mathfrak{X}'}$ to $B\vert_{\D-\mathfrak{W}'}$. Further, its quasisymmetric extension $\mathfrak{h}:\bS^1\to\bS^1$ would conjugate $\cH_G$ to $B$ on the circle. (We note here that the equivariance condition $\mathfrak{h}\circ\cH_G=B\circ\mathfrak{h}$ on $\partial\mathfrak{X}'$ is crucial to guarantee
that the lifting procedure yields a continuous map $\mathfrak{h}$.)
\begin{figure}[htp]
 \captionsetup{width=0.98\linewidth}
    \centering
    \includegraphics[width=1\textwidth]{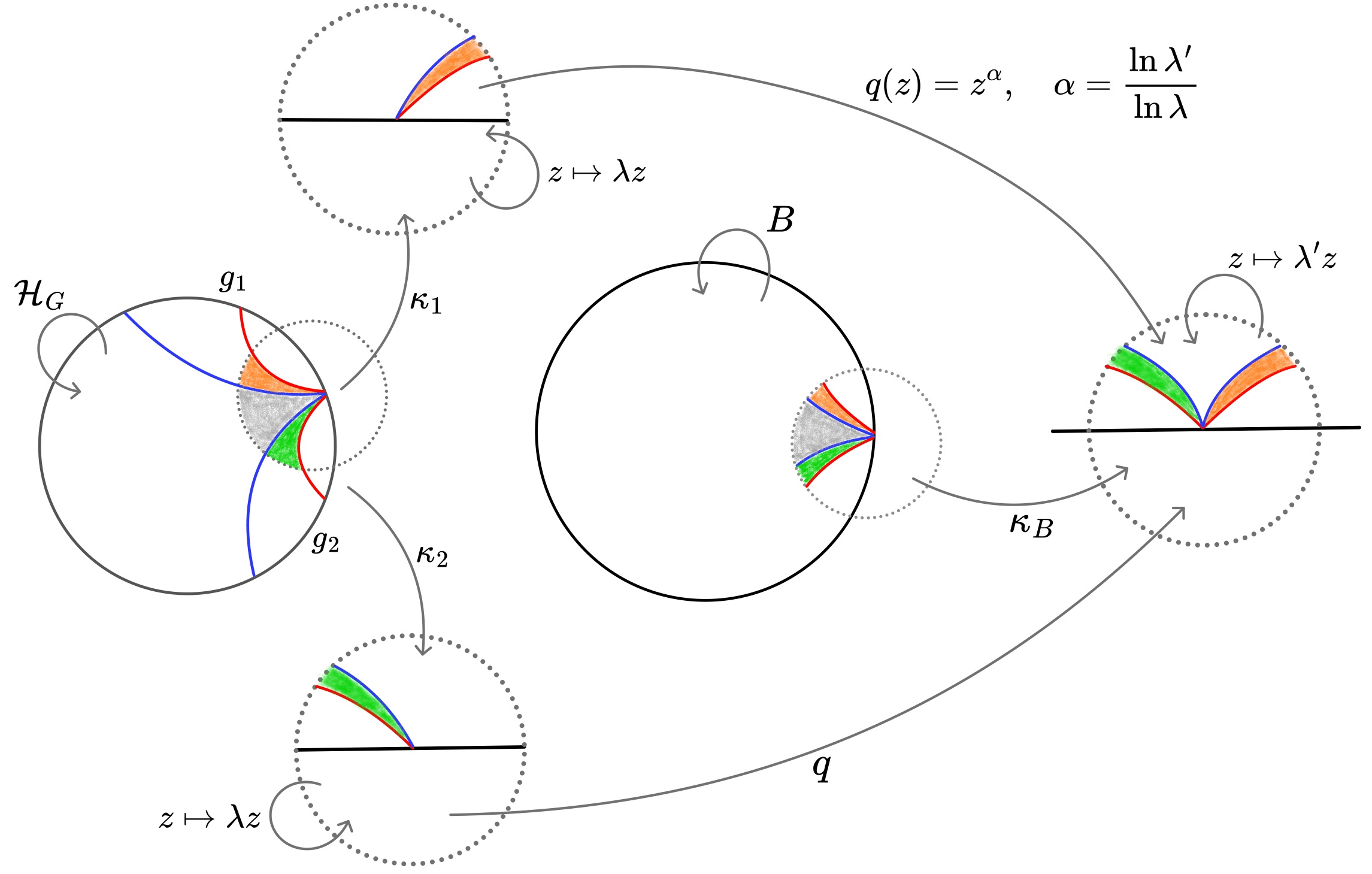}
    \caption{Illustrated is the quasiconformal interpolation step of Lemma~\ref{qc_david_conj_circle_lem}.}
    \label{koenigs_power_fig}
\end{figure}

We now proceed with the construction of the map $\mathfrak{h}:\overline{\mathfrak{X}'}\hookrightarrow\overline{\D}$. Let us fix a topological conjugacy $\mathfrak{f}:\bS^1\to\bS^1$ between $\cH_G$ and $B$. The sets $\mathfrak{f}(\mathfrak{p})$ and $\mathfrak{f}(\cH_G^{-1}(\mathfrak{p}))$ will turn out to be the vertex sets of the polygons $\mathfrak{W},\mathfrak{W}'$, respectively. 

Let $p\in\mathfrak{p}$ be a fixed point of $\cH_G$. Then, $\cH_G$ acts as $g_1, g_2\in G$ on the two sides of $p$ with $g_1(p)=g_2(p)=p$ and $g_1'(p)=g_2'(p)=:\lambda\in\R$ (cf. Theorem~\ref{thm:accuniformlattice}). Let $\kappa_i:N(p)\to\C$ be a local linearizing (Koenigs) coordinate for $g_i$ at $p$ carrying the arc of $\bS^1$ near $p$ to a real interval with $\kappa_i(p)=0$, $i\in\{1,2\}$. We may thus assume that $\kappa_1'(p)=\kappa_2'(p)$. Similarly, let $\kappa_B$ be a local linearizing (Koenigs) coordinate for $B$ at the fixed point $\mathfrak{f}(p)$ carrying the arc of $\bS^1$ near $\mathfrak{f}(p)$ to a real interval with $\kappa_B(\mathfrak{f}(p))=0$. Setting $\lambda':=B'(\mathfrak{f}(p))\in\R$, we note that an appropriate branch of $\mathfrak{q}(z)=z^\alpha$, $\alpha=\frac{\ln\lambda'}{\ln\lambda}$, conjugates $z\mapsto \lambda z$ to $z\mapsto\lambda' z$. We now define a quasiconformal map $\mathfrak{h}$ on $\left(\mathfrak{X}'-\Int{\mathfrak{X}}\right)\cap N(p)$ as $\kappa_B^{-1}\circ \mathfrak{q}\circ\kappa_i$, $i\in\{1,2\}$ (see Figure~\ref{koenigs_power_fig}, where the map $\mathfrak{h}$ is defined as $\kappa_B^{-1}\circ \mathfrak{q}\circ\kappa_1$ on the orange cuspidal region and as $\kappa_B^{-1}\circ \mathfrak{q}\circ\kappa_2$ on the green cuspidal region). By construction, the map $\mathfrak{h}$ thus defined, is quasisymmetric on the quasi-arc $\partial\mathfrak{X}\cap N(p)$ (in particular, $\mathfrak{h}\left(\partial\mathfrak{X}\cap N(p)\right)$ is a quasi-arc); indeed, thanks to the symmetrically hyperbolic assumption on $\cH_G$, the map $\mathfrak{h}$ is asymptotic to $\mathrm{const}\cdot z^\alpha$ near $p$. This allows us to extend $\mathfrak{h}$ quasiconformally to $\mathfrak{X}\cap N(p)$. A similar argument defines the map $\mathfrak{h}$ on the part of $\mathfrak{X}'$ near all the periodic points of $\mathfrak{p}$. On the part of $\mathfrak{X}'$ near the strictly pre-periodic points of $\mathfrak{p}$, we define $\mathfrak{h}$ by lifting the (already defined) map $\mathfrak{h}$ via iterates of $\cH_G$ and $B$, ensuring that $\mathfrak{h}(p)=\mathfrak{f}(p)$ for $p\in\mathfrak{p}$.

It is now straightforward to extend $\mathfrak{h}$ quasiconformally to $\mathfrak{X}$ such that $\mathfrak{h}(\partial\mathfrak{X})$ is a quasi-circle, $\mathfrak{h}:\partial\mathfrak{X}\to\mathfrak{h}(\partial\mathfrak{X})$ is a quasisymmetric homeomorphism, and $B(\mathfrak{h}(\partial\mathfrak{X}))\subset\Int{\mathfrak{h}(\mathfrak{X})}\cup\mathfrak{h}(\mathfrak{p})$. Then, lift $\mathfrak{h}$ to $\partial\mathfrak{X}'$ via $\cH_G$ and $B$ ensuring that $\mathfrak{h}, \mathfrak{f}$ agree on the set $\cH_G^{-1}(\mathfrak{p})$. By the conformality of $\cH_G$ (respectively, of $B$) away from $\mathfrak{p}$ (respectively, on $\bS^1$), the curve $\mathfrak{h}(\partial\mathfrak{X}')$ is a quasi-circle and the map $\mathfrak{h}:\partial\mathfrak{X}'\to\mathfrak{H}(\partial\mathfrak{X}')$ is a quasisymmetric homeomorphism. Finally, set $\mathfrak{W}=\mathfrak{h}(\mathfrak{X})$, and let $\mathfrak{W}'$ be the polygon bounded by $\mathfrak{h}(\partial\mathfrak{X}')$. As $\mathfrak{h}$ is already defined (as a quasiconformal map) on the cuspidal regions of $\mathfrak{X}'-\mathfrak{X}$, a standard application of Ahlfors-Beurling extension theorem permits one to quasiconformally interpolate $\mathfrak{h}$ on the remaining part of $\mathfrak{X}'-\mathfrak{X}$. Clearly, $\mathfrak{h}$ is the desired quasiconformal map that conjugates $\cH_{G}\vert_{\partial\mathfrak{X}'}$ to $B\vert_{\partial\mathfrak{W}'}$.}

(2) It is easily checked that the map $\cH_G$ satisfies the conditions of \cite[Theorem~4.13]{LMMN}, and hence there exists a circle homeomorphism conjugating $z^d$ to $\cH_G$, and admitting a David extension to $\D$ (cf. \cite[Proposition~2.17]{MM1}). Moreover, the map $B$ is quasisymmetrically conjugate to $z^d$; i.e., there exists a circle homeomorphism conjugating $B$ to $z^d$, and admitting a quasiconformal extension to $\D$. The result follows by composing these two conjugacies.    
\end{proof}

{Recall that by the \emph{signature} of a Fuchsian lattice, we mean the signature of the quotient orbifold $\Sigma=\D/G$.}

\begin{theorem}\label{uniform_latice_infused_mating_thm}
Let $G$ be a uniform lattice in $\PSL_2(\R)$ {with signature $(g,n) \notin \{(0,3), (1,1)\}$}, and $\cH_G$ be an associated accelerated Bowen--Series map with $d=\deg\left(\cH_G\vert_{\bS^1}\right)$. 
\begin{enumerate}[leftmargin=8mm]
\item Let $P$ be a degree $d$ polynomial with connected Julia set. Then, there exists a degenerate polynomial-like map that is a $P$-infused mating of $P$ and $\cH_G$.   
\item Let $P$ be a degree $d$ polynomial in the principal hyperbolic component of monic degree $d$ polynomials. Then, there exists a degenerate polynomial-like map that is a doubly infused mating of $P$ and~$\cH_G$.   
\end{enumerate}
\end{theorem}
\begin{proof}
(1) Recall that $\psi_P:\D\to\widehat{\C}-\cK(P)$ is a conformal map conjugating $z^d$ to $P$. Let $\mathfrak{h}:\widehat{\C}\to\widehat{\C}$ be the quasiconformal map from the first part of Lemma~\ref{qc_david_conj_circle_lem} conjugating $B(z)=z^d$ to $\cH_G$ on a pinched neighborhood of $\bS^1$.
Consider the map
\begin{align*}
&\widetilde{F}:\cK(P)\bigsqcup \psi_P\big(\mathfrak{h}^{-1}(\cD_G)\big)\\
&\widetilde{F}=\begin{cases}
    P,\quad \mathrm{on}\ \cK(P),\\
    \left(\mathfrak{h}\circ\psi_P^{-1}\right)^{-1}\circ \cH_G\circ  \left(\mathfrak{h}\circ\psi_P^{-1}\right)\quad \mathrm{on}\ \psi_P\big(\mathfrak{h}^{-1}(\cD_G)\big).
\end{cases}
\end{align*}
The equivariance properties of $\mathfrak{h}, \psi_P$ ensure that the second map in the piecewise definition of $\widetilde{F}$ agrees with $P$ on a relative pinched neighborhood of $\mathcal{J}(P)$ in $\overline{\mathcal{B}_\infty(P)}$. Since finitely many points are quasiconformally removable, it follows that $\widetilde{F}$ is a quasiregular map. 

Now, consider the Beltrami coefficient $\mu$ on $\widehat{\C}$ that is equal to the trivial Beltrami coefficient $\mu_0$ on $\cK(P)$, and equal to $\big(\psi_P\big)_*\left(\mathfrak{h}^*(\mu_0\vert_\D)\right)$ outside $\cK(P)$. As $\cH_G$ preserves $\mu_0\vert_\D$, it follows that $\mu$ is an $\widetilde{F}$-invariant essentially bounded Beltrami coefficient. Let $\Phi_P:\widehat{\C}\to\widehat{\C}$ be a quasiconformal map that solves the Beltrami equation with coefficient $\mu$. Standard arguments (involving Weyl's lemma) show that
$$
F:=\Phi_P\circ\widetilde{F}\circ\Phi_P^{-1}:\cD_F:=\Phi_P\big(\cK(P)\cup\psi_P(\mathfrak{h}^{-1}(\cD_G))\big)\longrightarrow\widehat{\C}
$$
is holomorphic. 

Note that the ideal boundary points of $\cD_G$ are periodic under $\cH_G$, and thus their preimages under $\mathfrak{h}$ are periodic under $z^d$. Since periodic dynamical rays of $P$ land on the Julia set, it is easily seen from the construction that $\cD_F$ is a pinched polygon. Further, $F:F^{-1}\big(\cD_F\big)\to\cD_F$ is a degree $d$ degenerate polynomial-like map with filled Julia set $\cK(F)=\Phi_P\big(\cK(P)\big)$.

Since $\mu=\mu_0$ on $\cK(P)$, we have that $\frac{\overline{\partial}\Phi_P}{\partial\Phi_P}=0$ a.e. on $\cK(P)$. It now follows from the construction of $\widetilde{F}$ that $\Phi_P$ is a hybrid conjugacy between a degenerate polynomial-like restriction of $P$ (with filled Julia set $\cK(P)$) and the degenerate polynomial-like restriction of $F$.

Finally, the map $\Phi_G:=\Phi_P\circ \psi_P\circ\mathfrak{h}^{-1}: \D\longrightarrow\cE(F)$ is quasiconformal, and preserves the standard complex structure. By Weyl's lemma, $\Phi_G$ is conformal; i.e., it is the required conformal conjugacy between $\cH_G$ and $F$.
\medskip

(2) Suppose that $P$ lies in the principal hyperbolic component of monic degree $d$ polynomials. A standard quasiconformal surgery argument gives us the existence of a quasiconformal map $\Psi:\widehat{\C}\to\widehat{\C}$ such that $\Psi\big(\mathcal{B}_\infty(P)\big)=\D$, $\overline{\partial}\Psi=0$ on $\mathcal{B}_\infty(P)$, and $\Psi$ conjugates the restriction of $P$ on a neighborhood of $\overline{\mathcal{B}_\infty(P)}$ to the restriction of $z^d$ on a neighborhood of $\overline{\D}$ (cf. \cite[\S 4.2]{BF14}). In particular, $\Psi$ conjugates $P\vert_{\mathcal{B}_\infty(P)}$ conformally to $z^d\vert_{\D}$. Hence, after possibly post-composing $\Psi$ with a M{\"o}bius map, we can assume that $\Psi\equiv \psi_P^{-1}$ on~$\overline{\mathcal{B}_\infty(P)}$.

Recall the quasiregular map $\widetilde{F}$ defined in part (1) of the theorem. By construction of $\Psi$ and $\mathfrak{h}$ (see Lemma~\ref{qc_david_conj_circle_lem}), we have that 
$$
\widetilde{F}\equiv \left(\mathfrak{h}\circ\Psi\right)^{-1}\circ \cH_G\circ  \left(\mathfrak{h}\circ\Psi\right)
$$
on $\overline{\psi_P\big(\mathfrak{h}^{-1}(\cD_G)\big)}$ as well as on a relative pinched neighborhood of $\mathcal{J}(P)$ in $\cK(P)$.
Thus, the quasiconformal map 
$$
\Phi_G:=\Phi_P\circ\Psi^{-1}\circ\mathfrak{h}^{-1}:\widehat{\C}\to\widehat{\C}
$$ 
is the required conjugacy between $\cH_G$, restricted to a pinched neighborhood of~$\overline{\D}$, and $F$, restricted to a pinched neighborhood of $\overline{\cE(F)}$.
\end{proof}

\begin{remark}\label{dilatation_rem}
Let $G_0$ be a co--compact Fuchsian lattice, $\cH_{G_0}$ be an accelerated Bowen--Series map. 
\begin{itemize}[leftmargin=*]
\item For a Fuchsian group $G\in\mathrm{Teich}(G_0)$, let $\cH_G$ be the accelerated Bowen--Series map of $G$ compatible with $\cH_{G_0}$ (see \S \ref{compatible_wm_subsubsec}).
Let $F_{P,G}$ be the doubly infused mating of $\cH_G$ with $P$ in the principal hyperbolic component of monic degree $d$ polynomials, where $d=\deg(\cH_{G_0}:\bS^1\to\bS^1)$ (produced in part (2) of Theorem~\ref{uniform_latice_infused_mating_thm}). Then, the dilatation of the hybrid conjugacy $\Phi_P$ between $P$ and $F_{P,G}$ depends only on $G\in\mathrm{Teich}(G_0)$. On the other hand, the dilatation of the hybrid conjugacy $\Phi_G$ depends both on $G\in\mathrm{Teich}(G_0)$ and~$P$.

\item Let $G'\in\cB(G_0)$ be a quasi--Fuchsian group, let $G$ be a Fuchsian group determined by the action of $G'\vert_{\Delta_0(G')}$, and $\cH_{G'}, \cH_{G}$ be the accelerated Bowen--Series maps of $G', G$ (respectively) compatible with $\cH_{G_0}$ (see \S \ref{compatible_wm_subsubsec}). With minor adjustments, the proof of Theorem~\ref{uniform_latice_infused_mating_thm} produces the doubly infused mating $F_{P,G'}$ between $\cH_{G'}\vert_{\Delta_0(G')}$ and a polynomial $P$ in the principal hyperbolic component of degree $d$ polynomials. In fact, the conformal mating $F_{P,G'}$ is the same as the doubly infused conformal mating $F_{P,G}$ between the polynomial $P$ and the accelerated Bowen--Series map $\cH_{G}\vert_{\D}$. Thus, by the discussion of the previous paragraph, the dilatation of the hybrid conjugacy $\Phi_P$ between $P$ and $F_{P,G'}$ depends only on $G\in\cB(G_0)$. However, the hybrid conjugacy between $F$ and $\cH_{G'}$ is not the same as the hybrid conjugacy between $F$ and $\cH_{G}$. Since $\cH_{G'}\vert_{\Delta_\infty(G')}$ is conformally conjugate to $\cH_{G_0}\vert_{\D^*}$, it follows that the dilatation of the hybrid conjugacy $\Phi_{G'}$ between $\cH_{G'}$ and $F_{P,G'}$ depends only on $P$, and is independent of $G'\in\cB(G_0)$. This contrast between Fuchsian and Bers slice groups will be relevant in \S \ref{bers_bdry_blaschke_mating_sec}.
\end{itemize}
\end{remark}

\begin{theorem}\label{non_uniform_latice_mating_thm}
Let $G$ be a non-uniform lattice in $\PSL_2(\R)$, and $\cH_G$ be an associated accelerated Bowen--Series map with $d=\deg\left(\cH_G:\bS^1\to\bS^1\right)$. Further, let $P$ be a degree $d$ hyperbolic polynomial with connected Julia set. Then, there exists a degenerate polynomial-like map that is a mating of $P$ and $\cH_G$.   
\end{theorem}
\begin{proof}
By the second part of Lemma~\ref{qc_david_conj_circle_lem}, there exists a topological conjugacy $\mathfrak{h}:\bS^1\to\bS^1$ from $z^d$ to $\cH_G$ that extends continuously as a David homeomorphism of $\D$. Hence, the proof of \cite[Theorem~3.2]{MM2} (involving a David surgery argument) applies verbatim to the current setting. 
\end{proof}

\subsection{Holomorphic dynamical systems generated by single-valued conformal matings}\label{holo_dyn_syst_from_conf_mating_subsec}

The single-valued conformal mating $F=P\mate\cH_G$ between a polynomial $P$ and an accelerated Bowen--Series map $\cH_G$ of a Fuchsian group $G$ constructed so far retains the full dynamics of $P\vert_{\cK(P)}$, but does not capture the full  dynamics of the group $G$ on $\D$. Indeed, the group $G$ appears in $F$ through the weakly mateable map $\cH_G$ (accelerated Bowen--Series map in the present setting).

We will now describe a process to promote $F$ to a \emph{holomorphic dynamical system} that encodes the group action $G\vert_\D$ completely. Further, passing to such a holomorphic dynamical system also provides a way to study algebraicity/transcendentality of the conformal mating~$F$.

\begin{defn}\label{holo_dyn_syst_def}\cite[\S 9.1]{McM96}, \cite[\S 4.1]{MS98}
Let $X$ be a Riemann surface. A holomorphic relation $R\subseteq X\times X$ is a countable union of
$1$-(complex-)dimensional analytic subsets of $X\times X$. A holomorphic dynamical system $(X,\mathscr{R})$ is a collection
$\mathscr{R}$ of holomorphic relations on~$X$.    
\end{defn}
To relate $(X,\mathscr{R})$ to usual iteration theory, one defines compositions of relation as
$$
R_2\circ R_1= \{(x,y) \in X\times X:\ \exists\ z\in X\ \mathrm{such\ that\ } (x,z) \in R_1\ \mathrm{and}\ (z,y) \in R_2 \}.
$$

{
\noindent A conjugacy $\phi$ between holomorphic dynamical systems $(X,\mathscr{R})$ and $(Y,\mathscr{S})$ is a homeomorphism $\phi:X\to Y$ (which may be quasiconformal/conformal) such that
$$
\mathscr{S}=\phi_*(\mathscr{R})=\{\phi_*(R)=(\phi\times\phi)(R): R\in\mathscr{R} \}.
$$
If $\phi_*(R)\in\mathscr{S}$ for some $R\in\mathscr{R}$, we say that $\phi$ carries a relation in $\mathscr{R}$ to a relation in $\mathscr{S}$. The above notion can be naturally extended to local conjugacies between holomorphic dynamical systems (i.e., conjugacies between restrictions of holomorphic dynamical systems to open subsets).}

\begin{defn}\label{holo_dyn_syst_from_conf_mating_def}
{Let $F=P\mate\cH_G$ be the conformal mating between a polynomial $P$ and an accelerated Bowen--Series map $\cH_G$ associated with a Fuchsian group $G$. (If $G$ is a uniform lattice, $F$ is a quasiconformally $P$-infused conformal mating of $P$ and $\cH_G$.)
\begin{enumerate}[leftmargin=*]
    \item We define $\mathrm{gr}(F\vert_{\Int{\cD_F}}^{-m}\circ F\vert_{\Int{\cD_F}}^n)$ to be the relation consisting of all pairs $(z,w)\in\widehat{\C}\times\widehat{\C}$ such that $F\vert_{\Int\cD_F}^m(z), F\vert_{\Int\cD_F}^n(w)$ are defined, and $F\vert_{\Int\cD_F}^m(z) = F\vert_{\Int\cD_F}^n(w)$.
    \item The holomorphic dynamical system $\mathscr{F}$ is defined as the collection of (normalizations of) the maximal analytic continuations of the analytic subsets $\mathrm{gr}(F\vert_{\Int\cD_F}^{-m}\circ F\vert_{\Int\cD_F}^n)$, for $m,n\geq 0$. We call $\mathscr{F}$ the holomorphic dynamical system associated with the conformal mating $F=P\mate\cH_G$.
\end{enumerate} }
\end{defn}

{We also define $\mathscr{P}$ to be the holomorphic dynamical system on $\widehat{\C}$ consisting of (normalizations of) the graphs of $P^{-m}\circ P^n$, $m,n\geq 0$. Note that $\mathscr{P}$ leaves the basin of infinity $\mathcal{B}_\infty(P)$, the filled Julia set $\cK(P)$, and the Julia set $\mathcal{J}(P)$ invariant.}

In what follows, we will identify a holomorphic map with its graph, when there is no scope of confusion.

\subsubsection{The uniform lattice case}\label{uniform_lattice_holo_dyn_syst_subsubsec}
{In this subsection, we consider matings between polynomials and uniform Fuchsian lattices, describing the holomorphic dynamical system $\mathscr{F}$ from Definition~\ref{holo_dyn_syst_from_conf_mating_def} as a combination of a polynomial action and a Fuchsian group action. Because $\mathscr{F}$ is constructed from the conformal mating $F$, which retains only the accelerated Bowen--Series map $\cH_G$ rather than the full group action $G\vert_{\D}$, we must leverage the presentation of the elements $g\in G$ in terms of $\cH_G$ provided by Proposition~\ref{countable_presentation_prop}. To this end, we recall the subsets $\mathscr{J}_g$ of $\bS^1$ associated with $g\in G$, which were constructed in Proposition~\ref{countable_presentation_prop} and are open, dense, and of full Lebesgue measure.}

\begin{theorem}\label{uniform_lattice_holo_dyn_syst_thm}
Let $F=P\mate\cH_G$ be the (quasiconformally) $P$-infused conformal mating between a polynomial $P$ and an accelerated Bowen--Series map $\cH_G$ associated with a uniform Fuchsian lattice $G$ (with signature $(g,n) \notin \{(0,3), (1,1)\}$). Then, the holomorphic dynamical system $\mathscr{F}$ associated with $F$ (see Definition~\ref{holo_dyn_syst_from_conf_mating_def}) is a \emph{quasiconformally $\mathscr{P}$-infused mating} of $\mathscr{P}$ and $G$ in the following sense. 
\begin{enumerate}[leftmargin=*]
\item $\mathscr{F}$ leaves $\cK(F)$ and $\cE(F)$ invariant.
\smallskip

\item {There exists a global quasiconformal map $\Phi_P$ with $\Phi_P(\cK(P))=\cK(F)$, $\overline{\partial}\Phi_P=0$ on $\cK(P)$, such that $\Phi_P^{-1}$ carries the graphs in $\mathscr{F}$, restricted to deleted pinched neighborhoods of~$\cK(F)$, onto the graphs in $\mathscr{P}$, restricted to deleted pinched neighborhoods of $\cK(P)$.}
\smallskip

\item The conformal map $\Phi_G:\D\to\cE(F)$ of Theorem~\ref{uniform_latice_infused_mating_thm} conjugates the holomorphic dynamical system $G\vert_{\D}$ (i.e., consisting of the graphs of $g\vert_{\D}$, where $g\in G$) to $\mathscr{F}\vert_{\cE(F)}$.

\item  Additionally, if $P$ lies in the principal hyperbolic component, then $\mathscr{F}$ is a doubly quasiconformally infused mating of $\mathscr{P}$ and $G$ in the following~sense:

\noindent {$\Phi_G$ extends to a global quasiconformal map that sends the graph of each $g\in G$, restricted to a neighborhood of $\D\cup\mathscr{J}_g$, to a graph in $\mathscr{F}$, restricted to a neighborhood of $\cE(F)\cup\Phi_G(\mathscr{J}_g)$.
Further, any graph $\mathrm{gr}(F\vert_{\Int\cD_F}^{-m}\circ F\vert_{\Int\cD_F}^n)$, restricted to a relative deleted pinched neighborhood of $\Lambda(F)$ in $\cK(F)$, arises this way.}
\end{enumerate}
\end{theorem}

\noindent The proof of the theorem will use the following general fact.

\begin{lem}\label{lem:+measure}
    Let $U \subseteq \C$ be a bounded simply connected domain. Let $a \in \partial U$ and $N$ be a disk containing $a$.
    Let $S$ be a connected component of $N \cap U$ with $a \in \partial S$.
    Then $\partial^I_U S \subseteq \partial^I U \cong \SS^1$ has positive Lebesgue measure.
\end{lem}
\begin{proof}
    Since $U, N$ are bounded and simply connected, $S$ is bounded and simply connected.
   Fix $x_0 \in S \subseteq U$.
    Suppose $\overline{B(a, 2r)} \subseteq N$. 
    Let $S' = \partial S \cap B(a, r)$. Then $S' \subseteq \partial U$.
    Since $S'$ is relatively open in $\partial S$, the harmonic measure $\omega_S^{x_0}(S')$ is positive.
    Since $S \subseteq U$, we have $\omega_U^{x_0}(S') \geq \omega_S^{x_0}(S') > 0$. By pulling back under the uniformization map of $U$, we conclude that the ideal boundary $\partial^I_U S$ has positive Lebesgue measure.
\end{proof}
\smallskip

\begin{proof}[Proof of Theorem~\ref{uniform_lattice_holo_dyn_syst_thm}]
We will use the notation $\Phi_P, \Phi_G$ of Theorem~\ref{uniform_latice_infused_mating_thm}. Recall also that $\cH_G:\cD_G\to\overline{\D}$ acts piecewise via the elements $g_j\in G$, $j\in\{1,\cdots,k\}$, where $\{g_j\}_{j=1}^k$ is a generating set of $G$.
\smallskip

\noindent (1) Note first that $F(\cK(F))=F^{-1}(\cK(F))=\cK(F)$. Hence, any local branch of $F\vert_{\Int{\cD}_F}^{-m}\circ F\vert_{\Int{\cD}_F}^n$ preserves $\cK(F)$. 
On the other hand, as $F\equiv \Phi_G\circ\cH_G\circ\Phi_G^{-1}$ on $\cD_F\cap\cE(F)$, it acts as $\Phi_G\circ g_j\circ\Phi_G^{-1}$ on various hyperbolic half--planes in $\cE(F)$. As $\Phi_G$ is conformal on $\D$, the conformal self--maps $\Phi_G\circ g_j\circ\Phi_G^{-1}$ of $\cE(F)$ provide analytic continuations of the restrictions of $F$ on the said hyperbolic half--planes. Thus, if a branch of $F\vert_{\Int{\cD}_F}^{-m}\circ F\vert_{\Int{\cD}_F}^n$ is defined on an open subset of $\cE(F)$, then this branch extends to $\cE(F)$ as the conformal self--map $\Phi_G\circ g\circ\Phi_G^{-1}:\cE(F)\to\cE(F)$, for some $g\in G$.

Let $f$ be an analytic continuation of $\Phi_G\circ g\circ\Phi_G^{-1}$ (for some $g\in G$) on an open disk $N$ with $N \cap \partial \cE(F) \neq \emptyset$. Then there exists a connected component $S$ of $N \cap \cE(F)$ on which $f = \Phi_G\circ g\circ\Phi_G^{-1}$. Let $a \in \partial S \cap \partial \cE(F) \cap N$. By Proposition~\ref{countable_presentation_prop}, the set $\mathscr{J}_g \subseteq \partial^I \cE(F) \cong \SS^1$ is open and has full measure. Thus, by Lemma~\ref{lem:+measure}, $\partial^I_{\cE(F)} S \cap \mathscr{J}_g$ has positive measure. Since the set of accessible points has full harmonic measure, there exist a point $s \in \partial^I_{\cE(F)} S \cap \mathscr{J}_g$ and an arc $\gamma:[0,1) \rightarrow S$ so that $\lim_{t \to 1}\Phi_G^{-1}(\gamma(t)) = s$. If $s$ lies in the component $J_{l,g}$ of $\mathscr{J}_g$, by Proposition~\ref{countable_presentation_prop}, there exist $p, q$ so that $g=\cH_G^{-p}\circ\cH_G^q$ on $J_{l,g}$. Evidently, the same relation holds in a neighborhood of $J_{l,g}$. Hence, $f(\gamma(t)) = F\vert_{\Int{\cD}_F}^{-p} \circ F\vert_{\Int{\cD}_F}^q(\gamma(t))$ for all $t$ sufficiently close to $1$. Thus, $f = F\vert_{\Int{\cD}_F}^{-p} \circ F\vert_{\Int{\cD}_F}^q$ on~$N$.

This shows that the maximal analytic continuations of $F\vert_{\Int{\cD}_F}^{-m}\circ F\vert_{\Int{\cD}_F}^n$ preserve $\cK(F), \cE(F)$.
\smallskip

\noindent (2) This follows from the definition of the hybrid conjugacy $\Phi_P$ between $P$ and $F$, and from the fact that all branches of $\mathscr{F}$ on $\cK(F)$ are of the form~$F\vert_{\Int{\cD}_F}^{-m}\circ F\vert_{\Int{\cD}_F}^n$.
\smallskip

\noindent (3) {We already know by the first part of the theorem that the analytic continuation of a branch of $F\vert_{\Int{\cD}_F}^{-m}\circ F\vert_{\Int{\cD}_F}^n$ to $\cE(F)$ is given by $\Phi_G\circ g\circ\Phi_G^{-1}$, for some $g\in G$. We claim that the graph of the conformal map $\Phi_G\circ g\circ\Phi_G^{-1}:\cE(F)\to\cE(F)$ lies in $\mathscr{F}$ for each $g\in G$. To this end, recall from Proposition~\ref{countable_presentation_prop} that each $g\in G$ can be expressed as compositions of forward/backward branches of $\cH_G$ on an open neighborhood of $\mathscr{J}_g$, where $\mathscr{J}_g\subset\mathbb{S}^1$ is an open set of full measure. This fact, combined with the arguments of part (1) of this theorem, shows that $\Phi_G\circ g\circ\Phi_G^{-1}\vert_{\cE(F)}$ agrees with some branch of $F\vert_{\Int{\cD}_F}^{-m}\circ F\vert_{\Int{\cD}_F}^n$ on an open subset of $\cE(F)$ (such that the boundary of the open set intersects $\partial\cE(F)$). Hence, the graph of $\Phi_G\circ g\circ\Phi_G^{-1}\vert_{\cE(F)}$ is an analytic continuation of $\mathrm{gr}(F\vert_{\Int{\cD}_F}^{-m}\circ F\vert_{\Int{\cD}_F}^n)$, and thus it belongs to $\mathscr{F}$. 
This proves that the restriction of $\mathscr{F}$ on $\cE(F)\times\cE(F)$ consists precisely of the graphs of the conformal self--maps $\Phi_G\circ g\circ\Phi_G^{-1}\vert_{\cE(F)}$, for $g\in G$.}
\smallskip

\noindent (4) According to Theorem~\ref{uniform_latice_infused_mating_thm}, when $P$ lies in the principal hyperbolic component, the conformal map $\Phi_G:\D\to\cE(F)$ admits a global quasiconformal extension that conjugates $\cH_G$ (restricted to a relative pinched neighborhood of $\bS^1$ in $\overline{\D^*}$) to $F$ (restricted to some relative pinched neighborhood of $\Lambda(F)$ in $\cK(F)$). This implies that $\Phi_G^{-1}$ conjugates the branches of $F\vert_{\Int{\cD}_F}^{-m}\circ F\vert_{\Int{\cD}_F}^n$, restricted to relative deleted pinched neighborhoods of $\Lambda(F)$ in $\cK(F)$, to branches of $\cH_G^{-m}\circ\cH_G^n$, restricted to relative deleted pinched neighborhoods of $\bS^1$ in $\overline{\D^*}$. As the branches of $\cH_G^{-m}\circ\cH_G^n$ are given by elements of the group $G$, it follows that $\Phi_G^{-1}$ conjugates the branches of $F\vert_{\Int{\cD}_F}^{-m}\circ F\vert_{\Int{\cD}_F}^n$, restricted to relative deleted pinched neighborhoods of $\Lambda(F)$ in $\cK(F)$, to elements of $G$, restricted to relative deleted pinched neighborhoods of $\bS^1$ in~$\overline{\D^*}$.

In the opposite direction, by Proposition~\ref{countable_presentation_prop}, any $g\in G$ can be expressed as compositions of forward/backward branches of $\cH_G$ on an open neighborhood of $\mathscr{J}_g$, where $\mathscr{J}_g\subset\mathbb{S}^1$ is an open set of full measure. In light of this fact and the proof of part (2) of Theorem~\ref{uniform_latice_infused_mating_thm}, we conclude that $\Phi_G\circ g\circ\Phi_G^{-1}$ can be expressed as compositions of forward/backward branches of $F\vert_{\Int{\cD}_F}$ on a relative neighborhood of $\Phi_G(\mathscr{J}_g)$ in $\cK(F)$. 
\end{proof}

\begin{remark}\label{closure_holo_dyn_syst_rem}
{Some comments about Items (2) and (4) of Theorem~\ref{uniform_lattice_holo_dyn_syst_thm} are in order. 
\smallskip

\noindent i) The boundary of the domain of definition $\cD_F$ of $F$ intersects the Julia set $\Lambda(F)$ at finitely many points, and the map $F$ may not admit analytic continuations to neighborhoods of these points (see Theorem~\ref{uniform_latice_infused_mating_thm}). As a result, for any fixed $m,n\geq 0$, the branches of $F\vert_{\Int{\cD}_F}^{-m}\circ F\vert_{\Int{\cD}_F}^n$ may not admit analytic continuations to neighborhoods of finitely many points of $\Lambda(F)$. Thus in general, these points do not appear as the first coordinates in the maximal analytic continuations of $\mathrm{gr}(F\vert_{\Int\cD_F}^{-m}\circ F\vert_{\Int\cD_F}^n)$. This is precisely the reason why we were forced to consider deleted pinched neighborhoods of the limit set in Theorem~\ref{uniform_lattice_holo_dyn_syst_thm}. However, the branches of $F\vert_{\Int\cD_F}^{-m}\circ F\vert_{\Int\cD_F}^n$ do extend continuously to these points of non-analyticity, so the closures of the analytic subsets $\mathrm{gr}(F\vert_{\Int\cD_F}^{-m}\circ F\vert_{\Int\cD_F}^n)$ in $\widehat{\C}\times\widehat{\C}$ contain these points of non-analyticity as first coordinates. As $\Phi_P$ conjugates $P\vert_{\mathcal{J}(P)}$ to $F\vert_{\Lambda(F)}$, it follows that $\Phi_P$ carries the graphs in $\mathscr{P}$, restricted to $\mathcal{J}(P)$, to the closures of $\mathrm{gr}(F\vert_{\Int\cD_F}^{-m}\circ F\vert_{\Int\cD_F}^n)$, restricted to~$\Lambda(F)$.
\smallskip

\noindent  ii) When $P$ lies in the principal hyperbolic component, the restriction of each $\Phi_G\circ g\circ\Phi_G^{-1}$, $g\in G$, to a neighborhood of $\cE(F)\cup\Phi_G(\mathscr{J}_g)$ is holomorphic. In fact, the corresponding graph belongs to $\mathscr{F}$. But $\Phi_G\circ g\circ\Phi_G^{-1}$ does not, in general, admit analytic continuations around the points of $\partial\Phi_G(\mathscr{J}_g)$; i.e., the points of $\partial\Phi_G(\mathscr{J}_g)$ do not appear as the first coordinates of the corresponding analytic subset of $\mathscr{F}$. However, the closure of the graph of the holomorphic restriction of $\Phi_G\circ g\circ\Phi_G^{-1}$ contains these points of non-analyticity as first coordinates, and $\Phi_G$ carries the graph of $g\vert_{\Lambda(G)}$ to this closure, restricted to~$\Lambda(F)\times\Lambda(F)$.
}
\end{remark}

 \begin{remark}
 In light of Theorem~\ref{uniform_lattice_holo_dyn_syst_thm}, it is legitimate to denote the holomorphic dynamical system $\mathscr{F}$ associated with the conformal mating $F=P\mate\cH_G$ as $\mathscr{P}\mate G$. We note that the holomorphic dynamical system $\mathscr{P}\mate G$ depends not only on the polynomial $P$ and the group $G$, but also on the choice of the accelerated Bowen--Series map $\cH_{G}$.  
\end{remark}

\subsubsection{The non-uniform lattice case}\label{non_uniform_lattice_holo_dyn_syst_subsubsec}

Recall that given a hyperbolic polynomial $P$ and an accelerated Bowen--Series map $\cH_G$ associated with a non--uniform Fuchsian lattice $G$ (with $\deg{(P)}=\deg{(\cH_G\vert_{\bS^1})}$), Theorem~\ref{non_uniform_latice_mating_thm} provides us with a degenerate polynomial-like map $F=P\mate\cH_G$ that is a mating of $P$ and $\cH_G$. Let $\mathscr{F}$ be the holomorphic dynamical system associated with $F$ (cf. Definition~\ref{holo_dyn_syst_from_conf_mating_def}). The proof of Theorem~\ref{uniform_lattice_holo_dyn_syst_thm} now applies verbatim to this setting and yields the following mating description of $\mathscr{F}$.

\begin{theorem}\label{non_uniform_lattice_holo_dyn_syst_thm}
Let $F=P\mate\cH_G$ be the conformal mating between a hyperbolic polynomial $P$ and an accelerated Bowen--Series map $\cH_G$ associated with a non--uniform Fuchsian lattice $G$.
Then, the holomorphic dynamical system $\mathscr{F}$ associated with $F$ is a mating of $\mathscr{P}$ and $G$ in the following sense. 
\begin{enumerate}[leftmargin=*]
\item $\mathscr{F}$ leaves $\cK(F)$ and $\cE(F)$ invariant.

\item There exists a homeomorphism $\Phi_P: \cK(P)\to \cK(F)$, conformal on the interior, such that $\Phi_P^{-1}$ carries the graphs in $\mathscr{F}$, restricted to $\cK(F)$, onto the graphs in $\mathscr{P}$, restricted to $\cK(P)$.
    
\item  There exists a conformal map $\Phi_G:\D\to\cE(F)$ that conjugates the holomorphic dynamical system $G\vert_{\D}$ to $\mathscr{F}\vert_{\cE(F)}$.
\end{enumerate}
\end{theorem}

As in the uniform lattice case, we will denote the holomorphic dynamical system $\mathscr{F}$ associated with the conformal mating $F=P\mate\cH_G$ as $\mathscr{P}\mate G$.

\subsection{Non-algebraicity of matings}\label{non_alg_mating_subsec}

In this subsection, we will show that the degenerate polynomial-like maps $F$ constructed in Theorems~\ref{uniform_latice_infused_mating_thm} and~\ref{non_uniform_latice_mating_thm} are typically not algebraic, although they are matings of two algebraic maps. We will also establish a general result to the effect that failure of algebraicity is not a specific feature of our construction of accelerated Bowen--Series maps. Passing to the broader framework of transcendental correspondences is unavoidable in order to realize the matings between polynomials and uniform Fuchsian lattices or higher-genus non-uniform Fuchsian lattices. In the light of the main theorems of \cite{MM1,MM2,LLM24}, this provides a complete classification in the torsion-free case.

\subsubsection{Non-algebraicity of matings arising from accelerated Bowen--Series maps}\label{non_alg_special_subsubsec}

The following two results establish failure of algebraicity of matings using certain specific properties of the accelerated Bowen--Series maps. The purpose of analyzing this special case (which is subsumed by a general certificate of transcendentality proved in \S \ref{non_alg_general_subsubsec}) is to illustrate two
different mechanisms for failure of algebraicity. 
\begin{itemize}[leftmargin=8mm]
    \item For a uniform lattice $G$, the transcendentality of $F=P\mate\cH_G$ is a consequence of the lack of strong orbit equivalence between $\cH_G$ and $G$.
    \item For a torsion-free non-uniform lattice $G$ of higher genus, although the accelerated Bowen--Series map $\cH_G$ is strongly orbit equivalent to $G$, a suitable monodromy action on the stabilizer of a point turns out to be aperiodic, destroying algebraicity. 
\end{itemize}

\begin{prop}\label{non_alg_prop_uniform}
Let $F=P\mate\cH_G$ be the $P$-infused conformal mating between a polynomial $P$ and an accelerated Bowen--Series map $\cH_G$ associated with a uniform Fuchsian lattice $G$. Then, the maximal analytic continuation of $F$ gives rise to a transcendental correspondence on an open subset of the Riemann sphere. 
\end{prop}
Thus, the holomorphic dynamical system $\mathscr{F}=\mathscr{P}\mate G$ constructed in Theorem~\ref{uniform_lattice_holo_dyn_syst_thm} is generated by an iterated transcendental correspondence.
\begin{proof}
By way of contradiction, let us suppose that $F$ is an algebraic function. Then, every germ defining $\mathscr{F}$ has at most finitely many distinct analytic continuations.

By Proposition~\ref{countable_presentation_prop}, there exist 
\begin{enumerate}[leftmargin=8mm]
    \item $g\in G$,
    \item an open, full-measure set $\displaystyle\mathscr{J}_g=\bigcup_{\ell=1}^\infty J_{\ell,g}$, and
    \item an infinite set $\mathscr{Z}_g=\big\{(m_\ell,n_\ell):\ell\in\N\big\}\subset\N_{\geq 0}\times\N_{\geq 0}$
\end{enumerate} 
such that $g=\cH_G^{-m_\ell}\circ\cH_G^{n_\ell}$ on $J_{\ell,g}$. By the proof of Theorem~\ref{uniform_lattice_holo_dyn_syst_thm}, the graph of $\Phi_G\circ g\circ\Phi_G^{-1}$ belongs to $\mathscr{F}$. It now follows that the maximal analytic continuation of this graph contains the graphs of infinitely many distinct holomorphic germs of the form $F^{-m_\ell}\circ F^{n_\ell}$, $(m_\ell,n_\ell)\in\mathscr{Z}_g$. In particular, this would force the analytic continuation to have infinitely many branch points. This is a contradiction. Hence, the maximal analytic continuation of $F$ defines a transcendental correspondence. 
\end{proof}

\begin{prop}\label{non_alg_prop_nonuniform}
Let $F=P\mate\cH_G$ be the conformal mating between a polynomial $P$ and an accelerated Bowen--Series map $\cH_G$ associated with a non-uniform Fuchsian lattice $G$ of genus $g \geq 1$ without torsion elements. Then, the maximal analytic continuation of $F$ gives rise to a transcendental correspondence on an open subset of the Riemann sphere.  
\end{prop}
\begin{proof}
    The failure of algebraicity comes from the fact that a certain monodromy around a branch point is infinite.
    Indeed, let $a \in [a_0, a_1]$ be the break-point where $\cH_G$ acts as $A_0$ on one side and $B_0A_0$ on the other. Note that $\cH_G(a) = \alpha_3$, where $\alpha_3$ is the repelling fixed point of $B_0$ (see Figure~\ref{acc_bs_fig}).

    By way of contradiction, suppose that $F$ is an algebraic function. Let $\widetilde{a}$ be the point in the domain of $F$ corresponding to $a$.
    Let $\gamma$ be a small loop around $\widetilde a$. 
    We show that analytic continuation of $F$ along $\gamma$ generates infinitely many distinct branches.

    For simplicity of the presentation, we assume that $\cK(F)$ is a Jordan disk. The general case can be proved with minor adjustments.
    Choose the basepoint $x_0\in \gamma$ on the limit set $\Lambda(F)$ so that the corresponding branch is 
    $$
    \cH_G = A_0.
    $$
    (Here and below, we abuse the notation, as the analytic branch should be conjugated by $\Phi_G$, i.e., it is $\Phi_G \circ A_0 \circ \Phi_G^{-1}$.)
    Orient the curve $\gamma$ so that it first passes through $\cK(F)$ and then through $\cE(F)$. We may assume that $\gamma$ intersects $\Lambda(F)$ at two points $x_0, x_1$.
    After analytic continuation along the first half of $\gamma$, the branch at $x_1$ is 
    $$
    \cH_G = B_0A_0.
    $$
    Continuing along the second half of $\gamma$ through $\cE(F)$, returns to $x_0$ with branch 
    $$
    B_0A_0 = \cH_G^2.
    $$
    Repeating the same argument, the next passage through $x_1$ gives
    $$
    \cH_G^2 = B_0(B_0A_0) = B_0^2A_0.
    $$
    Continuing along the second half of $\gamma$ gives the branch
    $$
    B_0^2A_0 = B_0 B_0 A_0 = \cH_G^3.
    $$
    Thus, we see that the analytic continuation along $\gamma$ at $x_0$ takes the form
    $$
    A_0 \to B_0A_0 \to B_0^2 A_0 \to \cdots \to B_0^nA_0 \to \cdots.
    $$
    Since $B_0$ has infinite order, these branches are distinct.
    This shows the monodromy around $\widetilde{a}$ is infinite, a contradiction. Therefore, $F$ is not algebraic.
\end{proof}

\subsubsection{Non-algebraicity of matings: the general case}\label{non_alg_general_subsubsec}

We now expound the key underlying mechanism in the proof of Proposition~\ref{non_alg_prop_nonuniform}: aperiodicity of induced monodromy actions on stabilizers of suitable points. This will, in fact, show that the non-algebraicity of matings for uniform lattices and higher-genus non-uniform lattices is unavoidable, irrespective of whether the construction of the matings uses accelerated Bowen--Series maps.

To this end, consider a Fuchsian lattice $G$ and a weakly mateable $(\bS^1,G)$-covering $A:\bS^1\to\bS^1$ in the sense of Definition~\ref{def-pwfm}. Let $P_0(z)=z^d$, where $d$ is the degree of $A\vert_{\bS^1}$. By Theorems~\ref{uniform_latice_infused_mating_thm} and~\ref{non_uniform_latice_mating_thm}, there exists a conformal mating $F=P_0\mate A$. We recall the conjugacies $\Phi_{P_0}:\overline{\D}=\cK(P_0)\to\cK(F)$ (between $P$ and $F$) and $\Phi_G:\overline{\D}\to\overline{\cE(F)}$ (between $A$ and $F$).
Further, let $\mathscr{F}$ be the holomorphic dynamical system generated by $F$ as in Definition~\ref{holo_dyn_syst_from_conf_mating_def}. 
Since $A$ is weakly mateable, the conclusions of Theorems~\ref{uniform_lattice_holo_dyn_syst_thm} and~\ref{non_uniform_lattice_holo_dyn_syst_thm} remain valid for $\mathscr{F}$.

\begin{prop}\label{alg_is_impossible_prop}
\noindent\begin{enumerate}[leftmargin=8mm]
    \item If $G$ is uniform, then $F$ is transcendental; i.e., the holomorphic dynamical system $\mathscr{F}=\mathscr{P_0}\mate G$ cannot be generated by an iterated algebraic correspondence.

    \item Let $G$ be non-uniform, and $F$ algebraic. Then,
    \begin{enumerate}[leftmargin=6mm]
        \item all break-points of $A$ are parabolic,
        \item the local monodromy of $F$ at each branch point has order two, and
        \item $\D/G$ has genus zero.
    \end{enumerate}
\end{enumerate}    
\end{prop}
\begin{proof}
As shown in Proposition~\ref{non_alg_prop_uniform}, failure of strong orbit equivalence between $A$ and $G$ would force $F$ to be transcendental. Hence, we may assume that $A$ is strongly orbit equivalent to $G$. Thus, for each $g\in G$, there exists a finite partition of $\bS^1$ into intervals on each of which $g\equiv A^{-\alpha}\circ A^\beta$ for some $(\alpha,\beta)$ depending on the interval.

Note that the homeomorphism $\mathfrak{h}=\Phi_G^{-1}\circ\Phi_P:\bS^1\to\bS^1$ conjugates $P_0(z)=z^d$ to $A$. It follows that
$$
\widehat{G}:=\mathfrak{h}^{-1}\circ G\circ\mathfrak{h}
$$
is a group of $d$-adic piecewise linear homeomorphisms acting on $\R/\Z\cong\bS^1$. For $u\in G$ and $t\in\bS^1$, we will use the notation 
$$
\widehat{u}=\mathfrak{h}^{-1}\circ u\circ\mathfrak{h},\hspace{4mm} \widehat{t}=\mathfrak{h}^{-1}(t),\quad \textrm{and}\quad \widetilde{u}=\Phi_G\circ u\circ\Phi_G^{-1},\hspace{4mm} \widetilde{t}=\Phi_G(t).
$$

Now let $g_0\in G$ be one of the maps in the piecewise definition of $A$, and $I\subset\bS^1$ the maximal interval such that $A\vert_I=g_0$. Let $q$ be the right endpoint of $I$ (here left/right are understood at the level of the universal cover), and $\gamma$ a sufficiently small loop around $\widetilde{q}$. We continue $\widetilde{g}_0$ analytically along $\gamma$, and observe that the continuation through $\cK(F)$ (which is given by $F$) comes back to $\cE(F)$ as some element $\widetilde{g}_1$ such that $\widetilde{g}_0(\widetilde{q})=\widetilde{g}_1(\widetilde{q})$ (and hence $g_0(q)=g_1(q)$). By the maximality of $I$, we have that
$$
g_1=h^{n_1}\circ g_0,\ \textrm{for some } n_1\in\Z-\{0\},
$$
where $h$ is a generator of the stabilizer of $g(q)$ in $G$. Similarly, as we analytically continue $\widetilde{g}_1$ along $\gamma$, the continuation through $\cK(F)$ (which is given by a branch of $F^{-\alpha_1}\circ F^{\beta_1}$) comes back to $\cE(F)$ as some element $\widetilde{g}_2$ with $\widetilde{g}_1(\widetilde{q})=\widetilde{g}_2(\widetilde{q})$. Thus, $g_2=h^{n_2}\circ g_0$, for some $n_2\in\Z$. Iterating the above procedure, we obtain an \emph{auxiliary dynamical system} on the stabilizer of $g_0(q)$ in $G$ given~by
$$
\Z\to\Z,\quad n_r\mapsto n_{r+1}.
$$
We will now find an explicit formula for the above monodromy action on $\Z$. Let the left and right slopes of $\widehat{g}_0$ at $\widehat{q}$ (respectively, of $\widehat{h}$ at $\widehat{g}_0(\widehat{q})$) be $d^{l_-}$ and $d^{l_+}$ (respectively, $d^{k_-}$ and $d^{k_+}$). In fact, as $A$ acts as $g_0$ on the maximal interval $I$ with $q\in\partial I$, it follows that $l_-=1$ and $l_+\neq 1$. Thus, the left slope of $\widehat{h}^{n_r}\circ\widehat{g}$ and  
and right slope of $\widehat{h}^{n_{r+1}}\circ\widehat{g}$ are $d^{1+n_rk_-}$ and $d^{l_++n_{r+1}k_+}$, respectively. However, near $\widehat{q}$, these two maps are restrictions of a common analytic germ of the form $P_0^{-\alpha_r}\circ P_0^{\beta_r}$. Hence, we have
$$
l_++n_{r+1}k_+=1+n_rk_-\ \implies\ n_{r+1}=\frac{k_-}{k_+}n_r+\frac{(1-l_+)}{k_+};
$$
i.e., the monodromy action on $\Z$ is given by the iteration of the affine map $L(x)=\frac{k_-}{k_+}x+\frac{(1-l_+)}{k_+}$.
\smallskip

\noindent\textbf{Case I: $h$ is hyperbolic.} In this case, $h$ is either contracting on both sides of $g(q)$ or expanding on both sides of $g(q)$. It follows that $k_\pm$ are simultaneously positive or simultaneously negative; i.e., $k_-/k_+>0$. Since $l_+\neq 1$, it now follows that the orbit of $0$ under $L$ is infinite. Hence, analytic continuation of $\widetilde{g}=\widetilde{h}^0\circ \widetilde{g}$ along $\gamma$ produces infinitely many distinct germs $\widetilde{h}^{n_r}\circ\widetilde{g}$, proving transcendentality of $F$. 
\smallskip

\noindent\textbf{Case II: $h$ is parabolic.} In this case, $h$ is contracting on one side of $g(q)$ and expanding on the other side, implying that $k_-/k_+<0$. So the orbit of $0$ under $L$ is finite if and only if $k_+=-k_-$; i.e., $L(x)=-x+\frac{(1-l_+)}{k_+}$ with $L^2=\textrm{id}$.

Now let $G$ be uniform. Then the element $h$ is necessarily hyperbolic, and hence the mating $F=P_0\mate A$ must be transcendental by Case~I.

On the other hand, if $G$ is non-uniform with $F$ algebraic, then Case~I implies that no break-point of $A$ can be hyperbolic. Further, Case~II implies that the local monodromy at a branch point of $F$ must have order two.

It remains to show that when $G$ is non-uniform and $F$ is algebraic, the quotient orbifold $\D/G$ has genus zero.
Let $S$ denote the underlying 2-manifold of $\D/G$.

By strong orbit equivalence between $A$ and $G$ (which is necessary for algebraicity of $F$), the M{\"o}bius maps appearing in the piecewise definition of $A$ form a generating set for $G$ (cf. \cite[Lemma~2.4]{MM1}). Since every break-point is parabolic by (2a) of this proposition, it follows that $G$ can be generated by finitely many parabolic elements along with at most one extra element. 

If $S$ has genus $g$ and $n$ punctures, then the first homology group $H_1(S)$ has rank $(2g+n-1)$. Since all conjugates of a based loop going around a single puncture of $S$ (representing a parabolic in $G$) are homologous, it follows that the punctures (and hence all parabolics in $G$) generate a subgroup of rank at most $n-1$ in $H_1(S)$. Since $G$ has at most one more generator beyond parabolics (by the previous paragraph), we must have $g=0$.
\end{proof}

\begin{remark}\label{rmk-gen0}
    A bit more work shows that the orbifold  in Case 2 of Proposition~\ref{alg_is_impossible_prop} can have at most two orbifold points (see \cite{LMM26c,BLL26} for algebraic correspondences realizing matings of cusped triangle groups with rational maps). We sketch the argument below. 
    
    Suppose not. Let $z_1, \cdots, z_k$, $k > 2$, denote the orbifold points. Let $\Delta\subset S$ denote a polygon with cyclically ordered vertices   $z_1, \cdots, z_k$. Let $* \in S$
    denote a base-point lying outside $\Delta$. Then any based loop at $*$ determines a Coxeter word in the sense of \cite[Section 4.2]{manning-ruffoni}. We choose based loop representatives 
    $\sigma_1,\cdots, \sigma_n$ lying \emph{outside} $\Delta$, and going around exactly  one of the $n$ punctures of $S$ 
    exactly once. Let $P_1, \cdots, P_n \in G$ denote the corresponding group elements. 
    Then the elements of $G$ representing (generators of) stabilizers of parabolics are all of the form $\gamma P_i \gamma^{-1}$. Each of these gives a trivial Coxeter word with respect to $\Delta$.
    On the other hand, a loop going round exactly one $z_j$ intersects
   the  two  edges of $\Delta$ incident on $z_j$ exactly once. Since the independent homotopy classes of the Coxeter words arising from such paths has cardinality at least $k-1 > 1$, one extra element of $G$ is inadequate to generate $G$ along with the conjugates of the parabolics.
\end{remark}

\begin{remark}\label{infinite_order_mono_rem}
In Proposition~\ref{non_alg_prop_uniform}, the transcendentality of the matings between polynomials $P$ and accelerated Bowen--Series maps $\cH_G$ of uniform lattices was deduced from the failure of strong orbit equivalence between $\cH_G$ and $G$. However, Proposition~\ref{alg_is_impossible_prop} shows that even strong orbit equivalence cannot salvage algebraicity in the uniform situation due to aperiodicity of monodromy actions on stabilizers of hyperbolic break-points.

On the other hand, the source of non-algebraicity in Propositions~\ref{non_alg_prop_nonuniform} and~\ref{alg_is_impossible_prop} is the same; namely, aperiodicity of monodromy actions on stabilizers of hyperbolic break-points. The only known examples of mateable $(\bS^1,G)$-coverings satisfying the necessary conditions of Proposition~\ref{alg_is_impossible_prop} (part 2) are Bowen--Series and higher Bowen--Series maps of punctured spheres (cf. \cite{MM1}). 
\end{remark}

\section{Bers slice of matings, and Carath{\'e}odory compactification}\label{compactification_sec}

Let $G_0$ be a Fuchsian lattice {with signature $(g,n) \notin \{(0,3), (1,1)\}$}, and let $\cH_{G_0}$ be an accelerated Bowen--Series map of $G_0$ with $d=\deg\left(\cH_{G_0}\vert_{\bS^1}\right)$. For $G\in\mathrm{Teich}(G_0)$, the accelerated Bowen--series map of $G$ compatible with $\cH_{G_0}$ is denoted by $\cH_G$.

Let $P$ be a degree $d$ monic, centered polynomial in the principal hyperbolic component. As before, we normalize the B{\"o}ttcher coordinate $\psi_P:\D\to\mathcal{B}_\infty(P)$ to be tangent to $1/z$ at $\infty$. The dynamical ray of $P$ at angle $0$ lands at a repelling fixed point of $P$ (cf. \cite[\S 18]{Mil06}), and we call it the marked repelling fixed point of $P$. Pick a conformal isomorphism $\cK(P)\longrightarrow\overline{\D^*}$ (where $\D^*=\widehat{\C}-\overline{\D}$) that sends the attracting fixed point of $P$ to $\infty$, and the marked repelling fixed point of $P$ to $1$. By conjugating $P\vert_{\cK(P)}$ by this conformal isomorphism, we obtain an expanding degree $d$ Blaschke product $B$ with an attracting fixed point at $\infty$. 
{
\begin{remark}
We have adopted a convention that differs slightly from the standard one; namely, we model the action of $P$ on the basin of infinity by the action of $z^d$ on $\D$, and the action of $P$ on the filled Julia set by the action of a Blaschke product $B$ on $\D^*$. This choice of convention will be convenient in \S \ref{bers_bdry_blaschke_mating_sec} (cf. Lemma~\ref{bers_bdry_group_blaschke_mv_mating_lem}).    
\end{remark}}

We define $\mathscr{B}$ to be the holomorphic dynamical system on $\widehat{\C}$ consisting of (normalizations of) the graphs of $B^{-m}\circ B^n$, $m,n\geq 0$. Evidently, $\mathscr{B}$ leaves $\bS^1$, $\D$, and $\D^*$ invariant; and $\mathscr{P}\vert_{\cK(P)}$ is topologically conjugate to $\mathscr{B}\vert_{\overline{\D^*}}$ where the conjugacy is conformal on $\Int{\cK(P)}$.

Recall that by Theorem~\ref{uniform_latice_infused_mating_thm} and Theorem~\ref{non_uniform_latice_mating_thm}, there exist degenerate polynomial-like maps $F=P\mate\cH_G$, for $G\in\mathrm{Teich}(G_0)$ and polynomials $P$ as above. We set
$$
\Delta_\infty(F):=\Int{\cK(F)}.
$$
Then, there exists a conformal ma $\Phi_B^G:\D^*\to\Delta_\infty(F)$ that conjugates $B$ to $F=P\mate\cH_G$; and hence, $F$ can be regarded as a conformal mating of $B$ and $\cH_G$ as well. Abusing notation, we will now denote $F$ by~$B\mate\cH_G$. The associated holomorphic dynamical system $\mathscr{F}$ of Definition~\ref{holo_dyn_syst_from_conf_mating_def} can be regarded as a mating between the holomorphic dynamical system $\mathscr{B}\vert_{\overline{\D^*}}$ and the group action $G\vert_{\D}$ (cf. Theorem~\ref{uniform_lattice_holo_dyn_syst_thm} and Theorem~\ref{non_uniform_lattice_holo_dyn_syst_thm}). Hence, we may denote $\mathscr{F}$ by $\mathscr{B}\mate G$.

Following \cite[Appendix~A]{LLM24} (cf. \cite[\S 4]{LMM26a}), we define the \emph{horizontal Bers slice through the polynomial $P$} (or equivalently, through its Blaschke model $B$) in the moduli space of degenerate polynomial-like maps as
$$
\mathfrak{B}(P)\equiv\mathfrak{B}(B):=\big\{F=B\mate\cH_G: G\in\mathrm{Teich}(G_0)\big\},
$$
up to M{\"o}bius conjugacy.
By Theorem~\ref{uniform_lattice_holo_dyn_syst_thm} and Theorem~\ref{non_uniform_lattice_holo_dyn_syst_thm}, this Bers slice can be identified with the subspace
{
$$
\faktor{\Big\{\mathscr{F}=\mathscr{B}\mate G:G\in\mathrm{Teich}(G_0)\Big\}}{\PSL_2(\C)}
$$}
of the moduli space of holomorphic dynamical systems on $\widehat{\C}$.

\begin{remark}
As mentioned in Remark~\ref{dilatation_rem}, minor modifications in the proofs of Theorem~\ref{uniform_latice_infused_mating_thm} and Theorem~\ref{non_uniform_latice_mating_thm} show that the elements of $\mathfrak{B}(B)$ can be regarded as degenerate polynomial-like maps/holomorphic dynamical systems realizing matings between the Blaschke product $B$ and groups $G$ lying in the Bers slice $\cB(G_0)$.
\end{remark}

After possibly conjugating the conformal mating $F=B\mate\cH_{G}$ by a M{\"o}bius map, we may assume that $\Phi_B^G$ satisfies 
$$
\Phi_B^G(\zeta)=\zeta+O(1/\zeta)\quad \mathrm{as}\quad \zeta\to\infty;
$$
i.e., the map $\Phi_B^G$ lies in the compact space $\Sigma_0$ of normalized univalent functions (cf. \cite[Theorem~1.7]{Pom75}).
With such normalization, the map 
$$
\E:\mathscr{F}=\mathscr{B}\mate G\mapsto \Phi_B^G
$$
embeds $\mathfrak{B}(B)$ into the space $\Sigma_0$, and the closure  
$$
\overline{\E(\mathfrak{B}(B))}=\overline{\big\{\Phi_B^G:G\in\cB(G_0)\big\}}\subset\Sigma_0
$$
is compact. 

One can associate a holomorphic dynamical system to each point of $\partial\E(\mathfrak{B}(B))$. To this end, consider a sequence of holomorphic dynamical systems $\mathscr{F}_n=\mathscr{B}\mate G_n$ in $\mathfrak{B}(B)$ (with associated conformal matings $F_n=B\mate\cH_{G_n}$) such that the sequence $G_n$ leaves all compact sets in $\cB(G_0)$. After passing to a subsequence, we may assume that the normalized conformal maps $\Phi_B^{G_n}$
converge locally uniformly to the conformal map $(\Phi_B^\infty:\D^*\to\widehat{\C}) \in \partial\E(\mathfrak{B}(B))$. As 
$$
F_n=\Phi_B^{G_n}\circ B\circ \left(\Phi_B^{G_n}\right)^{-1}\quad \mathrm{on} \quad \Delta_\infty(F_n),
$$
by setting
$$
F_\infty:=\Phi_B^{\infty}\circ B\circ \left(\Phi_B^{\infty}\right)^{-1}\quad \mathrm{on} \quad \Delta_\infty(F_\infty):=\Phi_B^\infty(\D^*),
$$
we conclude that the pointed disk $\left(\Delta_\infty(F_\infty),\infty\right)$ is the limit of $\left(\Delta_\infty(F_n),\infty\right)$, and $F_\infty\vert_{\Delta_\infty(F_\infty)}$ is the limit of $F_n\vert_{\Delta_\infty(F_n)}$ in the Carath{\'e}odory topology on pointed disks and holomorphic functions, respectively (cf. \cite{McM94}).
{Now define $\mathscr{F}_\infty$ to be the holomorphic dynamical system given by (normalizations of) the maximal analytic continuations of the analytic subsets $\mathrm{gr}(F_\infty\vert_{\Delta_\infty(F_\infty)}^{-m}\circ F_\infty\vert_{\Delta_\infty(F_\infty)}^n)$, for $m,n\geq 0$.} By declaring that the sequence $\mathscr{F}_n$ in $\mathfrak{B}(B)$ converges to $\mathscr{F}_\infty$, we obtain a compactification of $\mathfrak{B}(B)$ that is naturally homeomorphic to $\overline{\E(\mathfrak{B}(B))}$.

However, a priori, the dynamical systems $\mathscr{F}_\infty$ on $\partial\mathfrak{B}(B)$ may not be defined on the whole sphere, and as such may not arise as matings. In the special case when $\D/G_0$ is a $k$-times punctured sphere, the holomorphic dynamical systems in $\mathfrak{B}(B)$ are given by iterated algebraic correspondences of bi-degree $(2k-3):(2k-3)$  on the Riemann sphere \cite{MM2,LLM24}. It was shown in \cite{LMM26a} that in this case, the boundary $\partial\mathfrak{B}(B)$ also consists of algebraic correspondences of the same bi-degree. We note that ``absence of degree drop'' of the algebraic correspondences does not automatically follow from Carath{\'e}odory convergence of the conformal matings $F_n\vert_{\Delta_\infty(F_n)}$ and requires additional work (e.g., rescaling limit arguments).

In the next section, we will take a closer look at $\partial\mathfrak{B}(B)$ in the case when $\D/G_0$ is a closed surface.

\begin{remark}\label{geom_conv_rem}
The space $\mathfrak{B}(B)$ of holomorphic dynamical systems can also be endowed with the geometric topology (cf. \cite[\S 9]{McM96}). Inspired by analogous phenomena in Kleinian groups, in a future work we will construct sequences of dynamical systems in $\mathfrak{B}(B)$ whose geometric limits are strictly larger than their limits in the ``Carath{\'e}odory topology'' defined above.
\end{remark}

\section{Matings of Bers boundary groups with Blaschke products}\label{bers_bdry_blaschke_mating_sec} 

Let $G_0$ be a uniform Fuchsian lattice {with signature $(g,n) \notin \{(0,3), (1,1)\}$} having an accelerated Bowen--Series map 
$$
\cH_{G_0}:\cD_{G_0}\cup\cD_{G_0}^*\to\widehat{\C}
$$ 
with $d=\deg\left(\cH_{G_0}:\bS^1\to\bS^1\right)$; where $\cD_{G_0}^*=\eta(\cD_{G_0})$, and $\eta(z)=1/\overline{z}$, cf. \S \ref{weak_mateable_subsec}. Recall also the notation $\D^*=\eta(\D)$.
Let $\cB(G_0)$ be the Bers slice of $G_0$. For each $G\in\overline{\cB(G_0)}$, the domain of discontinuity $\Omega(G)$ has a $G$-invariant component, denoted by $\Delta_\infty(G)$, on which the $G$-action is conformally conjugate to $G_0\vert_{\D^*}$. After possibly conjugating $G$ by a M{\"o}bius map, we assume that $\infty\in \Delta_\infty(G)$. There exists a conformal isomorphism
$$
\xi_G:(\D^*,\infty)\to(\Delta_\infty(G),\infty)
$$
conjugating $G_0\vert_{\D^*}$ to $G\vert_{\Delta_\infty(G)}$. The conformal map $\xi_{G}:\D^*\to\Delta_\infty(G)$ extends continuously to the boundary, yielding a topological semi-conjugacy $\xi_{G}:\bS^1\to\Lambda(G)$ between the limit set actions of $G_0$ and $G$ (see \cite{Mj14,MS13}). This defines a continuous map 
\begin{equation}
\xi_G:\left(\overline{\D^*},\infty\right)\to\left(\overline{\Delta_\infty(G)},\infty\right),
\label{ct_eqn}
\end{equation}
and its restriction $\xi_G:\bS^1\to\Lambda(G)$  is the \emph{Cannon--Thurston map} for $G$. 
We will call the compact set 
$$
\mathscr{K}(G):=\widehat{\C}-\Delta_\infty(G)
$$ 
the \emph{filled limit set} of $G$.

Recall that for $G\in\mathcal{B}(G_0)$, the conformal map $\xi_G$ extends to a quasiconformal homeomorphism of $\widehat{\C}$ conjugating $G_0$ to $G$. We denote this global quasiconformal map by $\xi_G$, and note that 
$$
\cH_G:\cD_G\cup\cD_G^*:=\xi_G(\cD_{G_0})\cup\xi_G(\cD_{G_0}^*)\longrightarrow\widehat{\C}
$$
is the accelerated Bowen--Series map of $G$ compatible with $\cH_{G_0}$ (see~\S \ref{compatible_wm_subsubsec}).

\subsection{Multi--valued accelerated Bowen--Series maps}\label{multi_valued_abs_bers_bdry_subsec}

Fix $G\in\partial\mathcal{B}(G_0)$. The map $\xi_{G}$ semi-conjugates $\cH_{G_0}:\cD_{G_0}^*\to\overline{\D^*}$ to a \emph{multi--valued} piecewise $G$-map 
$$
\cH_{G}:\xi_{G}\left(\cD_{G_0}^*\right)\to\overline{\Delta_\infty(G)}.
$$ 
Observe that $\xi_{G}\left(\cD_{G_0}^*\right)$ is a relative pinched neighborhood of $\Lambda(G)$ in $\overline{\Delta_\infty(G)}$ (see Figure~\ref{mvm_fig}, cf. \cite[Lemma~7.6]{MM1}). If $G$ has accidental parabolics, then $\cH_{G}$ admits an analytic continuation to the complement of a collection of ideal polygons in $\Omega(G)-\Delta_\infty(G)$ as a single--valued piecewise $G$-map. We denote the domain of definition of this analytically continued multi--valued map $\cH_{G}$ by $\cD_{G}$, and call $\cH_{G}:\cD_{G}\to\widehat{\C}$ the \emph{multi--valued accelerated Bowen--Series map} of $G$ compatible with~$\cH_{G_0}$.

Clearly, $\cH_{G}$ is single--valued on $\cD_{G}-\Lambda(G)$. In fact, it turns out that the multi--valuedness of the map is $\cH_{G}$ supported on a strict subset of $\Lambda(G)$. To see this, we define the set
$$
\RT_G:=\bigcup_{\substack{j_1, j_2=1,\\ j_1\neq j_2}}^k \Big(\xi_{G}\left(I_{j_1}\right)\bigcap\xi_{G}\left(I_{j_2}\right)\Big)\subsetneq\Lambda(G),
$$
where $\{I_j\}_{j=1}^k$ is a partition of $\mathbb{S}^1$ such that $\cH_{G_0}$ acts as $g_j$ on $I_j$.

\begin{lem}\label{supp_mv_abs_lem}
The multi--valued accelerated Bowen--Series map $\cH_{G}$ is single--valued, holomorphic on $\cD_{G}-\RT_G$.
\end{lem}
\begin{proof}
Let $y\in\Lambda(G)$, and $\xi_{G}^{-1}(y)=\{x_1,\cdots,x_r\}$ for some $r\geq 1$. If $y\in\Lambda(G)-\RT_G$, then there exists $j(y)\in\{1,\cdots,k\}$ such that $\xi_{G}^{-1}(y)\subset I_{j(y)}$. Hence, $\cH_{G_0}(x_i)=g_{j(y)}(x_i)$, for all $i\in\{1,\cdots,r\}$. By definition, there exists an element $\pmb{g}_{j(y)}\in G$ corresponding to $g_{j(y)}\in G_0$ satisfying $\pmb{g}_{j(y)}\vert_{\Delta_\infty(G)}=\xi_{G}\vert_{\D^*}\circ g_{j(y)}\vert_{\D^*}\circ\left(\xi_{G}\vert_{\D^*}\right)^{-1}$. By continuity, 
$$
\cH_{G}(y)=\xi_{G}\left(g_{j(y)}(x_i)\right)=\pmb{g}_{j(y)}(y).
$$
Thus, $\cH_{G}=\xi_{G}\circ\cH_{G_0}\circ\xi_{G}^{-1}$ acts as $\pmb{g}_{j(y)}$ in a neighborhood of $y$ (recall that $\RT_G$ is closed), and hence $\cH_{G}$ is single--valued, holomorphic on $\Lambda(G)-\RT_G$. The same is clearly true on $\cD_{G}-\Lambda(G)$ as well.
\end{proof}

\begin{remark}\label{abs_discont_rem}
The map $\cH_{G}$ is never single--valued on the entire limit set for a Bers boundary group $G$ without accidental parabolics. However, if $G_0$ is a Fuchsian punctured sphere group (in which case the accelerated Bowen--Series map $\cH_{G_0}$ is a standard Bowen--Series map), there exist geometrically finite Bers boundary groups $\ G$ such the $\cH_{G}=\xi_{G}\circ\cH_{G_0}\circ\xi_{G}^{-1}$ is single--valued on the limit set $\Lambda(G)$ (see \cite[\S 7.2]{MM1}). 
\end{remark}

\subsection{Holomorphic dynamical systems combining Bers boundary groups with Blaschke products}\label{bers_bdry_group_blaschke_mating_subsec}

We continue with the above setup. 

\subsubsection{Construction of multi--valued matings} The next result is an analog of Theorem~\ref{uniform_latice_infused_mating_thm} for Bers boundary groups.

\begin{lem}\label{bers_bdry_group_blaschke_mv_mating_lem}
Let $G\in\partial\mathcal{B}(G_0)$, let $\cH_{G}$ be the multi--valued map of \S \ref{multi_valued_abs_bers_bdry_subsec}, and let $B$ be an expanding degree $d$ Blaschke product. Then, there exists a multi--valued map $F=B\mate\cH_{G}:\cD_F\to\widehat{\C}$ such that  
\begin{enumerate}[leftmargin=*]
 \item there exists an $F$-invariant partition  $\widehat{\C}=\mathscr{K}(F)\bigsqcup\Delta_\infty(F)$, 
    
\item there exists a conformal map $\Phi_B:\D^*\rightarrow\Delta_\infty(F)$ that conjugates $B$ to $F$, and

\item there exists a global quasiconformal map $\Phi_{G}$, conformal on $\mathscr{K}(G)$, such that $\Phi_{G}$ conjugates $\cH_{G}$, restricted to a pinched neighborhood of $\mathscr{K}(G)$ (wherever defined), to $F$, restricted to a pinched neighborhood of $\mathscr{K}(F)$ (wherever defined).
    \end{enumerate}
Further, $F$ is a single--valued, holomorphic map outside $\Phi_{G}(\RT_G)$. 
\end{lem}
\begin{proof}

\begin{figure}[htp]
 \captionsetup{width=0.98\linewidth}
    \centering
    \includegraphics[width=1\textwidth]{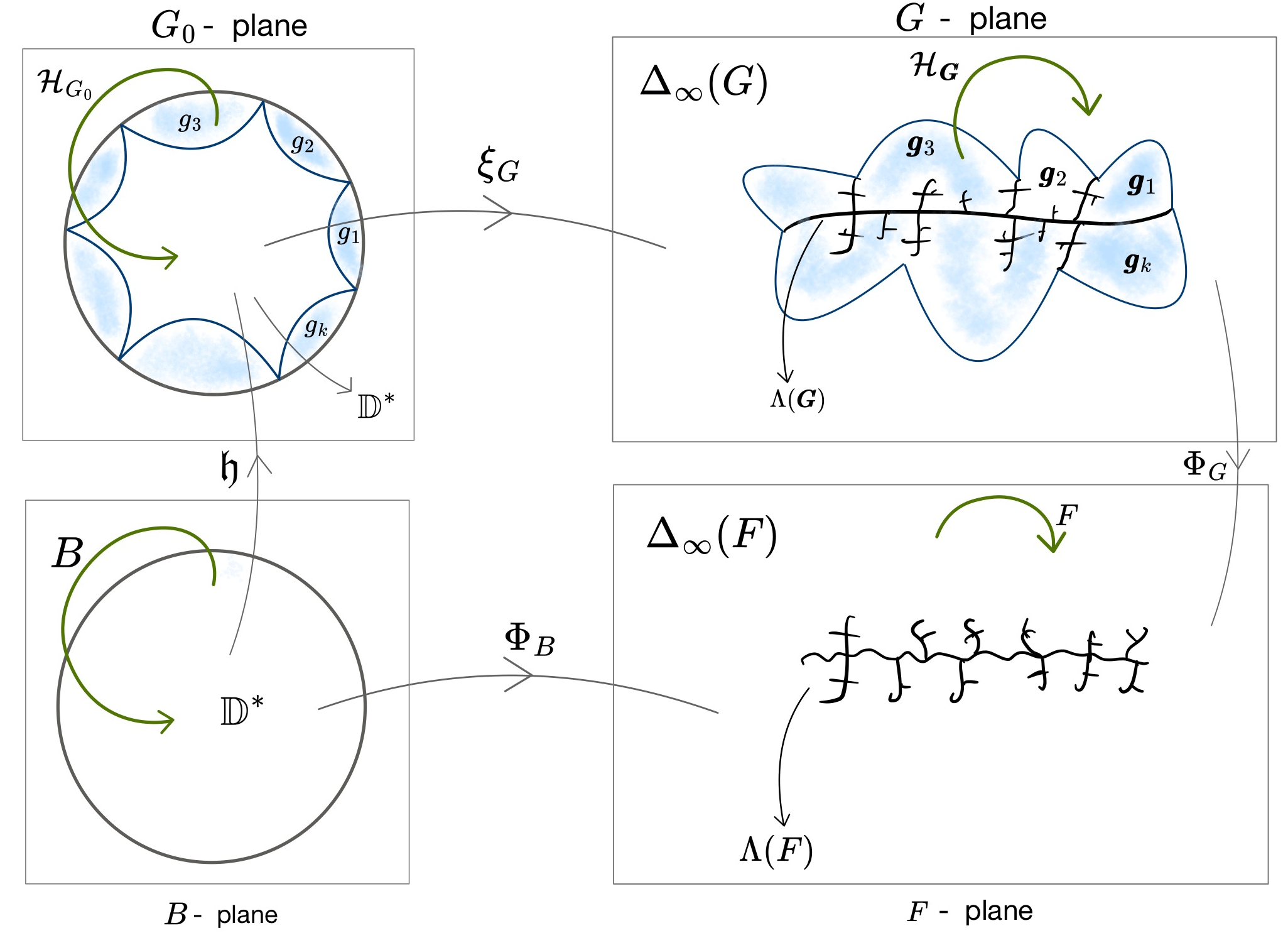}
    \caption{Illustrated is the construction of the multi--valued map $F$ as a conformal mating of a Blaschke product $B$ and the multi-valued accelerated Bowen--Series map $\cH_{G}$ of a Bers boundary group $G$. When $G$ is simply degenerate, the domain of definition $\cD_F$ is the entire sphere, and $F$ is single-valued on $\widehat{\C}-\Phi_G(\mathscr{T}_G)$, where $\Phi_G(\mathscr{T}_G)$ is a subtree of the limit set~$\Lambda(F)$.}
    \label{mvm_fig}
\end{figure}

According to the first part of Lemma~\ref{qc_david_conj_circle_lem}, there exist 
\begin{itemize}[leftmargin=*]
    \item a relative pinched neighborhood $N$ of $\bS^1$ in $\overline{\D^*}$, and
    \item a global quasiconformal map $\mathfrak{h}$ preserving $\D^*$,
\end{itemize} 
such that $\mathfrak{h}$ conjugates the restriction $B\vert_N$ to $\cH_{G_0}\vert_{\mathfrak{h}(N)}$.
We define a multi--valued map
\begin{align*}
&\widetilde{F}:\Delta_\infty(G)\bigcup\cD_{G}\longrightarrow\widehat{\C}\\
&\widetilde{F}=\begin{cases}  \big(\xi_{G}\circ\mathfrak{h}\big)\circ B\circ\big(\xi_{G}\circ\mathfrak{h}\big)^{-1}\quad \mathrm{on}\ \Delta_\infty(G),\\
    \cH_{G}  \hspace{3.92cm} \mathrm{on}\ \cD_{G}\cap\mathscr{K}(G).
\end{cases}
\end{align*}
The equivariance properties of $\mathfrak{h}$ and $\xi_{G}$ implies that $\widetilde{F}$ agrees with $\cH_{G}$ on $\xi_{G}(\mathfrak{h}(N))\cap\Delta_\infty(G)$. 
Hence, $\widetilde{F}$ acts as the single-valued holomorphic map $\cH_{G}$ on $\Int{\left(\Lambda(G)\cup\xi_{G}(\mathfrak{h}(N))\right)}-\RT_G$ (see Lemma~\ref{supp_mv_abs_lem}), and hence is quasiregular there. Further, $\widetilde{F}$ is quasiregular on $\Delta_\infty(G)$, since $\big(\xi_{G}\circ\mathfrak{h}\big)\circ B\circ\big(\xi_{G}\circ\mathfrak{h}\big)^{-1}$ is so. As finitely many points are quasiconformally removable, it follows that $\widetilde{F}$ is quasiregular outside $\RT_G$.

Now consider the Beltrami coefficient $\mu$ on the $G$-plane that is trivial on the filled limit set of $G$, and equal to $\big(\xi_{G}\circ\mathfrak{h}\big)_*\big(\mu_0\vert_{\D^*})\big)$ on $\Delta_\infty(G)$, where $\mu_0$ is the trivial Beltrami coefficient. Clearly, $\mu$ is $\widetilde{F}$-invariant.
Let $\Phi_{G}$ be a quasiconformal map straightening the Beltrami coefficient $\mu$. We define 
$$
F:= \Phi_{G}\circ\widetilde{F}\circ\Phi_{G}^{-1}:\cD_F:=\Phi_G\left(\Delta_\infty(G)\cup\cD_G\right)\to\widehat{\C}.
$$
As $\widetilde{F}$ is a quasiregular outside $\RT_G$, it follows by the Weyl's lemma that $F$ is a single--valued holomorphic map outside $\Phi_{G}(\RT_G)\subseteq \Phi_{G}(\Lambda(G))$, wherever defined. We also define 
$$
\Delta_\infty(F):=\Phi_{G}(\Delta_\infty(G)),\ \Lambda(F):=\Phi_{G}(\Lambda(G)),\ \mathscr{K}(F):=\Phi_{G}(\mathscr{K}(G))
$$
Further, the quasiconformal map 
$$
\Phi_B:=\Phi_{G}\circ\xi_{G}\circ\mathfrak{h}:\D^*\to\Delta_\infty(F)
$$
preserves the standard complex structure, and hence is conformal by the Weyl's lemma. By construction, $\Phi_B$ conjugates $B\vert_{\D^*}$ to $F\vert_{\Delta_\infty(F)}$.  

Similarly, as $\mu=0$ on $\mathscr{K}(G)$, the map $\Phi_{G}$ is conformal on $\mathscr{K}(G)$. Consequently, when $\Int{\mathscr{K}(G)}\neq\emptyset$, the map $\Phi_{G}$ conformally conjugates $\cH_{G}$ to the single--valued holomorphic map $F$ on the interior of $\mathscr{K}(F)$ (wherever defined). On the other hand, as $\widetilde{F}$ agrees with $\cH_{G}$ on $\xi_{G}(\mathfrak{h}(N))\cap\Delta_\infty(G)$, it follows that $F$ agrees with $\Phi_{G}\circ\cH_{G}\circ\Phi_{G}^{-1}$ on $\Phi_{G}\big(\xi_{G}\big(\mathfrak{h}\left(N\right)\big)\big)\cap\Delta_\infty(F)=\Phi_B(N\cap\D^*)$. Thus, $\Phi_{G}$ conjugates $\cH_{G}$, restricted to a pinched neighborhood of $\mathscr{K}(G)$ (wherever defined), to $F$, restricted to a pinched neighborhood of $\mathscr{K}(F)$ (wherever defined).
\end{proof}

\begin{figure}
 \captionsetup{width=0.96\linewidth}
\includegraphics[width=0.98\linewidth]{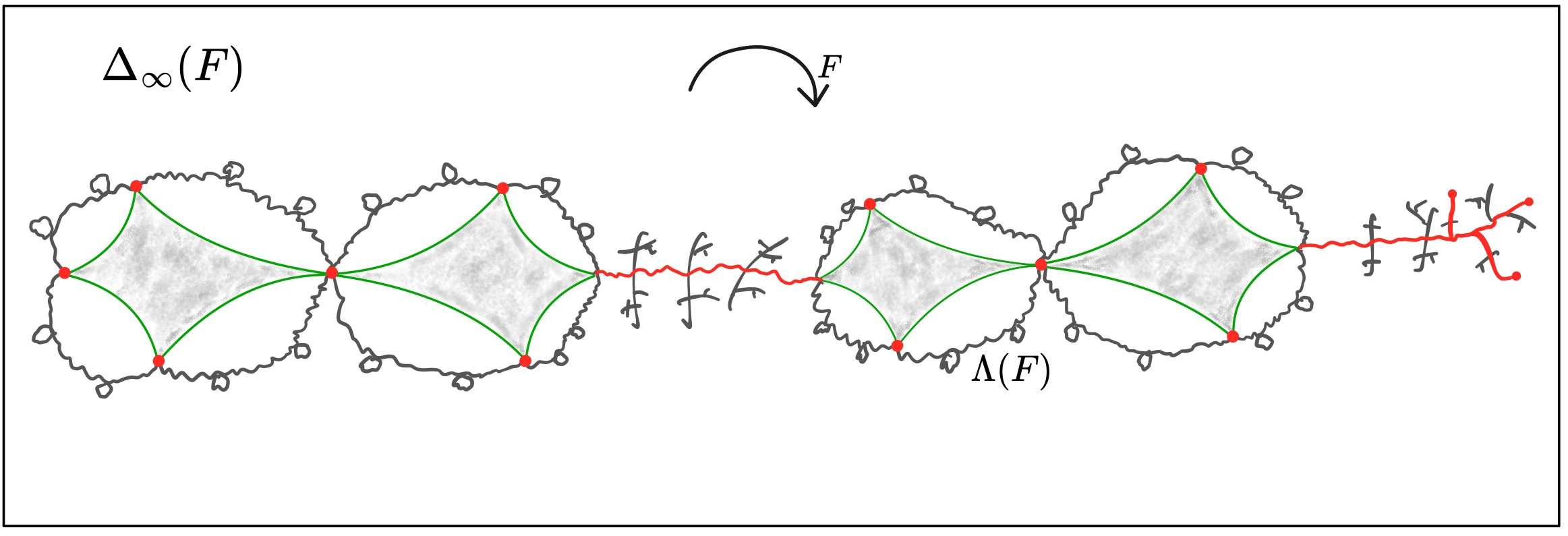}
    \caption{Depicted is the dynamical plane of a multi--valued mating $F$ between a Blaschke product $B$ and the multi-valued accelerated Bowen--Series map $\cH_{G}$ of a Bers boundary group $G$ with accidental parabolics. The set $\Phi_G(\mathscr{T}_G)$, that is a union of trees and finitely many points, is shown in red. The limit set $\Lambda(F)$ is the union of the black and the red sets. The shaded grey region is the complement of the domain of definition $\cD_F$. The interior of $\cD_F-\Phi_G(\mathscr{T}_G)$ is a simply connected domain.}
    \label{mv_mating_dom_fig}
\end{figure}
\begin{remark}\label{Dom_F_struct_rem}
 If $G$ is simply degenerate, then the domain of definition $\cD_F$ of the multi-valued map $F$ is the entire Riemann sphere $\widehat{\C}$. In this case, $\Phi_G(\mathscr{T}_G)$ is a sub-tree of the dendrite $\Lambda(G)$. In particular, $\cD_F-\Phi_G(\mathscr{T}_G)$ is a simply connected domain (see Figure~\ref{mvm_fig}). 
 
 On the other hand, if $G$ has accidental parabolics, then $\cD_F$ is a closed subset of $\widehat{\C}$ whose complement is a disjoint union of Jordan disks. Specifically, these Jordan disks are the interiors of ideal polygons contained in $\Phi_G\left(\Omega(G)-\Delta_\infty(G)\right)$. In this case, $\Phi_G(\mathscr{T}_G)$ is a forest contained in $\Lambda(G)$ (some of whose components may be points). Further, $\Int{\left(\cD_F-\Phi_G(\mathscr{T}_G)\right)}$ is a simply connected domain (see Figure~\ref{mv_mating_dom_fig}).
\end{remark}

\begin{remark}\label{unif_dil_rem}
1) Note that the construction of Lemma~\ref{bers_bdry_group_blaschke_mv_mating_lem} applies to Bers slice groups $G\in\cB(G_0)$ as well yielding single-valued conformal matings $F=B\mate\cH_G$. In fact, this procedure gives the same conformal mating as the one produced in Theorem~\ref{uniform_latice_infused_mating_thm}, where the polynomial and the Fuchsian group are canonically determined by the Blaschke product and the Bers slice group, respectively. Further, the dilatation of the global quasiconformal map $\Phi_{G}$ of Lemma~\ref{bers_bdry_group_blaschke_mv_mating_lem} depends only on $G_0$ and $B$, but not on $G\in\partial\cB(G_0)$ (cf. Remark~\ref{dilatation_rem}). 
\smallskip

\noindent 2) By the proof of Lemma~\ref{bers_bdry_group_blaschke_mv_mating_lem}, the multi--valued map $\widetilde{F}$ acts as a single group element on each interior component (i.e., component of the interior) of $\xi_{G}(\mathfrak{h}(N))-\RT_G$. Consequently, $F$ acts as $\Phi_{G}\circ \pmb{g}\circ\Phi_{G}^{-1}$, for some $\pmb{g}\in G$, on each interior component of $\Phi_{G}\left(\xi_{G}(\mathfrak{h}(N))-\RT_G\right)=\Phi_B(N)-\Phi_{G}(\RT_G)$. 
\end{remark}

\subsubsection{Promotion to holomorphic dynamical systems}
We now promote the above multi--valued mating $F=B\mate\cH_{G}$ to a dynamical system that captures the complete action of the group $G$. While this result is similar to Theorem~\ref{uniform_lattice_holo_dyn_syst_thm} in spirit, there is an important distinction between the two setups. Specifically, when $G$ is a Fuchsian group and $F$ is the conformal mating between $\cH_{G}$ and a polynomial $P$, one can reconstruct a conformal copy of the group $G$ from analytic continuations of the branches of $F$ on $\cE(F)$. However, for a Bers boundary group $G\in\partial\cB(G_0)$ without accidental parabolics, this is no longer possible as one of the components of the domain of discontinuity of $G_0$ completely degenerates in the limit. 

Note that each $\pmb{g}\in G$ corresponds to a unique element $g\in G_0$ (via a faithful representation $G_0\to G$). For $g\in G_0$, let $\mathscr{J}_g$ be the open, dense (full-measure) subset  of $\bS^1$ constructed in Proposition~\ref{countable_presentation_prop}. Then, there exists a relative neighborhood $\widehat{\mathscr{J}}_g$ of $\mathscr{J}_g$ in $\overline{\D^*}$, on which $g$ can be expressed as compositions of forward/backward branches of $\cH_{G_0}$. Also recall the global quasiconformal map $\Phi_{G}$ of Lemma~\ref{bers_bdry_group_blaschke_mv_mating_lem} that is conformal on $\mathscr{K}(G)$. 

We now define a set $\mathscr{N}_{\pmb{g}}$ that will be useful in the statement of the main theorem of this subsection. If $G$ is simply degenerate (in which case $\Int{\mathscr{K}(G)}=\emptyset$; i.e., $\mathscr{K}(G)=\Lambda(G)$), we set 
$$
\mathscr{N}_{\pmb{g}}:=\Lambda(G)\cup \xi_G(\widehat{\mathscr{J}}_g).
$$
Otherwise, $\Int{\mathscr{K}(G)}\neq\emptyset$; let $\{\mathscr{U}_\beta\}$ be the collection of the components $\mathscr{U}$ of $\Int{\mathscr{K}(G)}$ such that $\partial\mathscr{U}\cap\xi_G(\widehat{\mathscr{J}}_g)$ contains an arc. We define
$$
\mathscr{N}_{\pmb{g}}:=\Lambda(G)\cup \xi_G(\widehat{\mathscr{J}}_g)\cup \bigcup_\beta\mathscr{U}_\beta.
$$

\begin{remark}
    We do not know whether $\displaystyle\bigcup_\beta\mathscr{U}_\beta$ is equal to $\Int{\mathscr{K}(G)}$ when $G$ has accidental parabolics, but we expect it to be the case.
\end{remark}

\begin{theorem}\label{bers_bdry_group_blaschke_hds_mating_thm}
{Let $G_0$ be a uniform Fuchsian lattice with signature $(g,n) \notin \{(0,3), (1,1)\}$.} Let $G\in\partial\mathcal{B}(G_0)$, let $B$ be an expanding degree $d$ Blaschke product, and $\mathscr{B}$ be the holomorphic dynamical system consisting of the graphs of $B^{-m}\circ B^n$, $m,n\geq 0$. Then, there exists a holomorphic  dynamical system $\mathscr{F}$ such that $\mathscr{F}$ is a quasiconformally $\ G$-infused mating of $\mathscr{B}$ with the group $G$ in the following sense. 
\begin{enumerate}[leftmargin=*]
 \item There exists an $\mathscr{F}$-invariant partition  $\widehat{\C}=\mathscr{K}(F)\bigsqcup\Delta_\infty(F)$.
    
\item There exists a conformal map $\Phi_B:\D^*\rightarrow\Delta_\infty(F)$ such that $\Phi_B$ conjugates the holomorphic dynamical system consisting of the graphs of $B\vert_{\D}^{-m}\circ B\vert_{\D}^n$, to $\mathscr{F}\vert_{\Delta_\infty(F)}$.

\item $\Phi_G$ carries the graph of each $\pmb{g}\in G$, restricted to each interior component of $\mathscr{N}_{\pmb{g}}$, to a graph in $\mathscr{F}$, restricted to the corresponding interior component of $\Phi_G(\mathscr{N}_{\pmb{g}})$.
Further, any graph $\mathrm{gr}(F\vert_{\Int\cD_F}^{-m}\circ F\vert_{\Int\cD_F}^n)$, restricted to a relative deleted pinched neighborhood of $\Lambda(F)$ in $\overline{\Delta_\infty(F)}$, arises this way.
\end{enumerate}
\end{theorem}
\begin{proof}
Let $F=B\mate\cH_{G}$ be the multi--valued conformal mating between $B$ and $\cH_{G}$ constructed in Lemma~\ref{bers_bdry_group_blaschke_mv_mating_lem}. 
We define $\mathscr{F}$ to be the dynamical system given by (normalizations of) the maximal analytic continuations of $\mathrm{gr}(F\vert_{\Delta_\infty(F)}^{-m}\circ F\vert_{\Delta_\infty(F)}^n)$, for $m,n\geq 0$.

Let $\pmb{g}\in G$, and $g$ be the corresponding element of $G_0$. By definition, $\Phi_G\circ\pmb{g}\circ\Phi_G^{-1}$ is a quasiconformal homeomorphism of the sphere with 
$$
\overline{\partial}\left(\Phi_G\circ\pmb{g}\circ\Phi_G^{-1}\right)=0\quad \mathrm{on}\quad \mathscr{K}(F).
$$
Further, $\pmb{g}$ can be expressed as $\cH_{G}\vert_{\Delta_\infty(G)}^{-m}\circ \cH_{G}\vert_{\Delta_\infty(G)}^n$ on each component of $\xi_{G}(\widehat{\mathscr{J}}_g-\mathscr{J}_g)\subset\Delta_\infty(G)$ (with $m,n$ depending on the component). After possibly shrinking $\widehat{\mathscr{J}}_g$, we may assume that these branches of $\cH_{G}^{-m}\circ \cH_{G}^n$ agree with branches of $\widetilde{F}^{-m}\circ\widetilde{F}^n$ (see Lemma~\ref{bers_bdry_group_blaschke_mv_mating_lem} for the definition of $\widetilde{F}$).
Hence, the 
global quasiconformal map $\Phi_{G}\circ \pmb{g}\circ\Phi_{G}^{-1}$ can be written as local branches of $F\vert_{\Delta_\infty(F)}^{-m}\circ F\vert_{\Delta_\infty(F)}^n$ on $\Phi_{G}\left(\xi_{G}(\widehat{\mathscr{J}}_g-\mathscr{J}_g)\right)$.
In particular, the restriction of $\Phi_{G}\circ \pmb{g}\circ\Phi_{G}^{-1}$ on $\Phi_{G}\left(\xi_{G}(\widehat{\mathscr{J}}_g-\mathscr{J}_g)\right)$ is conformal. It now follows by the Weyl's lemma  that $\Phi_{G}\circ \pmb{g}\circ\Phi_{G}^{-1}$ is conformal on each interior component of $\mathscr{K}(F)\cup\Phi_G(\xi_G(\widehat{\mathscr{J}}_g))$. 

If $G$ is simply degenerate, then $\mathscr{K}(F)=\Lambda(F)$, and hence each interior component of $\Phi_G(\mathscr{N}_{\pmb{g}})=\Lambda(F)\cup\Phi_G(\xi_G(\widehat{\mathscr{J}}_g))$ intersects $\Phi_{G}\left(\xi_{G}(\widehat{\mathscr{J}}_g-\mathscr{J}_g)\right)$ non-trivially. On the other hand, if $G$ has accidental parabolics, then by definition, each $\mathscr{U}_\beta$ is a component of $\Int{\mathscr{K}(G)}$, and $\partial\mathscr{U}_\beta\cap\xi_G(\widehat{\mathscr{J}}_g)$ contains an arc. Hence, an interior component of $\Phi_G(\mathscr{N}_{\pmb{g}})=\Lambda(F)\cup\Phi_G\left(\xi_G(\widehat{\mathscr{J}}_g)\cup_\beta \mathscr{U}_\beta\right)$ intersects $\Phi_{G}\left(\xi_{G}(\widehat{\mathscr{J}}_g-\mathscr{J}_g)\right)$ non-trivially. It follows that in both cases, the graph of the restriction of $\Phi_{G}\circ \pmb{g}\circ\Phi_{G}^{-1}$ on each interior component of $\Phi_G(\mathscr{N}_{\pmb{g}})$ is an analytic continuation of the graph of a branch of $F\vert_{\Delta_\infty(F)}^{-m}\circ F\vert_{\Delta_\infty(F)}^n$, and thus it belongs to $\mathscr{F}$. The above discussion also implies that any analytic continuation of such a conformal restriction of $\Phi_{G}\circ \pmb{g}\circ\Phi_{G}^{-1}$ is given by $F\vert_{\Delta_\infty(F)}^{-m}\circ F\vert_{\Delta_\infty(F)}^n$, and it preserves $\mathscr{K}(F)$ and $\Delta_\infty(F)$.

On the other hand, $\widetilde{F}\vert_{\Delta_\infty(G)}^{-m}\circ\widetilde{F}\vert_{\Delta_\infty(G)}^n$ agrees with $\cH_G^{-m}\circ\cH_G^n$ on a relative deleted pinched neighborhood of $\Lambda(G)$ in $\overline{\Delta_\infty(G)}$. As $\cH_G$ acts piecewise via the elements of $G$, the branches of $F\vert_{\Delta_\infty(F)}^{-m}\circ F\vert_{\Delta_\infty(F)}^n$ act as $\Phi_G\circ g\circ\Phi_G^{-1}$, $g\in G$, on relative deleted pinched neighborhoods of $\Lambda(F)$ in $\overline{\Delta_\infty(F)}$. This implies the desired properties of the dynamical system $\mathscr{F}$.
\end{proof}

\begin{remark}
i) We will denote the dynamical system $\mathscr{F}$ of Theorem~\ref{bers_bdry_group_blaschke_hds_mating_thm} (associated with the conformal mating $F=B\mate\cH_{G}$) by~$\mathscr{B}\mate G$.

ii) In general, the branches of $F\vert_{\Delta_\infty(F)}^{-m}\circ F\vert_{\Delta_\infty(F)}^n$  do not admit analytic continuations around the points of $\Phi_G(\xi_G(\bS^1-\mathscr{J}_g))$. So these points may not arise as the first coordinates of the analytic subsets of $\mathscr{F}$. However, as $\mathscr{J}_g$ is dense in $\bS^1$ and $\xi_{G}:\overline{\D^*}\to\widehat{\C}$ is continuous, we have that $\Lambda(F)$ is contained in the closure of $\Phi_{G}\left(\xi_{G}(\widehat{\mathscr{J}}_g-\mathscr{J}_g)\right)$. It follows that the graphs of $\Phi_{G}\circ \pmb{g}\circ\Phi_{G}^{-1}\big\vert_{\Lambda(F)}$ are contained in $\overline{\mathscr{F}}$, and $\Phi_G$ carries the graphs of $g\vert_{\Lambda(G)}$ into $\overline{\mathscr{F}}$ (cf. Remark~\ref{closure_holo_dyn_syst_rem}).
\end{remark}

\subsection{Continuity of matings on the Bers boundary}\label{mating_cont_bers_bdry_subsec}

We now analyze continuity of matings of Bers closure groups with a fixed Blaschke product as the groups go from the interior to the boundary of the Bers slice. This will shed partial light on the compactification of the horizontal Bers slice $\mathfrak{B}(B)$ introduced in \S \ref{compactification_sec}.

\subsubsection{Background on Cannon-Thurston maps}\label{CT_subsubsec}
The following result, which is a special case of \cite[Theorem A]{MS17} and \cite[Theorem 4.1]{Mj-cmotions} will play an important role in the analysis below.

\begin{theorem}\label{thm-strong=unif}  
Let $\Gamma$ be a  closed surface group. Let $\rho_n: \Gamma  \to G_n$ be a sequence of quasi-Fuchsian representations converging strongly to a   Kleinian surface group $G_{\infty} = \rho_{\infty} (\Gamma )$. Then the sequence of Cannon-Thurston maps
$\xi_n: \Lambda(\Gamma) \to  \Lambda(G_n)$ converges uniformly to $\xi_{\infty}: \Lambda(\Gamma) \to  \Lambda(G_{\infty})$.
\end{theorem}

In the current setup (see Equation~\eqref{ct_eqn}), Theorem~\ref{thm-strong=unif} is the same as the uniform convergence of $\xi_n:\overline{\D^*}\to\overline{\Delta_\infty(G_n)}$ to $\xi_\infty:\overline{\D^*}\to\overline{\Delta_\infty(G_\infty)}$  
(under the assumption of strong convergence) viewed on the Riemann sphere.

In fact, the proofs of \cite[Theorem A]{MS17} and \cite[Theorem 4.1]{Mj-cmotions} furnish a bit more in our context. Let 
$\mathcal{L}$ denote the ending lamination for $G_{\infty}$ when the latter is simply degenerate. More generally, if $\Delta_\infty(G_{\infty})$ is the only $G_{\infty}$-invariant component of the domain of discontinuity, let $\mathcal{L}$ be the unique maximal unrealizable lamination given by the union of accidental parabolics and ending laminations. Let
$$
\partial \mathcal{L} = \{l^\pm : l \ \textrm{is a leaf of}\ \mathcal{L} \}
$$ 
denote the pairs of end-points of leaves of $\mathcal{L}$. Let $R_\mathcal{L}$ denote the equivalence relation on $\bS^1 = \partial \Hyp^2$ induced by $\mathcal{L}$ (by introducing all diagonals in complementary ideal polygons).

\begin{lem}\label{lem-dia0lamn}
Under the assumption of strong convergence,    $d_{\mathrm{sph}}\big(\xi_n(x),\xi_n (y)\big) \to 0$
 if and only if $x, y \in R_\mathcal{L}$, where $d_{\mathrm{sph}}$ denotes the spherical metric.
\end{lem}

\begin{proof}
The main theorem of \cite{mahan-elct} says that 
the Cannon-Thurston map for $G_{\infty}$ identifies $x, y \in \bS^1$ if and only if $x, y \in R_\mathcal{L}$. The conclusion now follows from
   Theorem~\ref{thm-strong=unif}. 
\end{proof}

Since $\infty\in\Delta_\infty(G_n)$ for
$n \in \mathbb{N} \cup \{\infty\}$, we can take a stereographic projection of the Riemann sphere to the plane $\mathbb C$, and assume that all the limit sets lie in $\mathbb C$. For $x, y \in  \Lambda(G_n)$, let $[x,y]^e_n$ denote the Euclidean geodesic in $\mathbb C$ joining $\xi_n(x), \xi_n(y)$. Then Lemma~\ref{lem-dia0lamn} says that the Euclidean length of $[x,y]^e_n$ converges to zero if and only if $x, y \in R_\mathcal{L}$. Since $G_n$ is quasi-Fuchsian it has two invariant components, the one containing $\infty$ is $\Delta_\infty(G_n)$, while the other one is denoted by $\Delta_0(G_n)$. 

We would like to show further, that the Euclidean geodesics $[x,y]^e_n$ in the previous paragraph can be replaced by hyperbolic geodesics in (the uniformization of) $\Delta_0(G_n)$.
To explain  the basic idea behind the proof, suppose that $\Delta_0(G_n)$ is convex in 
$\mathbb C$.
For $n \in \mathbb{N}$, let $\mathcal{L}_n^e$ denote the collection of Euclidean geodesics in $\Delta_0(G_n)$ joining $\xi_n(x), \xi_n(y)$ for all $(x, y) \in R_\mathcal{L}$. Then, by convexity of 
$\Delta_0(G_n)$, these Euclidean geodesics lie in 
$\Delta_0(G_n)$ and approximate the corresponding hyperbolic geodesics in the following sense:
Let $[x,y]^h_n$ denote the hyperbolic geodesic joining $x, y$ in $\Delta_0(G_n)$. Then the Hausdorff distance between $[x,y]^h_n$ and $[x,y]^e_n$
tends to zero as $n\to \infty$. It follows that the diameter of $[x,y]^h_n$ tends to zero as $n\to \infty$.

Next, for arbitrary $w, z\in\bS^1$ with $(w, z) \notin R_\mathcal{L}$, 
the Euclidean (respectively, hyperbolic) geodesic $[w,z]^e_n$ (respectively, $[w,z]^h_n$) is  transverse to $\mathcal{L}_n^e$ for all $n \in \mathbb{N}$. 

However, $\Delta_0(G_n)$ is typically not convex in 
$\mathbb C$. In general, even when $\Delta_0(G_n)$ is not convex, we can map $\mathcal{L}$ to $\Delta_0(G_n)$ via the uniformization map inherited from the simultaneous uniformization defining $G_n$
to obtain a hyperbolic realization $\mathcal{L}_n^h$. 
Further $\mathcal{L}_n^h$ has the property that for $(x, y) \in R_\mathcal{L}$, the leaf 
(or diagonal of complementary polygon) 
$[x,y ]^h_n$ in (the diagonal completion of) 
$\mathcal{L}_n^h$ has diameter going to zero as
$n\to \infty$.

\begin{remark}\label{diam_shrink_rem}
    The last statement needs some justification.
    Let $M_n = \Hyp^3/G_n$. Then, for $n \in \mathbb{N}$, $\mathcal{L}$ is realized as a geodesic lamination $\mathcal{L}_n^r$ in 
    the 3-manifold $M_n$.
    Choose base-points $o_n\in M_n$ for the  geometric convergence of $M_n$ to $M_\infty$.
    Then the proof of the main theorem of \cite{MS17,Mj-cmotions} shows that the distance of 
    $\mathcal{L}_n^r$ from $o_n$ tends to $\infty$ as $n \to \infty$. Hence for $(x, y) \in R_\mathcal{L}$,  the geodesic
$[x,y ]_n^h$ in (the diagonal completion of) 
$\mathcal{L}_n^h$ has diameter going to zero as
$n\to \infty$.

Note that the last step can also be obtained from the Gehring--Hayman theorem \cite[Theorem 4.2]{Pom} or \cite[Theorem 5.1]{OT10}. Indeed, we have already noted that $\xi_n:\overline{\D^*}\to\overline{\Delta_\infty(G_n)}$ converges uniformly to $\xi_\infty:\overline{\D^*}\to\overline{\Delta_\infty(G_\infty)}$  (under the assumption of strong convergence). Suppose $x, y \in R_\mathcal{L}$. Then by uniform convergence, for each $\varepsilon$, there exists $r> 1$ so that for all sufficiently large $n$, $\xi_n(e^{2\pi i x}), \xi_n(e^{2\pi i y}) \in B(\xi_\infty(e^{2\pi i x}), \varepsilon)$ and $B(\xi_\infty(e^{2\pi i x}), \varepsilon) \cap \xi_n(\D^*_r) = \emptyset$, where $\D^*_r = \{z:|z|> r\}$. Let $P_n: \C - \xi_n(\D^*_r) \rightarrow \C$ be a map that sends $\xi_n(se^{2\pi it})$ to $\xi_n(e^{2\pi it})$ for each $1\leq s \leq r$ and extends to the identity on  $\Delta_0(G_n) = \C - \xi_n(\overline{\D^*})$.
Let $\widetilde{\gamma}_n$ be the Euclidean line segment connecting $\xi_n(e^{2\pi i x}), \xi_n(e^{2\pi i y})$, and let $P_n(\widetilde{\gamma}_n)$ be its projection. 
Note that by uniform convergence, we may assume that the diameter of $S_n(\widetilde{\gamma}_n)$ is bounded by $2\varepsilon$.
By perturbing the points on $S_n(\widetilde{\gamma}_n)$ slightly, we obtain an arc $\gamma_n$ connecting $\xi_n(e^{2\pi i x}), \xi_n(e^{2\pi i y})$ whose interior points are contained in $\Delta_0(G_n)$ and whose diameter is bounded by $4\varepsilon$. By the Gehring--Hayman theorem \cite[Theorem 4.2]{Pom}, we conclude that the hyperbolic geodesic $[x,y]_n^h$ of $\Delta_0(G_n)$ connecting $\xi_n(e^{2\pi i x}), \xi_n(e^{2\pi i y})$ has diameter bounded by $K\varepsilon$ for some universal constant $K$. Since $\varepsilon$ is arbitrary, we conclude that the diameter of $[x,y]_n^h$ goes to zero.
\end{remark}

So far, we have considered the general case where the end-invariant is given as  follows.
Let $\Sigma$ denote the topological surface underlying $\Delta_\infty(G_\infty)/G_\infty $. There is a finite collection of (closed) essential subsurfaces with disjoint interiors $\Sigma_1, \cdots, \Sigma_k $
such that 
\begin{enumerate}[leftmargin=8mm]
\item a finite subcollection of $\Sigma_1, \cdots, \Sigma_k $ have conformal structures and correspond to geometrically finite subgroups of $G_\infty$,
\item the complementary subcollection consists of degenerate subsurfaces with ending laminations supported on them,
\item $\partial \Sigma_i$ consists of accidental parabolics,  
\item the ending lamination $\cL$ for $G_\infty$
consists of the union of the ending laminations for the degenerate subsurfaces and the union of the accidental parabolics,  
\item the end-invariant for $G_\infty$ consists of $\cL$ and the conformal structures on the 
finite subcollection in item (1).
\end{enumerate}
We specialize first to the principal case of interest where 
$G_\infty$ is simply degenerate, and the limit set $\Lambda(G_{\infty})$ is a dendrite.

Let $[w,z]_\infty$ be the unique arc in $\Lambda(G_{\infty})$ joining 
$\xi_\infty(w), \xi_\infty(z)$.

\begin{lem}\label{arc_in_dendrite_lem}
  Under the above assumptions,  $[w,z]^h_n$
  converges to $[w,z]_\infty$.
\end{lem}

\begin{proof}
    The main theorem of \cite{mahan-elct} guarantees that $\xi_\infty$ identifies a pair of points on $\bS^1$ if and only if the pair belongs to $R_\mathcal{L}$. As a map of the closed unit disk $\overline{\D}$ to the Riemann sphere, $\xi_\infty$ thus collapses leaves of 
    $\mathcal L$ and complementary ideal polygons to points.

   By Lemma~\ref{lem-dia0lamn} and Remark~\ref{diam_shrink_rem}, leaves of 
    $\mathcal{L}_n$ and complementary ideal polygons have diameter shrinking to zero
    as $n \to \infty$. 
     Recall that $[w,z]_n^h$ 
     is transverse to $\mathcal L$ 
    for all $n \in \mathbb{N}$.
    The Lemma now follows from Lemma~\ref{lem-dia0lamn}.
\end{proof}

    Lemma~\ref{arc_in_dendrite_lem} deals with one extreme case, where the end-invariant is a single ending lamination (necessarily with support all of $\Sigma$). Another extreme case is where all the essential subsurfaces $\Sigma_1, \cdots,\Sigma_k$ are geometrically finite and the ending lamination $\cL$ consists of the simple closed curves $\bigcup_i \partial \Sigma_i$, all corresponding to accidental parabolics. 
    In this case, Lemma~\ref{lem-dia0lamn} and the succeeding discussion till before Lemma~\ref{arc_in_dendrite_lem} show that $\cL$ can be embedded as $\cL_n$ in $\Delta_0(G_n)$ with diameter of each leaf going to zero. Further,  the Carath{\'e}odory limits of the Jordan domains $\Delta_0(G_n)$ is a countable collection of Jordan disks, and the closures of these disks are arranged via a tree. Note that this tree can be identified with the Bass-Serre tree $\mathbb T$ of the splitting of $\pi_1(\Sigma)$ along the cyclic groups corresponding to the boundary curves of the~$\Sigma_i$'s.

    The main theorem of \cite{MS13} provides the special case of Theorem~\ref{thm-strong=unif} in this case, as $G_\infty$ is geometrically finite.
    From this, it follows that the limiting tree 
   of Jordan disks (i.e., the closure of the union of the limiting Jordan disks) in the previous paragraph can be identified with the filled limit set 
    $\mathscr{K}(G_\infty)$.

   \begin{remark}\label{rmk-arcingeofinite}
    Now, let $w, z$  be a pair of arbitrary points in $\bS^1$ as before and $[w,z]_n^h$ be the hyperbolic geodesic in $\Delta_0(G_n)$ as in the discussion preceding Lemma~\ref{arc_in_dendrite_lem}. Let $w_\infty, z_\infty \in
    \Lambda(G_\infty)$ be the Cannon-Thurston images of $w, z$ for $G_\infty$. Then there is a unique (finite, semi-infinite or bi-infinite) geodesic in $\mathbb T$ encoding the sequence of Jordan disks crossed by a path between  
    $w_\infty, z_\infty$. Replace the intersection of this minimal path with any Jordan disk (in $\Int{\mathscr{K}(G_\infty)}$) by the hyperbolic geodesic in that Jordan disk. Let 
    $[w,z]_\infty^h$ be the resulting concatenation of hyperbolic geodesics. A reprise of the argument in Lemma~\ref{arc_in_dendrite_lem} shows that $[w,z]^h_n$ converges to $[w,z]^h_\infty$.
\end{remark}

\noindent {\bf The general case:}
     Finally, we deal with the hybrid case, where the end-invariant consists of more than one geometrically finite surface and/or ending laminations. We refer the reader to Figure~\ref{mv_mating_dom_fig}. We combine
the conclusions of Lemma~\ref{arc_in_dendrite_lem} and 
Remark~\ref{rmk-arcingeofinite}. We continue with the notation in the geometrically finite case, and start with the 
Bass-Serre tree $\mathbb T$. Construct the geodesic in $\mathbb T$ corresponding to the hyperbolic geodesic joining $w, z$. Call this geodesic $[w,z]_\mathbb{T}$. 

We first construct an intermediate auxiliary geometrically finite group $G_\infty'$ where all the
subsurfaces $\Sigma_i$'s are geometrically finite.
Let $[w,z]^{gf}_\infty$ be the corresponding piecewise geodesic path constructed in Remark~\ref{rmk-arcingeofinite} above. The $G_\infty'$-orbits of a finite collection
of Jordan domains in $\Int{\mathscr{K}(G_\infty')}$ correspond to subgroups that are degenerate in $G_\infty$. Let $\cD$ denote the union of $G_\infty'$-orbits of these finite collection of Jordan domains. 

The limit set of $G_\infty$ can be obtained from the 
limit set of $G_\infty'$ by degenerating the Jordan disks $D \in \cD$ to  dendrites
$\cT_D$ by collapsing  leaves (and complementary ideal polygons) of the
corresponding ending laminations (by \cite[Theorem 6.1]{mahan-survey}, or \cite[Theorem 5.6]{mahan-kl}, see also\cite[Theorem 5.3]{mahan-icm}). 

Let $D \in \cD$ be one of the disks that are collapsed, and let $[w,z]^{D}_\infty:=[w,z]^{gf}_\infty \cap D$ denote the bi-infinite hyperbolic geodesic in $D$ that occurs as one of the pieces of $[w,z]^{gf}$.
Replace $[w,z]^{D}_\infty$ by the unique arc $[w,z]^{\cT_D}_\infty$ in 
$\cT_D$ joining its end-points. Finally, let $[w,z]_\infty$ denote the concatenation of $[w,z]^{D}_\infty$'s (for geometrically finite surfaces in the end invariant) and $[w,z]^{\cT_D}_\infty$'s (for degenerate surfaces in the end invariant). The conclusion that $[w,z]^h_n$ converges to $[w,z]_\infty$ now follows  in this general case, too. 
  
\subsubsection{Convergence of matings}
Suppose that $G_n\in\cB(G_0)$ converge strongly to $G_\infty\in\partial\cB(G_0)$. We lay out the setup for the main convergence result of this section.
\begin{itemize}[leftmargin=6mm]
    \item For $n\in\N\cup\{\infty\}$, we denote by $\cH_n:=\cH_{G_n}$ the accelerated Bowen--Series map associated to $G_n$, compatible with an accelerated Bowen--Series map $\cH_{G_0}$ of the base group $G_0$. Note that $\cH_\infty$ is generally multi-valued, but single-valued outside $\mathscr{T}_{G_\infty}$ (cf. \S \ref{multi_valued_abs_bers_bdry_subsec}).
    \item Let $B$ be an expanding Blaschke product with $\deg{\left(B\right)}=\deg{\left(\cH_{G_0}\vert_{\bS^1}\right)}$. 
    \item The conformal mating between $\cH_n$ and $B$ is denoted by $F_n$, and the associated holomorphic dynamical system is denoted by $\mathscr{F}_n$, for $n\in\N\cup\{\infty\}$.
    \item The domain of definition of $F_n$ is denoted by $\cD_{F_n}$, $n\in\N\cup\{\infty\}$. For $n\in\N$, $\cD_{F_n}$ is a closed Jordan disk (see Theorem~\ref{uniform_latice_infused_mating_thm}). On the other hand, $\cD_{F_\infty}$ is the complement of a disjoint union of Jordan domains in $\widehat{\C}$ (see Lemma~\ref{bers_bdry_group_blaschke_mv_mating_lem} and Remark~\ref{Dom_F_struct_rem}).
    \item For $n\in\N\cup\{\infty\}$, the action of $F_n$ on $\Delta_\infty(F_n)$ is conformally conjugate to $B\vert_{\D^*}$. When $n\in\N$, we have $\Delta_\infty(F_n)=\Int{\cK(F_n)}$, where $\cK(F_n)$ is the filled Julia set of the degenerate polynomial-like map $F_n$ (cf. \S \ref{compactification_sec}). Further, $\Delta_\infty(F_\infty)=\Phi_{G_\infty}(\Delta_\infty(G_\infty))$ (cf. Lemma~\ref{bers_bdry_group_blaschke_mv_mating_lem}).
    \item $\mathfrak{B}(B)$ is the horizontal Bers slice through $B$ (see \S \ref{compactification_sec}).
\end{itemize}

Recall that the conformal mating between a Blaschke product $B$ and a Bers slice group $G_n\in\cB(G_0)$ can be constructed by quasiconformal surgery either by starting with the dynamical plane of $B$ (as in Theorem~\ref{uniform_latice_infused_mating_thm}) or by starting with the dynamical plane of $G_n$ (as in Lemma~\ref{bers_bdry_group_blaschke_mv_mating_lem}). To study continuity of the mating operation as $G_n\to G_\infty$, it will be convenient to appeal to the latter recipe, and this is the approach we will adopt in the proposition below.

\begin{prop}\label{mating_cont_bers_bdry_prop}
     With the above setup, the restriction $F_n\vert_{\Int{\cD_{F_n}}}$ converges to  $F_\infty\vert_{\Int{\left(\cD_{F_\infty}-\Phi_{G_\infty}(\mathscr{T}_{G_\infty})\right)}}$ in the Carath{\'e}odory topology. In particular, the holomorphic dynamical systems $(\widehat{\C},\mathscr{F}_n)\in\mathfrak{B}(B)$ converge to the dynamical system $(\widehat{\C},\mathscr{F}_\infty)\in\partial\mathfrak{B}(B)$.
\end{prop}
\begin{proof}
    Let us denote the Cannon--Thurston map associated with $G_n$ by $\xi_n$, $n\in\cup\{\infty\}$. By Theorem~\ref{thm-strong=unif}, we have uniform convergence of the 
maps $\xi_n:\overline{\D^*}\to\overline{\Delta_\infty(G_n)}$  to $\xi_\infty:\overline{\D^*}\to\overline{\Delta_\infty(G_\infty)}$. 

Let 
$$
\Pi_n:=\widehat{\C}-\left(\Delta_\infty(G_n)\cup\cD_{G_n}\right),
$$ 
and $v_{1,n}, \cdots , v_{k,n}$ denote the vertices of $\Pi_n$.
By Theorem~\ref{thm-strong=unif}, $v_{j,n} \to 
v_{j,\infty}$, for points $v_{j,\infty}\in \Lambda({G_\infty})$, $j\in\{1, \cdots, k\}$. Note that the points $v_{j,\infty}$ are distinct by Lemma~\ref{lem-dia0lamn}. By Lemma~\ref{arc_in_dendrite_lem} and the discussion following it, we conclude that $\overline{\Pi}_n$ converges in the Hausdorff topology to the closure of $\left(\widehat{\C}-\cD_{G_\infty}\right)\cup\mathscr{T}_{G_\infty}$.

Let $\widetilde{F}_n$ be the quasiregular map appearing in the construction of $F_n=B\mate\cH_n$, $n\in\N\cup\{\infty\}$ (see the proof of Lemma~\ref{bers_bdry_group_blaschke_mv_mating_lem}).
The discussion of the previous two paragraphs implies that $\widetilde{F}_n$, restricted to $\Int{\mathrm{Dom}(\widetilde{F}_n)}$, converges to $\widetilde{F}_\infty$, restricted to $\Int{\left(\mathrm{Dom}(\widetilde{F}_\infty)-\mathscr{T}_{G_\infty}\right)}$, in the Carath{\'e}odory topology.

To complete the proof of the proposition, it now suffices to show that the quasiconformal map $\Phi_{G_n}$ hybrid conjugating $\widetilde{F}_n$ to $F_n$ converges to the quasiconformal map $\Phi_{G_\infty}$ hybrid conjugating $\widetilde{F}_\infty$ to $F_\infty$.

Let $\mu_n$ be the Beltrami coefficient constructed in Lemma~\ref{bers_bdry_group_blaschke_mv_mating_lem} associated with $G_n$ (in Figure~\ref{mvm_fig}, the top right corner shows the $G_n$-plane). Specifically, $\mu_n$~is
\begin{enumerate}[leftmargin=*]
    \item trivial on the filled limit set
 $\mathscr{K}(G_n)$, 
 and 
\item equal to $(\xi_n)_*(\mathfrak{h}_*(\mu_0\vert_{\D^*}))$ on $\Delta_\infty(G_n)$ (where $\mu_0$ is the trivial Beltrami coefficient).
\end{enumerate}
Finally, to emphasize the dependence on $n$, we denote the conformal conjugacy between $B\vert_{\D^*}$ and $F_n\vert_{\Delta_\infty(F_n)}$ by $\Phi_{B,n}=\Phi_{G_n}\circ\xi_n\circ\mathfrak{h}$.
After appropriate normalizations and after passing to a subsequence, we can assume that $\Phi_{B,n}$ converges locally uniformly to a conformal map $\widehat{\Phi}_{B,\infty}:\D^*\hookrightarrow\widehat{\C}$.

The maps $\Phi_{G_n}$ are $K$-quasiconformal, for some fixed $K$ (see Remark~\ref{unif_dil_rem}). After possibly passing to a further subsequence, we can assume that $\Phi_{G_n}$ converges uniformly to a $K$-quasiconformal map $\Phi_\infty$. By the uniform convergence of $\xi_{n}$ to $\xi_\infty$ on $\overline{\D^*}$, the map $\Phi_\infty$ is conformal on $\Int{\mathscr{K}(G_\infty})$, and the map $\Phi_\infty\circ\xi_\infty\circ\mathfrak{h}$ is conformal on $\D^*$. Further, by the Ahlfors' measure zero conjecture proven by Brock--Canary--Minsky in \cite{BCM12}, the limit set $\Lambda(G_\infty)$ has zero area. Hence, 
\begin{enumerate}[leftmargin=8mm]
    \item $\overline{\partial}\Phi_\infty=0$ on $\mathscr{K}(G_\infty)$, and
    \item $\Phi_\infty$ sends the complex structure defined by the Beltrami coefficient $\big(\xi_{\infty}\circ\mathfrak{h}\big)_*\big(\mu_0\vert_{\D^*})\big)$ to the standard complex structure.
\end{enumerate}
By the proof of Lemma~\ref{bers_bdry_group_blaschke_mv_mating_lem}, the quasiconformal map $\Phi_{G_\infty}$ which conjugates $\cH_{G_\infty}$ to the conformal mating $F_\infty=B\mate \cH_\infty$ has the same Beltrami coefficient. Hence, up to a M{\"o}bius change of coordinates, we have $\Phi_\infty=\Phi_{G_\infty}$.
\end{proof}

By the Bers Density Conjecture established in \cite{BCM12}, each $G_\infty\in\partial\cB(G_0)$ is the strong limit of a sequence $\{G_n\}\subset\cB(G_0)$. In fact, according to \cite{AC96}, if $\cB(G_0)\ni G_n\rightarrow G_\infty$ algebraically, and $G_\infty\in\partial\cB(G_0)$ has no parabolics, then the convergence is strong. Thus, Proposition~\ref{mating_cont_bers_bdry_prop} implies that $\partial\mathfrak{B}(B)$ contains the dynamical systems $\mathscr{F}=\mathscr{B}\mate G$, for all $G\in\partial\cB(G_0)$. However, we do not know if 
$$
\mathfrak{B}(B)-\big\{\mathscr{F}=\mathscr{B}\mate G:G\in\partial\cB(G_0)\big\}
$$
is empty.

\section{Dependence of matings and their parameter spaces on surface topology}\label{compare_sec}

Both the dynamical structure and the parameter space of the degenerate polynomial-like maps (or holomorphic dynamical systems) arising as matings of polynomials and Fuchsian lattices depend crucially on the topology of the surface.

\subsection*{The case of punctured spheres}

Holomorphic dynamical systems realizing matings between polynomials and Fuchsian punctured sphere groups were investigated in \cite{MM1,MM2,LLM24}. A special feature of Fuchsian groups uniformizing $S_{0,k}$ is that the associated Bowen--Series map of \S \ref{ss:idealcorpolygon} is already a (continuous) degree $2k-3$ expansive $C^1$ covering maps of the circle; i.e., no acceleration (cf. \S \ref{acceleration_sec}) is necessary to mate such Bowen--Series maps with polynomials. Further, the degenerate polynomial-like maps arising as matings of polynomials and these Bowen--Series maps are algebraic functions: the so-called \emph{B-involutions} \cite{LLM24}. Hence, the associated holomorphic dynamical systems $\mathscr{F}$ are realized as iterated algebraic correspondences. 

However, this mating construction involves David surgery (due to the lack of quasisymmetric compatibility between expanding power maps and Bowen--Series maps with parabolics), and this yields matings of Bowen--Series maps of Fuchsian punctured sphere groups only with subhyperbolic / geometrically finite polynomials. The fact that the matings are algebraic comes to partial rescue. Specifically, algebraic correspondences of a given bi-degree serve as an ambient parameter space for these matings. A \emph{vertical fiber}; i.e., degenerate polynomial-like maps having the Bowen--Series map of a fixed Fuchsian punctured sphere group as their external map, is a compact subset of this ambient parameter space \cite{LMM26a}. This compactness property, combined with combinatorial continuity arguments, permits one to construct matings between Bowen--Series maps associated with $S_{0,k}$ and at most finitely renormalizable degree $2k-3$ polynomials with no non-repelling cycles. These matings are realized both as degenerate polynomial-like maps and as algebraic correspondences; see \cite{LLM24}. 

The algebraicity of these matings also allows one to obtain sharper information on the boundary of a horizontal Bers slice $\mathfrak{B}(B)$ in this space of matings (i.e., matings between a fixed expanding degree $2k-3$ Blaschke product $B$ and all Fuchsian groups uniformizing $S_{0,k}$; cf. \S \ref{compactification_sec}). According to \cite{LMM26a}, the holomorphic dynamical systems lying on the boundary $\partial\mathfrak{B}(B)$ are generated by algebraic correspondences of bi-degree $(2k-3,2k-3)$ on $\widehat{\C}$ (in other words, $\mathfrak{B}(B)$ is pre-compact in the space of  bi-degree $(2k-3,2k-3)$ algebraic correspondences). However, unlike Theorem~\ref{bers_bdry_group_blaschke_hds_mating_thm}, quasiconformal surgery techniques cannot be used to produce matings of $B$ with groups on the Bers boundary of $\mathrm{Teich}(S_{0,k})$ (due to the existence of parabolics). We conjecture the following:

\begin{conj}
    Correspondences lying on the boundary of the above horizontal Bers slice are precisely the matings of $B$ with Kleinian groups lying on $\partial_{\mathrm{Bers}}\mathrm{Teich}(S_{0,k})$.
\end{conj}

\subsection*{The case of closed surfaces}

In a certain sense, closed surfaces and punctured spheres lie on  opposite ends of the spectrum in our mating framework. The degenerate polynomial-like maps having accelerated Bowen--Series maps of uniform Fuchsian lattices as their external maps are not algebraic, and the associated holomorphic dynamical systems are realized as iterated transcendental correspondences (see Propositions~\ref{non_alg_prop_uniform} and ~\ref{alg_is_impossible_prop}).
Further, we do not have an explicit description for these matings or their ambient parameter spaces. 

Nonetheless, these accelerated Bowen--Series maps are quasisymmetrically conjugate to power maps, allowing us to mate them with all polynomials with connected Julia set (see Theorem~\ref{uniform_latice_infused_mating_thm} and Theorem~\ref{uniform_lattice_holo_dyn_syst_thm}). In fact, there is a natural bijection between vertical fibers and the connectedness loci of polynomials of appropriate degree.
(Here, the vertical fiber consists of degenerate polynomial-like maps having the accelerated Bowen--Series map of a fixed Fuchsian closed surface group as the external map.)

The same quasisymmetric compatibility is also responsible for the construction of dynamical systems as matings between Blaschke products $B$ and Kleinian groups $G\in\partial_{\mathrm{Bers}}\mathrm{Teich}(S_{g,0})$ (see  Theorem~\ref{bers_bdry_group_blaschke_hds_mating_thm}). By Proposition~\ref{mating_cont_bers_bdry_prop}, these dynamical systems lie on the boundary of the corresponding horizontal Bers slice $\mathfrak{B}(B)$; but we do not know if the boundary $\partial\mathfrak{B}(B)$ contains other (more exotic) dynamical systems.

\subsection*{The case of positive genus surfaces with punctures}

This is precisely the setting in which both the algebraicity of the mated holomorphic dynamical systems and the quasisymmetric compatibility between accelerated Bowen–Series maps and power maps break down. Thus, we are only able to construct degenerate polynomial-like maps and holomorphic dynamical systems as matings of such Fuchsian groups with subhyperbolic/geometrically finite polynomials (see Theorem~\ref{non_uniform_latice_mating_thm} and Theorem~\ref{non_uniform_lattice_holo_dyn_syst_thm}). Moreover, no precise information is available on the boundary of a horizontal Bers slice $\mathfrak{B}(B)$ in this setting, Specifically, we neither know whether $\partial\mathfrak{B}(B)$ contains matings between the Blaschke product $B$ and Kleinian groups lying on $\partial_{\mathrm{Bers}}\mathrm{Teich}(S_{g,n})$,~$g,n>0$, nor  whether the dynamical systems on $\partial\mathfrak{B}(B)$ are globalizable.

\appendix

\section{Orbit equivalence for co-compact orbifold groups}\label{sec:oeorbifold}
In this appendix, we prove orbit equivalence properties of pull-back Bowen--Series maps for general uniform lattices.

\begin{prop}[=Proposition~\ref{prop:nonclosedsurfaceOESt}]\label{prop:nonclosedsurfaceOE}
    Let $G$ be a uniform lattice for a compact hyperbolic surface of type $(g, n; \nu_1,\ldots,\nu_n)$ where $n \geq 1$ and 
    $$
    (g,n) \notin\{(0,3), (1,1)\}.
    $$
    There is a pull-back Bowen--Series Markov map $\mathcal{F}:\mathbb{S}^1 \longrightarrow \mathbb{S}^1$ and a unique exceptional $G$-orbit \[
  E\subset\mathbb S^1\setminus\mathrm{bk}(\mathcal{F})
\]
with the following properties.
\begin{enumerate}
\item If $x\in\mathbb S^1\setminus(\mathrm{bk}(\mathcal{F})\cup E)$, then
\[
  \mathcal O_G(x)=\operatorname{GO}_\mathcal{F}(x).
\]
\item There are points $x_v,x_s\in E$ in distinct $\mathcal{F}$-grand orbits such
that
\[
  E=\operatorname{GO}_\mathcal{F}(x_v)\cup\operatorname{GO}_\mathcal{F}(x_s).
\]
\end{enumerate}
In particular, $\mathcal{F}$ is {regular} orbit equivalent to $G$ away from the single $G$-orbit~$E$.
\end{prop}
 \begin{figure}[htp]
 \captionsetup{width=0.96\linewidth}
    \centering
    \includegraphics[width=0.75\textwidth]{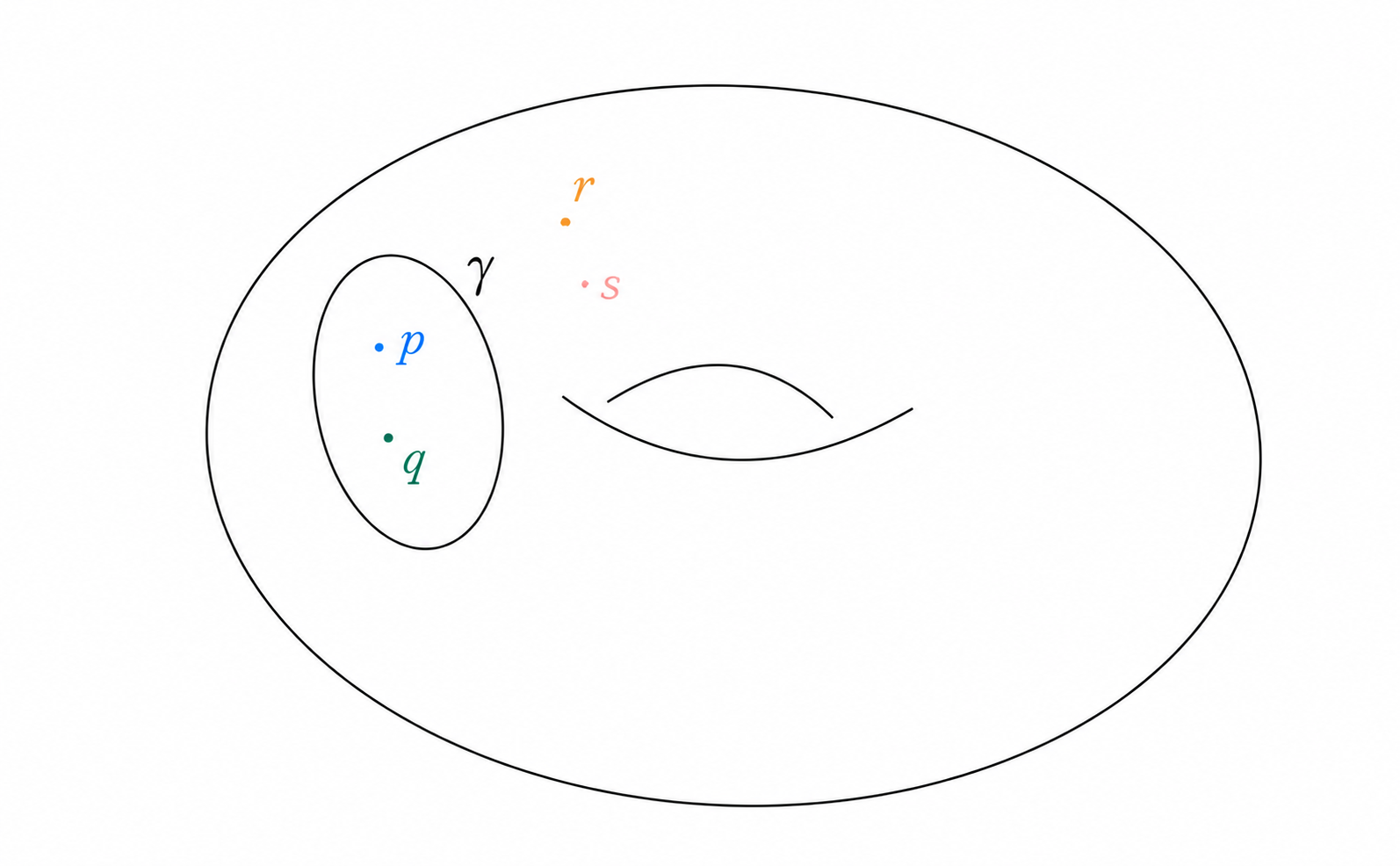}
    \includegraphics[width=1\textwidth]{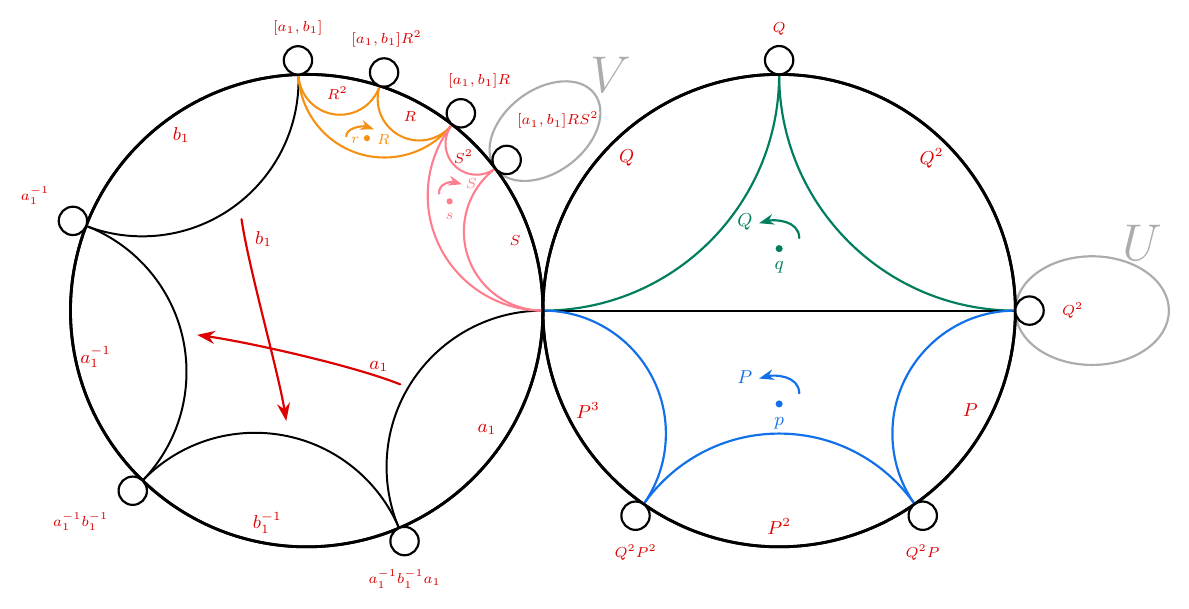}
    \caption{The Basilica Bowen--Series map for $\widetilde{G}$ on the boundary of the Bers slice of type $(1, 4; 3, 3, 3, 4)$.}
    \label{fig:(g,n>1)}
\end{figure}

\subsection{Case $n\geq 2$}
In this case, since $(g, n) \neq (0,3)$, there exists an essential simple closed curve $\gamma$ around two orbifold points, denoted by $p, q$. 
We obtain a nodal surface with two components after pinching $\gamma$. We have two sub-cases.
\begin{enumerate}
    \item\label{subcase:hyp} either both components are hyperbolic (see Figure~\ref{fig:(g,n>1)});
    \item\label{subcase:para} or exactly one component is parabolic (see Figure~\ref{fig:(1,2;2,2)}).
\end{enumerate}

We first consider the subcase~\eqref{subcase:hyp}.
Let $\widetilde{G}$ be the group on the boundary of the Bers slice obtained by pinching this curve $\gamma$. Note that the limit set $\Lambda(\widetilde{G})$ is homeomorphic to a Basilica. Let $\Delta_\infty\subset\Omega(\widetilde{G})$ be the totally invariant component. Then $\widetilde{G}|_{\Delta_\infty}$ is conformally conjugate to the Fuchsian group $G$.

The proof is similar to the closed case (Proposition~\ref{prop:closedsurfaceOE}), thus we only sketch the idea of the proof.

Let the cone points enclosed by \(\gamma\) be \(p,q\), and choose elliptic generators \(P,Q\) of orders \(\nu_p,\nu_q\), respectively, see Figure~\ref{fig:(g,n>1)}.
Orient the core polygon of this component so that the primitive element associated with the node is
\[
  W:=PQ.
\]
The same orientation is used for the node word on the complementary component.  Put \(m:=n-2\), the number of cone points on that component.
The hypothesis of subcase~\eqref{subcase:hyp} says precisely that
\[
  \frac{1}{\nu_p}+\frac{1}{\nu_q}<1
\]
and that the complementary one-cusped orbifold has negative orbifold
Euler characteristic.

There are two possibilities.  If \(m\geq1\), choose generators
\begin{equation}\label{eq:genenm>0}
  a_1,b_1,\ldots,a_g,b_g,\quad R_1,\ldots,R_{m-1},S
\end{equation}
for the complementary component, with \(S\) the terminal elliptic
generator adjacent to the node, and write
\begin{equation}\label{eq:orbifold-elliptic-node-word}
  A:=\prod_{i=1}^{g}[a_i,b_i]\,R_1\cdots R_{m-1} \text{ so that }   W=AS.
\end{equation}
Note that $A$ has word length
\[
  N_A:=4g+m-1.
\]
The hyperbolicity of the complementary component implies $N_A\geq1$.

If \(m=0\), then \(n=2\) and hyperbolicity of the complementary component
forces \(g\geq1\).  Put
\begin{equation}\label{eq:orbifold-handle-node-word}
  r:=a_g,\qquad s:=b_g,\qquad
  K:=\prod_{i=1}^{g-1}[a_i,b_i],\text{ so that }
  W=K[r,s],
\end{equation}
and define
\[
  \tau:=s^{-1}rs=rK^{-1}W,
  \qquad
  N_H:=4(g-1)+1=4g-3.
\]
Thus \(N_H\) is the word length of \(rK^{-1}\).

As in \S \ref{sss:closed-pinching-notation}, let \(I_L,I_R\) be the side intervals
associated with the two-cone and complementary components. 
Let $U, V$ be the limbs at ideal vertices of the $(\nu_p,\nu_q, \infty)$-orbifold polygon and the other core polygon which are mapped to the central node by $Q^{\nu_q-1}$ and $AS^2$ if $m\geq 1$ and $Ar^{-1}s^{-1}r$ if $m = 0$ respectively; see Figure~\ref{fig:g3}.  Denote the
corresponding intervals by $I_U$ and $I_V$. 
The element \(W\) fixes the common endpoints of \(I_L,I_R\), and hence preserves both side intervals.

As an analog of Lemma~\ref{lem:closed-most-generators-good}, we have
\begin{lem}\label{lem:goodgenN>1}
    Suppose $m \geq 1$. Then the generators
    \[ Q, a_i, b_i, R_j,i = 1,\cdots, g,  j = 1,\cdots, m-1\] are good. The generator $P$ is good away from $I_U \cup I_L$ and $S$ is good away from $I_V \cup I_R$.

    Suppose $m = 0$. Then the generators 
     \[ Q, a_i, b_i, r, i = 1,\cdots, g-1\] are good. The generator $P$ is good away from $I_U \cup I_L$ and $s$ is good away from $I_R$.    
\end{lem}

Recall that the stabilizer of the central node is generated by
\begin{align*}
    W = PQ = AS = \prod_{i=1}^{g}[a_i,b_i]\,R_1\cdots R_{m-1} &\text{ if } m \geq 1;\\
    W = PQ = K[r,s] = \prod_{i=1}^{g-1}[a_i,b_i][r, s] &\text{ if } m = 0.
\end{align*}

Similar to Lemma~\ref{lem:closed-node-branches},  we have the following by a straightforward computation.
\begin{lem}\label{lem:communonclosed}
If $x\in I_L\setminus\mathrm{bk}(\mathcal{F})$, then
\begin{equation}
  \mathcal{F}(Px)=Q^{\nu_q-1}P^{\nu_p-2}(Px)=Q^{-1}P^{-1}x = W^{-1}x\in I_L.
\end{equation}
If $x \in I_U\setminus\mathrm{bk}(\mathcal{F})$, then
\begin{equation}
  \mathcal{F}(x)=Q^{\nu_q-1}(x)=Q^{-1}P^{-1}(Px) = W^{-1}(Px)\in I_L.
\end{equation}

For $m \geq 1$, if $x\in I_R\setminus\mathrm{bk}(\mathcal{F})$, then
\begin{equation}\label{eqn:S}
  \mathcal{F}(Sx)=A(Sx)=Wx\in I_R.
\end{equation}

For $m \geq 0$, if $x\in I_R\setminus\mathrm{bk}(\mathcal{F})$, then
\begin{equation}\label{eqn:s}
  \mathcal{F}(sx)=Kr^{-1}s^{-1}r(sx)=Wx\in I_R.
\end{equation}
\end{lem}

As an analogue for Lemma~\ref{lem:closed-defect-transitions}, we have the following. The computation is similar to that in Lemma~\ref{lem:closed-defect-transitions}, and we omit it.
\begin{lem}\label{lem:iterateFnonclosed}
The following statements hold.
\begin{enumerate}
\item Suppose that $x\notin\mathrm{bk}(\mathcal{F})$ and $x,Wx\in I_R$.
Then either 
\begin{itemize}
    \item $Wx\sim_{\mathcal{F}} x$, or
    \item $Qx\in I_L$.
\end{itemize}
In the second case,
\[
  \mathcal{F}(x)=Q x,\qquad \mathcal{F}(Wx)=W^{-1}Q x,
\]
and both points belong to $I_L$.

\item Suppose that $m \geq 1$, $x\notin\mathrm{bk}(\mathcal{F})$ and $x,Wx\in I_L$.
Then either 
\begin{itemize}
    \item $Wx\sim_{\mathcal{F}} x$, or
    \item $x\in I_V$.
\end{itemize}
In the second case,
\[
  \mathcal{F}(x)=WS x,\qquad \mathcal{F}^{N_A}(Wx)=A^{-1}W x = Sx,
\]
and both points belong to $I_R$.

\item Suppose that $m = 0$, $x\notin\mathrm{bk}(\mathcal{F})$ and $x,Wx\in I_L$.
Then either 
\begin{itemize}
    \item $Wx\sim_{\mathcal{F}} x$, or
    \item $rsx = \cF^2(x)\in I_V$.
\end{itemize}
In the second case,
\[
  \mathcal{F}^3(x)=W\tau x,\qquad \mathcal{F}^{N_H}(Wx)= \tau x,
\]
and both points belong to $I_R$.
\end{enumerate}
\end{lem}

A similar argument as in Lemma~\ref{lem:closed-exceptional-points} gives the following.
\begin{lem}\label{lem:goodbcnonclosed}
The following statements hold.
\begin{enumerate}
    \item Suppose that $m \geq 1$. Let $H_R:=SW^{-1}Q$ and $H_L:=W^{-1}QS$, and let $x_R, x_L$ be the repelling fixed points of $H_R, H_S$ respectively.
Then for every regular $x$,
\[
  Wx\sim_{\mathcal{F}} x\quad\Longleftrightarrow\quad x\notin \{x_R, x_L\}.
\]
Consequently, 
\[
  Px\sim_{\mathcal{F}} x\Longleftrightarrow x\notin \{Qx_R, Qx_L\} \quad Sx\sim_{\mathcal{F}} x\Longleftrightarrow x\notin \{x_R, x_L\}.
\]
Their periods under $\mathcal{F}$ are
$\operatorname{per}_{\cF}(x_R)=N_A+1$ and $\operatorname{per}_{\cF}(x_L)=2$.
Their $\mathcal{F}$-grand orbits are distinct.  Finally,
\[
x_R=Sx_L,\qquad \cF^{N_A}(Wx_L)=x_R,\qquad \cF(Wx_R)=x_L.
\]
\item Suppose that $m = 0$. Let $H_R:=\tau W^{-1}Q$ and $H_L:=W^{-1}Q\tau$.
Let $x_R, x_L$ be the repelling fixed points of $H_R, H_S$ respectively.
Then for every regular $x$,
\[
  Wx\sim_{\mathcal{F}} x\quad\Longleftrightarrow\quad x\notin \{x_R, x_L\}.
\]
Consequently, 
\[
  Px\sim_{\mathcal{F}} x\Longleftrightarrow x\notin \{Qx_R, Qx_L\} \quad sx\sim_{\mathcal{F}} x\Longleftrightarrow x\neq x_R.
\]
Their periods under $\mathcal{F}$ are $\operatorname{per}_{\cF}(x_R)=N_H+1=4g-2$ and $\operatorname{per}_{\cF}(x_L)=4$. Their $\mathcal{F}$-grand orbits are distinct.  Finally,
\[
x_R=\tau x_L,\qquad \cF^{N_H}(Wx_L)=x_R,\qquad \cF(Wx_R)=x_L.
\]
\end{enumerate}
\end{lem}
    Let $x_L, x_R$ be the points as in Lemma~\ref{lem:goodbcnonclosed}. The same argument as in the proof of Proposition~\ref{prop:closedsurfaceOE} gives that if $x \notin \mathcal{O}_G(x_L)= \mathcal{O}_G(x_R)$, then $\mathcal{O}_G(x) = \mathcal{GO}_\mathcal{F}(x)$. Otherwise, $\mathcal{O}_G(x_L) = \mathcal{GO}_\mathcal{F}(x_L) \cup \mathcal{GO}_\mathcal{F}(x_R)$.
    This concludes the proof of Proposition~\ref{prop:nonclosedsurfaceOE} in the case $n \geq 2$ and subcase~\eqref{subcase:hyp}.

 \begin{figure}[htp]
\captionsetup{width=0.96\linewidth}
\centering
    \includegraphics[width=0.75\textwidth]{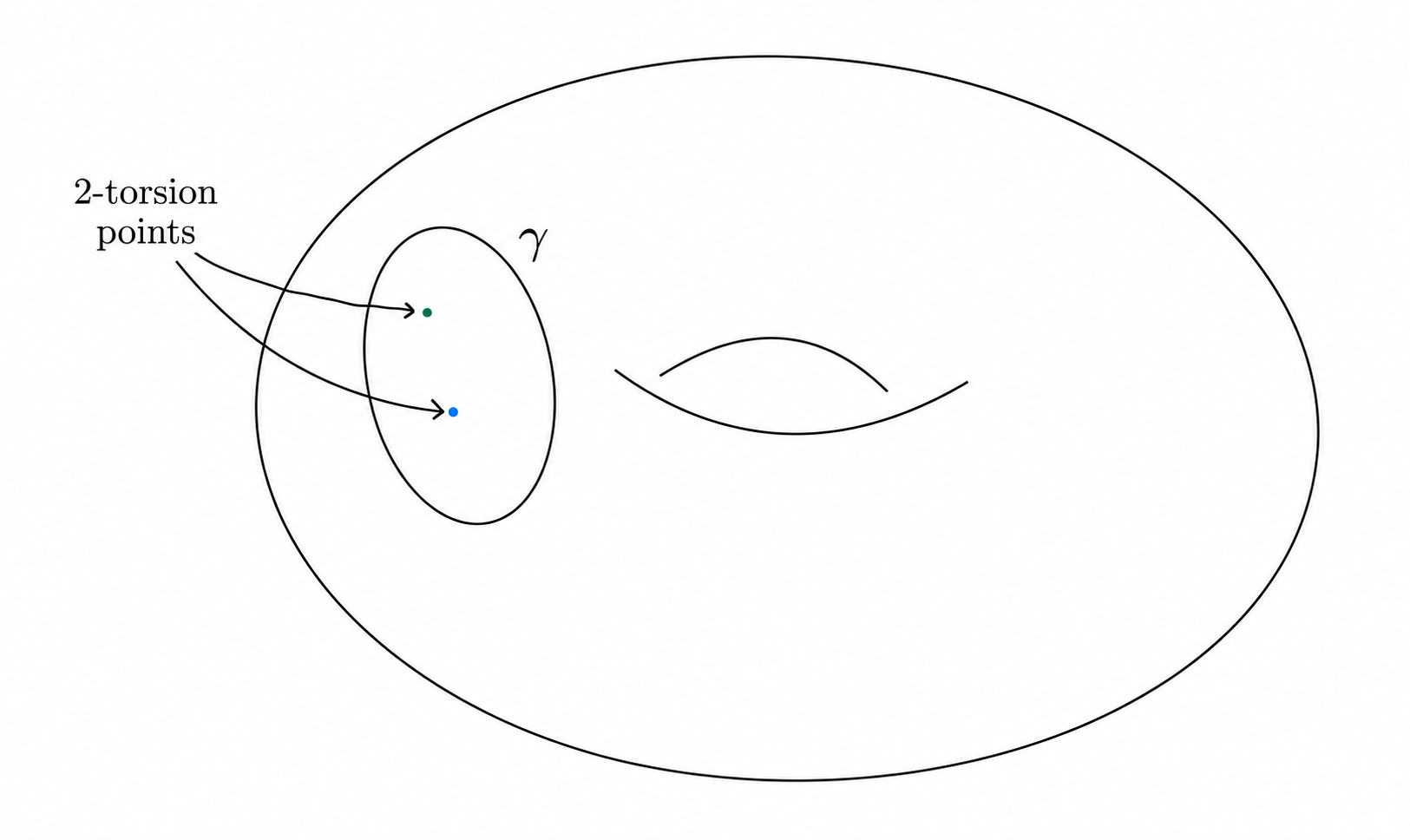}
    \includegraphics[width=1\textwidth]{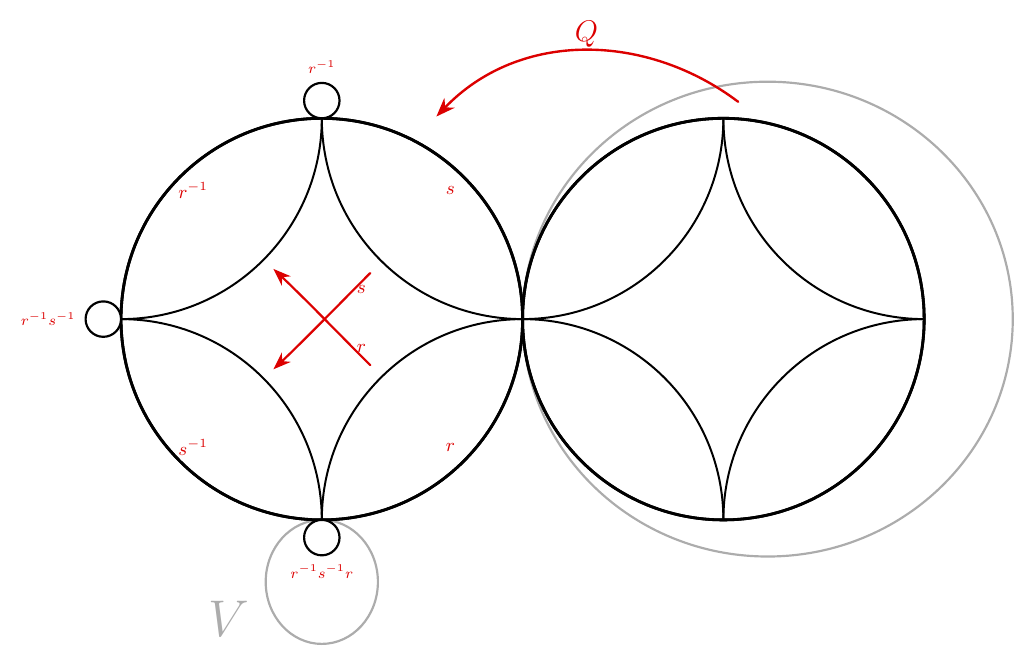}
    \caption{The Basilica Bowen--Series map for $\widetilde{G}$ on the boundary of the Bers slice of type $(1, 2; 2, 2)$.}
    \label{fig:(1,2;2,2)}
\end{figure}

\medskip

Now consider the subcase~\eqref{subcase:para}. 
A parabolic nodal component has signature $(0,3;2,2,\infty)$, so it is cut off by a curve enclosing two order-two cone points. Note the other nodal component is necessarily hyperbolic.

Let $P, Q$ be the two involutions associated to the two cone points, so that $P, Q$ have a common fixed point at the central node. We orient so that $W = PQ$ is a primitive generator of the parabolic subgroup at the central node.
As in \S \ref{sss:closed-pinching-notation}, let \(I_L,I_R\) be the side intervals of the central node. We define the Bowen--Series map so that $\cF = Q$ on $I_R$; see Figure~\ref{fig:(1,2;2,2)}.

Let $m = n-2$. Similar as in the subcase~\ref{subcase:hyp}, we choose generators in the complementary component as in \eqref{eq:genenm>0} and \eqref{eq:orbifold-elliptic-node-word} if $m \geq 1$, and as in \eqref{eq:orbifold-handle-node-word} if $m = 0$.

With the above notation and setup, the same computation gives that 
\begin{itemize}
    \item the generator $Q$ is good and the same conclusion as Lemma~\ref{lem:goodgenN>1} for the generators of the hyperbolic nodal component;
    \item the same conclusion of \eqref{eqn:S} and \eqref{eqn:s} in Lemma~\ref{lem:communonclosed};
    \item the same conclusion of Lemma~\ref{lem:iterateFnonclosed} and Lemma~\ref{lem:goodbcnonclosed}.
\end{itemize}
Therefore, the same argument gives Proposition~\ref{prop:nonclosedsurfaceOE} in the case $n \geq 2$ and subcase~\eqref{subcase:para}.

\subsection{Case $n=1$ and $g\geq2$}
 \begin{figure}[htp]
 \captionsetup{width=0.96\linewidth}
    \centering
    \includegraphics[width=0.75\textwidth]{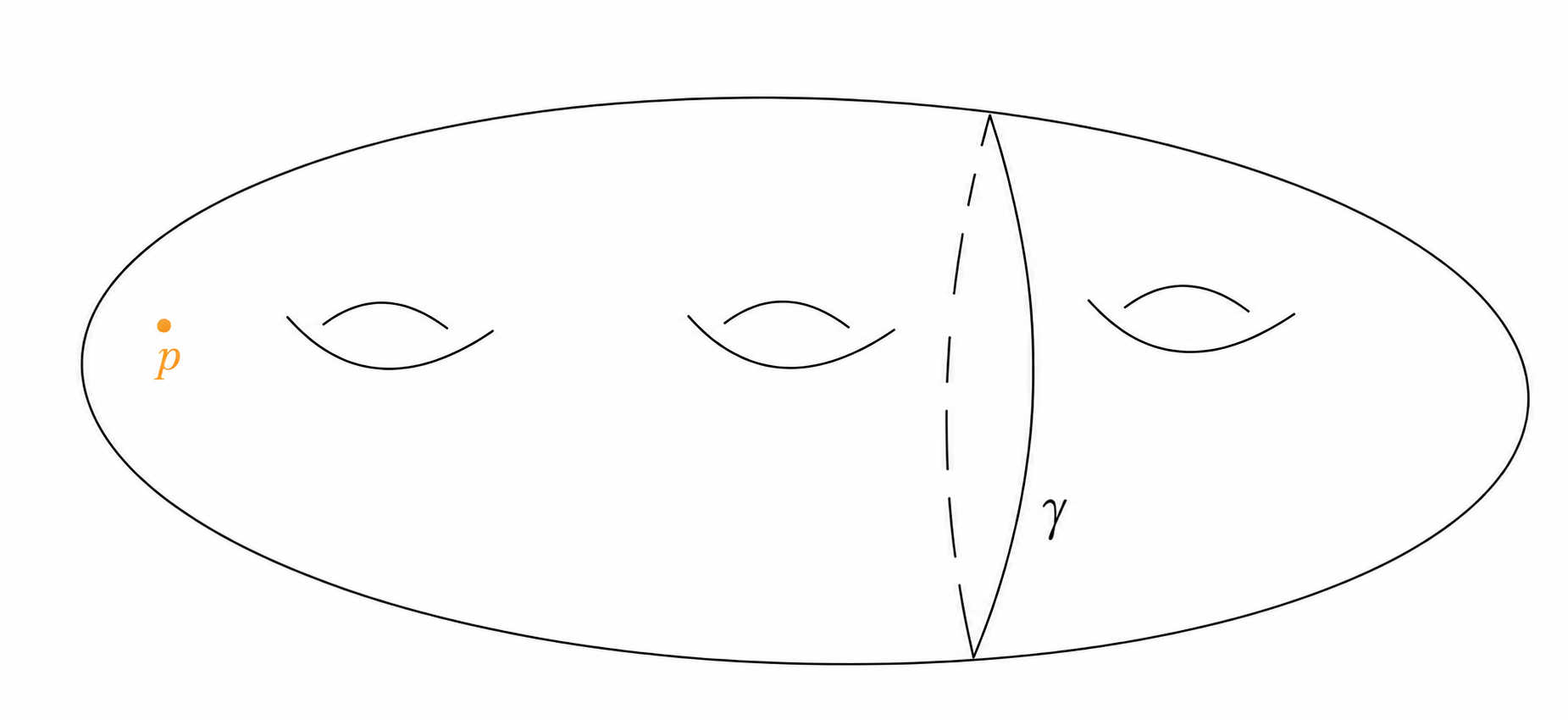}
    \includegraphics[width=1\textwidth]{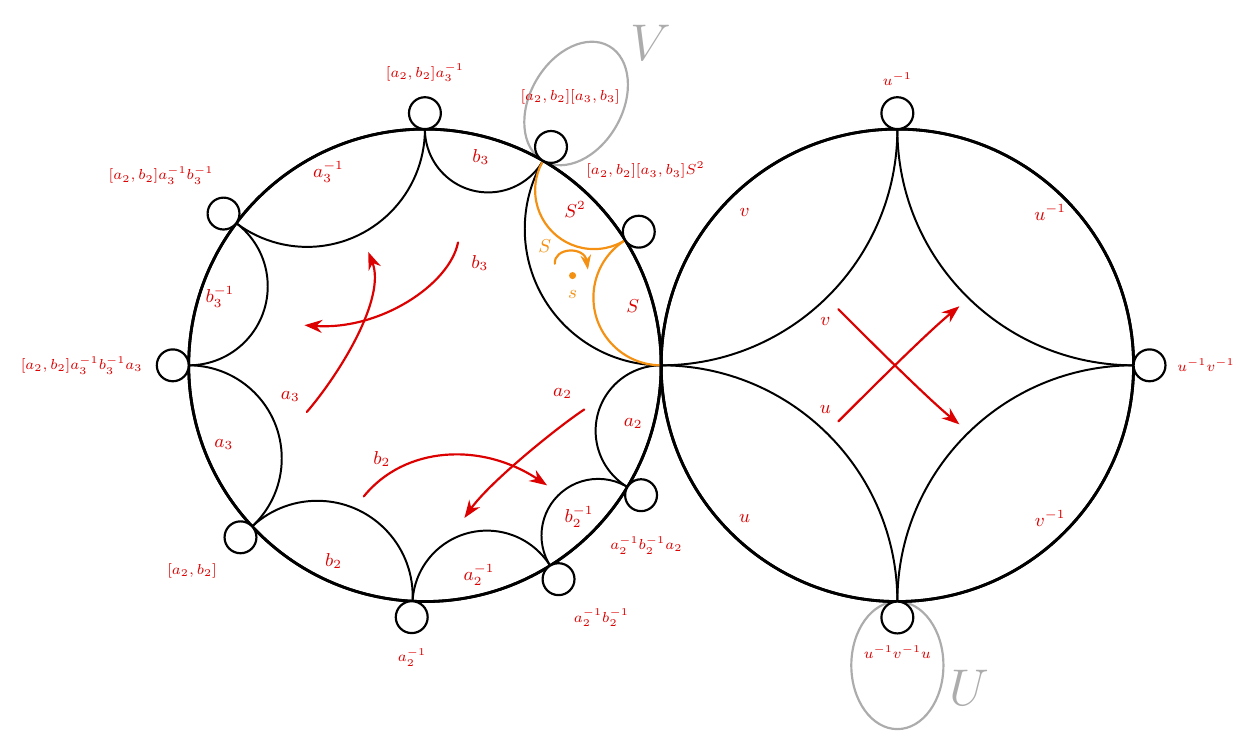}
    \caption{The Basilica Bowen--Series map for $\widetilde{G}$ on the boundary of the Bers slice of type $(3, 1; 3)$.}
    \label{fig:(g>1,1)}
\end{figure}

The case $n=0$ is Proposition~\ref{prop:closedsurfaceOE}, and the case $n\geq2$ was proved above.  Hyperbolicity rules out $(0,1)$, while $(1,1)$ is excluded, so this is the only remaining case. Let $s$ be the cone point, of order $\nu$, and pinch a separating curve disjoint from $s$ which cuts off a one-holed torus. We put the resulting once-punctured torus on the right and the component of signature $(g-1,2;\infty,\nu)$ on the left (see Figure~\ref{fig:(g>1,1)}).

Choose handle generators $u,v$ on the right, handle generators $a_j,b_j$ $(2\leq j\leq g)$ on the
left, and the elliptic generator $S$ of order $\nu$.  Orient them so that
\begin{equation}\label{eq:orbifold-one-cone-words}
\begin{aligned}
  A&:=\prod_{j=2}^{g}[a_j,b_j],&
  W&:=[u,v]=AS, &
  \sigma&:=v^{-1}uv.
\end{aligned}
\end{equation}
Let $N_A = 4(g-1)$ be the length of $A$.
Let $I_R,I_L$ be the two side intervals and let $I_U\subset I_R$, $I_V\subset I_L$ be the corresponding intervals of the limbs where the Basilica Bowen--Series map is $\sigma, AS^2$ respectively as shown in Figure~\ref{fig:(g>1,1)}. 

The same computations for $u, v$ in Lemma~\ref{lem:closed-most-generators-good} and for $a_j, b_j, S$ in the case of $m \geq 1$ of Lemma~\ref{lem:goodgenN>1} gives the following.
\begin{lem}\label{lem:g>=2-most-generators-good}
The generators
\[
  u\quad \text{and}\quad a_j,b_j,\ j\in\{1,\cdots, g-1\},
\]
are good.  The generator $v$ is good away from $I_L$, and $S$ is good away from $I_V \cup I_R$.  
\end{lem}

The same computation as the first case of Lemma~\ref{lem:closed-defect-transitions} and the second case of Lemma~\ref{lem:iterateFnonclosed} gives the following.
\begin{lem}\label{lem:g>=2-defect-transitions}
The following statements hold.
\begin{enumerate}
\item Suppose that $x\notin\mathrm{bk}(\mathcal{F})$ and $x,Wx\in I_R$.
Then either 
\begin{itemize}
    \item $Wx\sim_{\mathcal{F}} x$, or
    \item $uvx=\mathcal{F}^2(x)\in I_U$.
\end{itemize}
In the second case,
\[
  \mathcal{F}^3(x)=W\sigma x,\qquad \mathcal{F}(Wx)=\sigma x,
\]
and both points belong to $I_L$.

\item Suppose that $m \geq 1$, $x\notin\mathrm{bk}(\mathcal{F})$ and $x,Wx\in I_L$.
Then either 
\begin{itemize}
    \item $Wx\sim_{\mathcal{F}} x$, or
    \item $x\in I_V$.
\end{itemize}
In the second case,
\[
  \mathcal{F}(x)=WS x,\qquad \mathcal{F}^{N_A}(Wx)=A^{-1}W x = Sx,
\]
and both points belong to $I_R$.
\end{enumerate}
\end{lem}

A similar argument as in Lemma~\ref{lem:closed-exceptional-points} gives the following.
\begin{lem}\label{lem:g>=2exceptional-points}
Let $H_R:=S\sigma$ and $H_L:=\sigma S$, and let $x_R, x_L$ be the repelling fixed points of $H_R, H_L$ respectively.
Then for every regular $x$,
\[
  Wx\sim_{\mathcal{F}} x\quad\Longleftrightarrow\quad x\notin\{x_R, x_L\}.
\]
Consequently,
\[
vx\sim_{\mathcal{F}} x \Longleftrightarrow x \neq x_L \quad Sx\sim_{\mathcal{F}} x \Longleftrightarrow x \notin \{x_R, x_L\}
\]
Their periods under $\mathcal{F}$ are
\[
  \operatorname{per}_{\mathcal{F}}(x_R)=N_A+3,
  \qquad
  \operatorname{per}_{\mathcal{F}}(x_L)=2.
\]
Their $\mathcal{F}$-grand orbits are distinct.  Finally,
\[
  x_R = S x_L \qquad
  \mathcal{F}^{N_A}(Wx_L)=x_R,
  \qquad
  \mathcal{F}(Wx_R)=x_L.
\]
\end{lem}
Let $x_L, x_R$ be the points as in Lemma~\ref{lem:g>=2-defect-transitions}. The same argument as in the proof of Proposition~\ref{prop:closedsurfaceOE} gives that if $x \notin \mathcal{O}_G(x_L)= \mathcal{O}_G(x_R)$, then $\mathcal{O}_G(x) = \mathcal{GO}_\mathcal{F}(x)$. Otherwise, $\mathcal{O}_G(x_L) = \mathcal{GO}_\mathcal{F}(x_L) \cup \mathcal{GO}_\mathcal{F}(x_R)$.
This concludes the proof of Proposition~\ref{prop:nonclosedsurfaceOE} in the case $n=1, g \geq 2$, and Proposition~\ref{prop:nonclosedsurfaceOE} follows.

\section{Acceleration of Bowen–Series maps for non-uniform lattices with torsion}\label{sec:torsionproof}
In this section, we provide the details of the proof of Theorem~\ref{thm:accnon-uniformlattice} in the case where $G$ has torsion points. While the idea of the proof remains the same as in the torsion-free case, finding suitable break-points for carrying out the acceleration needs more book-keeping in the current setting.

\begin{proof}[Proof of Theorem~\ref{thm:accnon-uniformlattice} when $G$ has torsion points]
        Suppose that $G$ has signature $(g, n; \nu_1,\ldots, \nu_n)$ with $\nu_1 = \nu_2 = \ldots = \nu_k = \infty$, and $\nu_{k+1}, \ldots, \nu_{n} < \infty$.
    Let $m = n - k$ be the number of orbifold points.
    Note that in our case, $k \geq 1, m \geq 1$.
    Let $\mathcal{P}$ be a core polygon with corresponding Bowen--Series map $\mathcal{F}_G$.

    Note that the vertices $a_i$ form $k$ orbits under $G$. More precisely, they are 
    \begin{align*}
    &\{a_0,\ldots, a_{4g}, a_{4g+2k-2},a_{4g+2k-1}, \cdots., a_{N-1} = a_{4g+2(k-1)+(\sum_{s=k+1}^n(\nu_s-1)) - 1}\}, \\
    &\{a_{4g+1}, a_{4g+2k-3}\},\ldots, \{a_{4g+k-i}, a_{4g+k+i-2}\}, \ldots, \{a_{4g+k-2}, a_{4g+k}\}, \{a_{4g+k-1}\}.
    \end{align*}
    We pick one representative in each orbit: 
    $$
    S = \begin{cases}
        \{a_0, a_{4g+1}, a_{4g+2},\ldots, a_{4g+k-1}\} & \text{if } g = 0\\
        \{a_2, a_{4g+1}, a_{4g+2},\ldots, a_{4g+k-1}\} & \text{if } g > 0
    \end{cases}.
    $$
    Let $\Phi: \{a_0,\ldots, a_{N-1}\} \longrightarrow S$ be the map that sends $a_i$ to the unique element in $S$ that is in the same $G$ orbit as $a_i$.
    Let 
    $$
    S_1 := \Big\{a_i: i\in\big\{0,\cdots, a_{4g+2(k-1)}\big\}\Big\}\ \mathrm{and}\ S_2:= \Big\{a_{i,0}: i \in \big\{0,\cdots, n-k-1\big\}\Big\},
    $$ 
    where $a_{i,j}$ is defined in Formula~\eqref{a_i_j_formula}.
    Note that by construction 
    $$
    \cF_G(a_{i,j}^+) = a_{i, \nu_{k+1+i}-1} = a_{i+1, 0}\quad \mathrm{and}\quad \cF_G(a_{i,j}^-) = a_{i, 0}.
    $$
    Thus $\cF_G(a_i^\pm) \in S_1 \cup S_2$ for all $i \in \{0,\cdots, N-1\}$.
    
    We consider three subcases.
\smallskip

\noindent\textbf{Case~(a):} $g \geq 1$, or $g = 0, k \geq 2, m \geq 2$, or $g = 0, k = 1, m \geq 3$.

    Let $r = \begin{cases}
        \max\{2, m\} &\text{if } g+k-1 > 0,\\
        \max\{2, m-1\} &\text{if } g+k-1 = 0.
    \end{cases}$
    
    Note that for sufficiently small $\varepsilon$, $\cF_G^r|_{[a_{i,0}, a_{i,0}+\varepsilon]}$ is given by a loxodromic element in $G$ for $i\in \{0,\cdots, n-k-1\}$, and this element has a unique repelling fixed point $\alpha_j \in [a_{i,0}, a_{i,1}] = [a_j, a_{j+1}]$ for some $j$. Similarly, $\cF_G^r|_{[a_{i+1,0}-\varepsilon, a_{i+1,0}]}$ is loxodromic for $i\in \{0,\cdots, n-k-2\}$, and has a unique repelling fixed point $\alpha_j \in [a_{i,\nu_{k+1+i}-2}, a_{i+1,0}] = [a_j, a_{j+1}]$ for some $j$.

    The construction is now similar to the torsion-free case. We construct $\mathcal{H}_G$ as follows.
    Let $I^0:=[a_i, a_{i+1}]$. We define $\mathcal{H}_{G,0}|_{I^0} = \mathcal{F}_{G}|_{I^0}$.

    Suppose that $\mathcal{H}_{G,s-1}|_{I^{0}}$ is defined.
Let $I^s = \mathcal{H}_{G,s-1}|_{I^{0}}(I^0) = [x^{s,-}, x^{s,+}]$. Abusing notations slightly, we view $[x^{s,-}, x^{s,+}]$ as an interval in the universal cover of the circle, so that it may project to the circle, wrapping around multiple times. Suppose that $x^{s,\pm} \in S$. Then we set $\mathcal{H}_G|_{I^0} = \mathcal{H}_{G,s-1}|_{I^{0}}$.
Otherwise, suppose that $x^{s,-} \not\in S$. Let $[x^{s,-}, \widetilde{x}^{s,+}] \subseteq I^s$ so that $[x^{s,-}, \widetilde{x}^{s,+}]$ is one of the level $0$ intervals, i.e., $[x^{s,-}, \widetilde{x}^{s,+}] = [a_j, a_{j+1}]$ for some $j = j(s, -)$.
Note that $a_j \in S_1 \cup S_2$.
Thus, as constructed earlier, there exists some fixed point $\alpha_j \in [a_j, a_{j+1}]$ of $\cF_G$ or $\cF_G^2$ depending on whether $j < 4g+2(k-1)$ or $j \geq 4g+2(k-1)$.
Let $[a_i, \widetilde{\alpha}^{s,+}] = (\mathcal{H}_{G,s-1}|_{I^{0}})^{-1}([x^{s,-}, \alpha_j])$. Similarly, we can define $[\widetilde{\alpha}^{s,-}, a_{i+1}]$ with the understanding that $\widetilde{\alpha}^{s,-} = a_{i+1}$ if $x^{s,+} \in S$.
Let
$$
\delta(j):= \begin{cases}
    1 & \text{ if }0 \leq j < 4g+2(k-1),\\
    r & \text{ if }j \geq 4g+2(k-1).
\end{cases}
$$
We define
$$
\mathcal{H}_{G,s}|_{I^{0}}(x) = \begin{cases}
    \mathcal{F}_G^{\delta(j(s,-))} \circ \mathcal{H}_{G,s-1}|_{I^{0}}(x)& \text{ if } x\in [a_i, \widetilde{\alpha}^{s,+}],\\
    \mathcal{F}_G^{\delta(j(s,+))} \circ \mathcal{H}_{G,s-1}|_{I^{0}}(x)& \text{ if } x\in [\widetilde{\alpha}^{s,-}, a_{i+1}],\\
    \mathcal{H}_{G,s-1}|_{I^{0}}(x)& \text{ if } x\in [\widetilde{\alpha}^{s,+}, \widetilde{\alpha}^{s,-}].
\end{cases}
$$
One can verify that $\cH_{G, s}$ stabilizes after finitely many iterations as in the torsion-free case. Let $\mathcal{H}_G = \mathcal{H}_{G,s}$ for $s$ sufficiently large. One can also verify that the additional break-points introduced are symmetric, and $\cH_G$ satisfies the conditions required by the theorem.
\smallskip

\noindent\textbf{Case~(b):} $g = 0, k \geq 2, m = 1$. 

Let 
    $$
    B = \bigcup_{j=0}^{\nu_{k+1}-2} [a_{0,j}, a_{0,j+1}] = [a_{0,0}, a_{0, \nu_{k+1}-2}] = [a_{2(k-1)}, a_{0}]
    $$ 
    be the union of all intervals associated to the unique orbifold point. 
    
    Recall that $\alpha_i \in [a_i, a_{i+1}]$ is the fixed point of $\mathcal{F}_G$ for $i=0,..., 2(k-1)$.
    Note that for $i = 1,..., k-2$, the vertices $a_i, a_{2(k-1)-i}$ are exchanged by $\mathcal{F}_G$. To make $a_i$ fixed, define
    $$
    d_i^- = (\mathcal{F}_G|_{[a_{i-1}, a_i]})^{-1}(\alpha_{C-i-1}), \qquad d_i^+ = (\mathcal{F}_G|_{[a_{i}, a_{i+1}]})^{-1}(\alpha_{C-i}).
    $$
    Then we have
    $$
    a_i \in [d_i^-, d_i^+], \;\; \mathcal{F}^2|_{G}(d_i^\pm) = \mathcal{F}|_{G}(d_i^\pm) \text{ and } \mathcal{F}^2|_{G}(a_i) = a_i.
    $$
    Let $J =[a_{0,j}, a_{0,j+1}] \subseteq B$ be an interval. We define $\alpha_J\coloneqq (\mathcal{F}_G|_J)^{-1}(\alpha_{2(k-1)-1})$, and
    $$
    h_J(x) = \begin{cases}
        \mathcal{F}_G^2(x) & \text{if } x\in [a_{0,j}, \alpha_J],\\
        \mathcal{F}_G(x) & \text{if } x\in [\alpha_J, a_{0, j+1}].
    \end{cases}
    $$
    Note that $h_J(\partial J) = a_0$ and is continuous by construction.
    Let $h_B: B \rightarrow \bS^1$ be the union of all such maps. Then $h_B$ is continuous and sends $a_{0,j}$ to~$a_0$.

    Let $\beta_0  = (\mathcal{F}_G|_{[a_0,a_1]})^{-1}(a_0)= P_0(a_0)$. Then $\mathcal{F}_G([a_0, \beta_0]) = B$.
    We define
    $$
    \mathcal{H}_G(x) = \begin{cases}
        \mathcal{F}_G^2(x) & \text{if } x\in \bigcup_{i=1}^{k-2}[d_i^-, d_i^+],\\
        h_B \circ \mathcal{F}_G(x) & \text{if } x\in [a_0, \beta_0],\\
        h_B(x) & \text{if } x\in B,\\
        \mathcal{F}_G(x) & \text{otherwise.}
    \end{cases}
    $$
    One can verify that $\cH_G$ satisfies the conditions required by the theorem.
    \smallskip
    
\noindent\textbf{Case~(c):} $g = 0, k = 1, m = 2$.

    Let $A = \bigcup_{i=0}^{\nu_2-2} I_i$ and $B = \bigcup_{i=0}^{\nu_3-2} J_i$ be the intervals associated with the two orbifold points.
    Note that $\bS^1 = A \cup B$ and at least one of the orbifold points has signature $\geq 3$. Without loss of generality, we assume that $\nu_2 \geq 3$.

    Recall $I_{i} = [a_{0,i}, a_{0,i+1}]$, with $a_0 = a_{0,0}$.
    Let $\alpha$ be the fixed point for $\mathcal{F}^4_G|_{[a_{i,1}, a_{i,1}+\varepsilon]}$. Note that $\alpha$ has periodic itinerary $(I_1, J_0, I_0, J_0)$ of period $4$ under $\mathcal{F}_G$, i.e., 
    $$
    \alpha \in I_1, \;\; \mathcal{F}_G(\alpha) \in J_0, \;\; \mathcal{F}^2_G(\alpha) \in I_0, \;\; \mathcal{F}^3_G(\alpha) \in J_0 \;\; ...
    $$
    Let $\beta\coloneqq I_0\cap I_1 = \{a_{0,1}\}$.
    By construction, we have $\cF_G^4([\beta, \alpha]) = [a_0, \alpha]$.
    
    Let $J_i \subseteq B$ be an interval in $B$. We define 
    $$
    \alpha_i\coloneqq (\cF_G|_{J_i})^{-1}(\alpha), \text{ and } \beta_i\coloneqq (\cF_G|_{J_i})^{-1}(\beta),
    $$
    and
    $$
    h_{J_i}(x) = \begin{cases}
        \cF_G(x) & \text{if } x \in [a_{1,i}, \beta_i],\\
        \cF_G^5(x) & \text{if } x \in [\beta_i, \alpha_i],\\
        \cF_G^2(x) & \text{if } x \in [\alpha_i, _{1,i+1}].
    \end{cases}
    $$
    Note that $h_{J_i}(\partial J_i) = a_0$ and is continuous by construction. Let $h_B: B \rightarrow \bS^1$ be the union of all such maps. Then $h_B$ is continuous and sends $a_{1,j}$ to $a_0$. We define
    $$
    \cK_G(x) = \begin{cases}
        h_B \circ \cF_G(x) & \text{if } x\in A,\\
        h_B(x) & \text{if } x\in B.
    \end{cases}
    $$
    Note that $\cK_G$ is continuous and every periodic break-point is symmetric except $\beta$.
    Note that $\beta$ has period $2$ and there exists $W\in G$ such that $\cK_G|_{[\beta, \beta+\varepsilon)} = W$ and $\cK_G|_{(\beta-\varepsilon, \beta]}$.
    
    Choose regular periodic points $y_\pm$ near $\cK_G(\beta)$ of distinct odd prime periods $r_\pm$. Let $x_\pm$ be the inverse images so that $\cK_G(x_\pm) = y_\pm$ and $\beta \in [x_-, x_+]$. We also require the cycles of $y_\pm$ to avoid the closed interval $[x_-, x_+]$ and this interval to avoid $\cK_G(\beta)$ and any break-points.
    By the Chinese remainder theorem, we can choose $L$ so  that $r_- \equiv 1 \mod 2L$ and $r_+ \equiv 1 \mod 4L$.
    Now define
    $$
    \cH_G(x) = \begin{cases}
        \cK_G^{2L}(x) & \text{if } x\in [x^-, \beta],\\
        \cK_G^{4L}(x) & \text{if } x\in [\beta, x^+],\\
        \cK_G(x) & \text{otherwise.}
    \end{cases}
    $$
    One can verify that $\cH_G$ satisfies the conditions required by the theorem.
\end{proof}
 
\begin{remark}
For a different construction of conformal Markov (but not piecewise M{\"o}bius) coverings of the circle which are \emph{virtually orbit equivalent} to the genus-zero orbifold groups treated in Cases~(b) and~(c) of the above proof, see \cite{MM2,LLM24,LMM26c}.
\end{remark}

\subsection*{AI Disclosure}
AI tools were used for text editing and figure generation from hand-drawn figures.
AI tools were not used for any of the mathematical content, including proofs or proof ideas.


\begin{thebibliography}{LMMN25}

\bibitem[AC96]{AC96}
J.~W.~Anderson and R.~D.~Canary.
\newblock Cores of hyperbolic 3-manifolds and limits of Kleinian groups.
\newblock {\em Amer. J. Math.}, 118:745--779, 1996.

\bibitem[AIM09]{AIM09}
	K.~Astala, T.~Iwaniec, and G.~Martin.
	\newblock {\em Elliptic partial differential equations and quasiconformal
		mappings in the plane}, volume 148 of {\em Princeton Mathematical Series}.
	\newblock Princeton Univ. Press, Princeton, NJ, 2009.

\bibitem[BL20]{BL1}
S.~Bullett and L.~Lomonaco.
\newblock Mating quadratic maps with the modular group {II}.
\newblock {\em Invent. Math.}, 220:185--210, 2020.

\bibitem[BL22]{BL2}
S.~Bullett and L.~Lomonaco.
\newblock Dynamics of Modular Matings.
\newblock {\em Adv. Math.}, 410, Paper No. 108758, 43 pp, 2022.

\bibitem[BL24]{BL3}
S.~Bullett and L.~Lomonaco.
\newblock Mating quadratic maps with the modular group {III}: the modular Mandelbrot set.
\newblock {\em Adv. Math.}, 458, Paper No. 109956, 48 pp., 2024.

\bibitem[BLL26]{BLL26}
S.~Bullett, L.~Lomonaco, and M.~R.~Laude.
\newblock Matings between compositions of rational maps and free products of finite cyclic groups.
\newblock \url{https://arxiv.org/abs/2607.14364}, 2026.
	
\bibitem[BP94]{BP94}
S.~Bullett and C.~Penrose.
\newblock Mating quadratic maps with the modular group.
\newblock {\em Invent. Math.}, 115:483--511, 1994.

\bibitem[BF14]{BF14} 
B.~Branner and N.~Fagella, 
\newblock \emph{Quasiconformal Surgery in Holomorphic Dynamics}. 
\newblock Cambridge Studies in Advanced Mathematics, vol.\ 141, Cambridge University Press, Cambridge, 2014. 

\bibitem[BCM12]{BCM12}
J.F. Brock, R. Canary and Y. Minsky.
\newblock  The classification of Kleinian surface groups, II: The ending lamination conjecture.  
\newblock {\em Ann. of Math. (2)}, 176, no. 1,: 1--149, 2012.

\bibitem[BLLM26]{BLLM}
S.~Bullett, L.~Lomonaco, M.~Lyubich, and S.~Mukherjee.
\newblock Mating parabolic rational maps with Hecke groups.
\newblock {\em Proc. London Math. Soc. (3)}, vol. 132, no. 3, Paper no. e70132, 2026.

\bibitem[BS79]{BS79}
R.~Bowen and C.~Series.
\newblock Markov maps associated with {F}uchsian groups.
\newblock {\em Inst. Hautes \'{E}tudes Sci. Publ. Math.}, 50:153--170, 1979.

\bibitem[DH85]{DH85}
A.~Douady and J.~H. Hubbard.
\newblock On the dynamics of polynomial-like mappings.
\newblock {\em Ann. Sci. {\'E}c. Norm. Sup.}, 18:287--343, 1985.

\bibitem[Fat29]{Fatou29}
P.~Fatou.
\newblock Notice sur les travaux scientifiques de {M}. {P}. {F}atou.
\newblock Astronome titulaire de l'observatoire de Paris, 5--29, 1929,
\url{https://www.math.purdue.edu/~eremenko/dvi/fatou-b.pdf}.

\bibitem[LMM22]{LMM22}
K.~Lazebnik, N.~G. Makarov, and S.~Mukherjee.
\newblock Bers slices in families of univalent maps.
\newblock {\em Math. Z.}, 300:2771--2808, 2022.

\bibitem[LLMM21]{LLMM21}
S.-Y. Lee, M.~Lyubich, N.~G. Makarov, and S.~Mukherjee.
\newblock Schwarz reflections and anti-holomorphic correspondences.
\newblock {\em Adv. Math.}, 385:Paper No. 107766, 88, 2021.

\bibitem[LLMM23]{LLMM23}
S.-Y. Lee, M.~Lyubich, N.~G. Makarov, and S.~Mukherjee.
\newblock Dynamics of Schwarz reflections: the mating phenomena.
\newblock {\em Ann. Sci. {\'E}cole Norm. Sup. (4)}, 56:1825--1881, 2023. 

\bibitem[LLM24]{LLM24}
Y.~Luo, M.~Lyubich, and S.~Mukherjee.
\newblock A general dynamical theory of Schwarz reflections, B-involutions, and algebraic correspondences.
\newblock \url{https://arxiv.org/abs/2408.00204}, 2024.

\bibitem[LM26b]{LM25b}
L.~Lomonaco and S.~Mukherjee.
\newblock Algebraic correspondences and Schwarz reflections: Where rational dynamics meets Kleinian groups.
\newblock {\em Proceedings of the ICM 2026}, vol. 5,
pages 168–190, 2026. 

\bibitem[LMM26a]{LMM26a}
Y.~Luo, M.~Mj, and S.~Mukherjee.
\newblock Teichm{\"u}ller spaces, polynomial loci, and degeneration in spaces of algebraic correspondences.
\newblock {\em J. Reine Angew. Math.}, 839:151--223, 2026.

\bibitem[LMM26b]{LMM26b}
Y.~Luo, M.~Mj, and S.~Mukherjee.
\newblock Universality of the Basilica.
\newblock \url{https://arxiv.org/abs/2601.13553}, 2026.

\bibitem[LMM26c]{LMM26c}
Y.~Luo, M.~Mj, and S.~Mukherjee.
\newblock Combining cusped triangle groups with Blaschke products: commensurable matings.
\newblock \url{https://arxiv.org/abs/2607.07886}, 2026

\bibitem[LMM24]{LMM24}
M.~Lyubich, J.~Mazor, and S.~Mukherjee.
\newblock Antiholomorphic correspondences and mating I: realization theorems.
\newblock {\em Comm. Amer. Math. Soc.}, 4:495--547, 2024.

\bibitem[LM26a]{LM25a}
M.~Lyubich and S.~Mukherjee.
\newblock Mirrors of conformal dynamics: Interplay between anti-rational maps, reflection groups, Schwarz reflections, and correspondences.
\newblock In {\em Algebraic, Complex, and Arithmetic
Dynamics}, pp 341–476, edited by Laura DeMarco and Mattias Jonsson, Simons Symposia, Springer, 2026.



\bibitem[LMMN25]{LMMN}
M.~Lyubich, S.~Merenkov, S.~Mukherjee, and D.~Ntalampekos.
\newblock David extension of circle homeomorphisms, welding, mating, and removability.
\newblock {\em Mem. Amer. Math. Soc.}, vol. 313, no. 1588, pp. v+110, 2025.

\bibitem[LN24]{LN24}
Y.~Luo and D.~Ntalampekos.
\newblock Uniformization of gasket Julia sets.
\newblock To appear in {\em J. Eur. Math. Soc. (JEMS)}, \url{https://arxiv.org/abs/2411.17227}, 2024.

\bibitem[LN26]{LN26}
Y.~Luo and D.~Ntalampekos.
\newblock Piecewise quasiconformal dynamical systems of the unit circle.
\newblock {\em Ergodic Theory Dynam. Systems}, 46(12):3082-3122, 2026.

\bibitem[LP97]{LP97}
G.~Levin and F.~Przytycki.
\newblock When do two rational functions have the same {J}ulia set?
\newblock {\em Proc. Amer. Math. Soc.}, 125:2179–2190, 1997.

\bibitem[McM94]{McM94}
C. McMullen.
\newblock {\em Complex Dynamics and Renormalization (AM-135)}.
\newblock Princeton University Press, 1994.

\bibitem[McM96]{McM96}
C.~T. McMullen.
\newblock {\em Renormalization and $3$-manifolds which fiber over the circle}.
\newblock Ann. of Math. Stud., 142.
Princeton University Press, Princeton, NJ, 1996. x+253 pp.

\bibitem[MR26]{manning-ruffoni}
J.F.~Manning and L.~Ruffoni.
\newblock {\em Incubulable hyperbolic 3-pseudomanifold groups}.
\newblock preprint, arxiv:2601.02655 2026.




\bibitem[MS98]{MS98}
C.~T.~McMullen and D.~P.~Sullivan.
\newblock Quasiconformal homeomorphisms and dynamics. {III}. The {T}eichm{\"u}ller space of a holomorphic dynamical system.
\newblock {\em Adv. Math.}, 135:351--395, 1998.

\bibitem[Mil06]{Mil06}
J.~Milnor.
\newblock {\em Dynamics in one complex variable}, volume 160 of {\em Annals of
  Mathematics Studies}.
\newblock Princeton University Press, Princeton, NJ, third edition, 2006.

\bibitem[Mj14a]{Mj14}
M.~Mj.
\newblock {Cannon-Thurston maps for surface groups}.
\newblock {\em Ann. of Math. (2)}, 179:1--80, 2014.

\bibitem[Mj14b]{mahan-elct}
M.~Mj.
\newblock { Ending laminations and Cannon-Thurston maps,} With an appendix by Shubhabrata Das and Mj. \newblock {\em Geom. Funct. Anal. 24, no. 1}, 297--321, 2014. 

\bibitem[Mj17]{mahan-kl}
M.~Mj.
\newblock { Cannon-Thurston maps for Kleinian groups.}
\newblock {\em Forum Math. Pi 5, e1, 49 pp., 2017}

\bibitem[Mj18]{mahan-icm}
M.~Mj.
\newblock {  Cannon-Thurston maps.}
\newblock {\em Proceedings of the International Congress of Mathematicians—Rio de Janeiro 2018. Vol. II. Invited lectures, 885--917, World Sci. Publ., Hackensack, NJ, 2018.}

\bibitem[Mj22]{mahan-survey}
M.~Mj.
\newblock { Cannon-Thurston maps in Kleinian groups and geometric group theory.}
\newblock {\em  Surveys in differential geometry} 2020. Surveys in 3-manifold topology and geometry, 281--318, Surv. Differ. Geom., 25, Int. Press, Boston, MA, 2022. 

\bibitem[Mj25]{Mj-cmotions}
M.~Mj.
\newblock{Motions of Limit Sets: A Survey}
\newblock {\em Teichm\"uller Theory and Grothendieck-Teichm\"uller Theory}, edited by Athanase Papadopoulos and Weixu Su.
\newblock {\em Advanced Lectures in Mathematics}, Vol. 49, pgs  129--168.
\newblock International Press, Somerville, Massachusetts, U.S.A., and  Higher Education Press, Beijing, China, 2025.

\bibitem[MM23]{MM1}
M.~Mj and S.~Mukherjee.
\newblock Combining rational maps and {K}leinian groups via orbit equivalence.
\newblock \emph{Proc. Lond. Math. Soc (3)}, 126:1740--1809, 2023.

\bibitem[MM25]{MM2}
M.~Mj and S.~Mukherjee.
\newblock Matings, holomorphic correspondences, and a Bers slice. 
\newblock {\em J. {\'E}c. polytech. Math.}, 12:1445--1502, 2025.

\bibitem[MS13]{MS13} 
M.~Mj and C.~Series.
\newblock Limits of limit sets I.
\newblock {\em Geom. Dedicata}, 167:35--67, 2013.

\bibitem[MS17]{MS17}
M.~Mj and C.~Series.
\newblock Limits of Limit Sets II: Geometrically infinite groups. 
\newblock {\em Geom. Topol.}, 21:647--692, 2017.

\bibitem[OT10]{OT10}
L.~G.~Oversteegen and E.~D.~Tymchatyn.
\newblock Extending isotopies of planar continua.
\newblock {\em Ann. of Math. (2)}, 172(3): 2105--2133, 2010.

\bibitem[Pom75]{Pom75}
C.~Pommerenke.
\newblock \emph{Univalent functions}, with a chapter on quadratic differentials by Gerd Jensen.
\newblock \emph{Studia Math./Math. Lehrb{\"u}cher}, Band XXV,
Vandenhoeck \& Ruprecht, G{\"o}ttingen, 376 pp., 1975.

\bibitem[Pom92]{Pom}
C.~Pommerenke.
\newblock {\em Boundary behaviour of conformal maps}, volume 299 of {\em Grundlehren Math. Wiss.}
\newblock Springer--Verlag, Berlin, 1992, x+300 pp.

\bibitem[Sta04]{Sta04}
M. Stadlbauer.
The return sequence of the Bowen--Series map for punctured surfaces.
\newblock {\em Fund. Math.}, 182:221--240, 2004.
\end{thebibliography}
\end{document}